\documentstyle[epsfig]{article}

\newcommand{\qed}{{$\Box$}}

\newcommand{\ba}{\begin{array}}
\newcommand{\ea}{\end{array}}

\newtheorem{lemma}{Lemma}
\newtheorem{hypo}{Assumption}
\newcommand{\bhypo}{\begin{hypo}}
\newcommand{\ehypo}{\end{hypo}}
\newtheorem{defi}{Definition}
\newcommand{\ble}{\begin{lemma}}
\newcommand{\ele}{\end{lemma}}
\newcommand{\bde}{\begin{defi}}
\newcommand{\ede}{\end{defi}}

\newtheorem{prop}{Proposition}
\newcommand{\epr}{\end{prop}}
\newcommand{\bpr}{\begin{prop}}
\newtheorem{teo}{Theorem}
\newcommand{\bth}{\begin{teo}}
\newcommand{\eth}{\end{teo}}
\newtheorem{rema}{Remark}
\newcommand{\bre}{\begin{rema}}
\newcommand{\ere}{\end{rema}}
\newcommand{\ee}{\end{equation}}
\newcommand{\be}{\begin{equation}}

\begin{document}

\centerline{ \huge \bf Matching Procedure for the}
\vskip 0.3 cm
\centerline{ \huge \bf   Sixth Painlev\'e  Equation (May 2006)}

\vskip 0.3 cm
\centerline{\Large Davide Guzzetti}

\begin{abstract}   

\vskip 0.3 cm
 In the framework of the isomonodromy deformation method, we 
 present a  constructive procedure 
 to obtain  the critical
 behavior of Painlev\'e VI transcendents and solve the connection
 problem. This procedure yields two and  one parameter families of solutions,
 including trigonometric and logarithmic behaviors, and three classes of 
solutions with Taylor expansion
 at a critical point.
\end{abstract}

\vskip 1 cm

\section{Introduction}

This paper appeared in May 2006. I put it on the archive now, with
more that four years of delay, for completeness sake. The paper is
published in  J.Phys.A: Math.Gen. {\bf 39} (2006), 11973-12031, with
some modifications. 
The  sixth Painlev\'e equation is:   
$$
{d^2y \over dx^2}={1\over 2}\left[ 
{1\over y}+{1\over y-1}+{1\over y-x}
\right]
           \left({dy\over dx}\right)^2
-\left[
{1\over x}+{1\over x-1}+{1\over y-x}
\right]{dy \over dx}
$$
$$
+
{y(y-1)(y-x)\over x^2 (x-1)^2}
\left[
\alpha+\beta {x\over y^2} + \gamma {x-1\over (y-1)^2} +\delta
{x(x-1)\over (y-x)^2}
\right]
,~~~~~\hbox{(PVI)}.
$$ 
The generic solution  has essential 
singularities and/or branch points in 0,1,$\infty$. It's behavior at
these points will be called {\it critical}. 
   The other singularities, 
 which depend on the initial conditions, are  poles. 
 A solution of PVI   
 can be analytically continued to a meromorphic function on the universal
   covering of ${\bf P}^1\backslash \{ 0,1 ,\infty \}$. 
 For generic values of the integration constants and of the parameters $\alpha$,$\beta$,$\gamma$,$\delta$,
 it cannot be expressed via elementary or classical
   transcendental functions. For this reason, it  is called a {\it
   Painlev\'e transcendent}.  Solving (PVI) means: ~i) Determine 
the critical behavior of the transcendents 
at the {\it critical points} $x=0,1,\infty$. Such a behavior must depend on 
two integration constants. ~ii) Solve the {\it connection
  problem}, namely: find the relation between couples of integration
constants  at $x=0,1,\infty$.

\vskip 0.2 cm
We use a  matching procedure to study the above two problems. The
 procedure  allows us     
 to compute the first leading terms   
of the critical behavior at a critical point and the associated
 monodromy data.  This 
procedure is essentially the isomonodromy deformation method. The
reason for our terminology is that we make particular use of  the
matching between local solutions of two different reductions of the linear
system of ODE,  associated to (PVI) by the isomonodromy deformation theory.  
This matching allows us to obtain the
leading term(s) of the 
asymptotic behavior of a corresponding  Painlev\'e transcendent
$y(x)$.  In
this sense, we say that our approach is constructive. Namely, we
don't assume any behavior of $y(x)$; rather, we obtain it from the
matching condition. This differs from other authors' approach, who
start by assuming a given asymptotics for $y(x)$ and then 
compute the corresponding monodromy data (and so they 
solve the connection problem). This kind of approach was successfully
used for some of the Painlev\'e equations and allowed many
progresses.  Our approach is developed to tackle with the cases when 
 we don't know - or we are not able to  guess - the
asymptotic behavior. In the case of (PVI), we may say that most of the
solutions are known. But for some points in the space of  monodromy
data, we 
still don't know the corresponding critical behaviors. Our work is
motivated by the need to explore these remaining cases. 
   
 Once the local matching is done, we proceed with a global 
 description of the solutions of the associate linear system of ODE, in order
 to compute its {\it monodromy data}. These are the monodromy data
 associated to the solution $y(x)$, of which the asymptotic  behavior  
 has been obtained by the precedent step. Again, this computation  is 
  done by a (global)  matching, among solutions of the two reduced 
systems and that of the original one.  
This is the main powerful point of the isomonodromy deformation method.  The
 monodromy data are computed in terms of the coefficients of the
 linear system of ODE, which are elementary functions  of the
 parameters (namely, the integration constants) appearing in the leading
 term  of the asymptotic behavior of $y(x)$. The inversion of the
 formulae expressing the monodromy data, gives the leading term of
 $y(x)$ in term of the monodromy data. 

The procedure can be repeated at the other singularities
$x=1,\infty$. In case of (PVI), $x=0,1,\infty$ are equivalent by 
symmetry transformations. These  facts allow to  
solve the connection
 problem (\cite{Jimbo}, 
\cite{DM}, \cite{guz1}, \cite{guz3}, \cite{Boalch}). \footnote{ For reasons of space, we limit ourselves to the
computation of monodromy data, without explaining  how the connection
problem is practically 
 solved  once the monodromy data are computed, and how the analytic
 continuation is done. We refer the reader to  \cite{Jimbo}, 
\cite{DM}, \cite{guz1}, \cite{guz3}, \cite{Boalch}. 
 The  behaviors at 
$x=1$ and $x=\infty$, and the dependence of them on the 
monodromy data are deduced from the behavior at $x=0$ by symmetry 
transformations.
 PVI is invariant for the change of variables  $y(x)=1-\tilde{y}(t)$, $x=1-t$ and simultaneous permutation of $\theta_0$, $\theta_1$. This means that $y(x)$ solves PVI if and only if $\tilde{y}(t)$ solves PVI with permuted  parameters and independent variable $t$.  
Similarly, PVI is invariant for  $y(x)=1/\tilde{y}(t)$, $x=1/t$ and simultaneous permutation of $\theta_{\infty}$, $\theta_0$. It  is invariant for  $y(x)=(\tilde{y}(t)-t)/(1-t)$, $x=t/(t-1)$ and simultaneous permutation of $\theta_{0}$, $\theta_x$.  By composing the third, first and again third symmetries, we get 
$y(x)= \tilde{y}(t)/t$, $t=1/x$ with the permutation of $\theta_1$, $\theta_x$. 
Therefore, the critical points $0,1,\infty$ are equivalent.} 

\vskip 0.2 cm 

The work of Jimbo \cite{Jimbo} 
is the first on the subject.    
For  generic values of 
 $\alpha$, $\beta$, $\gamma$ $\delta$, PVI admits a 
2-parameter class of solutions, with the following critical behaviors:.     
\begin{equation}
y(x)= a x^{1-\sigma}(1+O(|x|^{\epsilon})),~~~~x\to 0,
\label{loc1introduzione}
 \end{equation}
\begin{equation}
y(x)= 1-a^{(1)}(1-x)^{1-\sigma^{(1)}} (1+O(|1-x|^{\epsilon})),~~~~x\to 1,
\label{loc2introduzione}
 \end{equation}
\begin{equation}
y(x)= a^{(\infty)}
 x^{\sigma^{(\infty)}}(1+O(|x|^{-\epsilon})),~~~~x\to \infty,
\label{loc3introduzione}
\end{equation}
where $\epsilon$ is a small positive number, $a^{(i)}$ and
 $\sigma^{(i)}$ are complex numbers such that $a^{(i)}\neq 0$ and 
  $ 0< \Re \sigma<1$, $0< \Re \sigma^{(1)}<1$, $0< \Re \sigma^{(\infty)}<1$.
We remark that  $x$ converges
 to the critical points {\it inside a sector} with vertex
 on the corresponding critical point.   
The  connection problem is to 
finding  the relation among the three pairs $(\sigma,a)$,
$(\sigma^{(1)},a^{(1)})$, $(\sigma^{(\infty)},a^{(\infty)})$.  
In \cite{Jimbo} the problem is solved by the 
isomonodromy deformation method. In particular, the exponents are
determined by the relations:
$$ 
 2 \cos(\pi \sigma)=\hbox{tr}(M_0M_x), ~~~
2 \cos(\pi \sigma^{(1)})=\hbox{tr}(M_1M_x),~~~
2 \cos(\pi \sigma^{(\infty)})=\hbox{tr}(M_0M_1).
$$
Here $M_0$, $M_x$, $M_1$ are monodromy matrices to be introduced
below. 

The above class of solutions  was enlarged in
\cite{Sh} and 
\cite{guz3}, to the values  $\sigma\in{\bf C}$, $\sigma\not \in
(-\infty,0]\cup[1,+\infty)$ (here we consider $x\to 0$).  When $\Re \sigma \geq 1$ or $\Re
  \sigma\leq 0$,   the critical behavior is like the above, but 
 it holds  for
  $x\to 0$ in a spiral-shaped domain in the universal covering of a
  punctured neighborhood of $x=0$, along a  
  paths joining a point $x_0$  to $x=0$. Along special paths which
   approach the movable poles, these solution may have
  behavior $y(x)\sim \sin^{-2}\bigl({i\sigma\over 2}\ln x
  +\varphi(x,a)\bigr)$, 
where
  $\varphi(x,a)$ is a phase depending on the parameter $a$. 
 The transformation 
  $\sigma \mapsto \pm\sigma+2N$, $N\in{\bf Z}$,  
leaves the identity   
tr$(M_0M_x)=
  2\cos(\pi \sigma)$ invariant. Its effect on the solutions is
  studied in \cite{guz3}. 
 As a result, one can reduce to the values $0\leq
  \Re\sigma\leq 1$, $\sigma\neq 0,1$. 
 The reader may find a synthetic
description  of
these results in the review paper \cite{guz4}.

 It is an open problem to determine the critical behavior, say
 at $x=0$, for $\sigma=0,1$. To be more precise, the problem is
 encountered when tr$(M_iM_j)=\pm 2$. These are precisely the points
 of the space of monodromy data mentioned above, in correspondence
 of which we do not know the critical behavior.  In addition, certain 
non-generic values of $\alpha,\beta,\gamma,\delta$ are not yet
 studied.  
The matching procedure is motivated by the need to explore these
 unknown  cases.

\vskip 0.2 cm
\noindent
As a result of  the matching procedure, we obtain: 

{\bf R1)} A two-parameter 
family of solutions, of the type found by Jimbo \cite{Jimbo}. Besides, 
we show that there
are solutions with trigonometric behavior. 

{\bf R2)}   
One-parameter families of  solutions, including a class of logarithmic
solutions.   

\noindent
Together with the results of \cite{Sh} and 
\cite{guz3}, R1) and R2)  will cover all cases tr$(M_iM_j)\neq - 2$, namely
$\sigma\neq 1$. By symmetry transformations, some of the cases  tr$(M_iM_j)=
- 2$ can be obtained from the above results (for example, the {\it
  Chazy solutions} \cite{M}).   

{\bf R3)} The solutions which admit a Taylor expansion at
$x=0$. 

{\bf R4)} We compute the corresponding monodromy
data.   

\noindent
In virtue of the symmetries of (PVI) 
(birational transformations of $(x,y(x))$), it can be shown  that the
solutions with Taylor expansion at $x=0$,  obtained by the matching
procedure,  are the representatives of three equivalent classes, which
include all the solutions admitting a Taylor expansion at a critical
point. If we define $\sigma$
through the relation tr$(M_0M_x)=
  2\cos(\pi \sigma)$, 
 the representatives of three equivalent classes correspond to 
values $\sigma=0$, $\sigma=
  \pm(\theta_1\pm/\mp\theta_\infty)$ and $\sigma=1$.

\vskip 0.2 cm
A further step in the study of PVI, is the problem of the systematic 
classification of all the solutions of (PVI) in terms of the monodromy
data of the associated linear system. As we discussed above,  the  matching 
procedure is 
effective to produce new solutions, associated to monodromy data 
for which the connection problem has not yet been   
 studied. Therefore,  it is a tool to study the 
classification problem. This classification will be done  in another
paper.

\vskip 0.2 cm 

  A matching procedure, to obtain asymptotic behaviors and monodromy data   
in the framework of the isomonodromy deformation method, 
was suggested by Its and Novokshenov in \cite{its}, 
 for the second and third Painlev\'e equations. 
 The work by Jimbo \cite{Jimbo} can be regarded as an 
implicit matching procedure.   
 This method was further developed and
used by Kapaev, Kitaev,  Andreev, and Vardanyan. Here we cite the case of 
 the fifth
Painlev\'e equation,  in \cite{KitaevAndreev}.
   An analogous matching scheme 
is used  in \cite{Kitaev}, for a different problem (limit PVI 
$\to$ PV).

\vskip 0.3 cm
\noindent
{\bf Acknowledgements (May 2006):} The author wishes to thank  Alexander Kitaev for 
introducing him to the matching procedure and for many
discussions.    This paper was written  in RIMS, Kyoto
University,   
 supported by the Kyoto Mathematics COE fellowship.

%%%%%%%%%%%%%%%%%%%%%%%%%%%%%%%%%%%%%%%%%%%%%%%%%%%%%%%%%%%%%%%%%%
%%%%%%%%%%%%%%%%%%%%%%%%%%%%%%%%%%%%%%%%%%%%%%%%%%%%%%%%%%%%%%%%%%
%%%%%%%%%%%%%%%%%%%%%%%%%%%%%%%%%%%%%%%%%%%%%%%%%%%%%%%%%%%%%%%%%%

\vskip 1 cm

\centerline{\large \bf PART I: Matching Procedure   and Results}
\section{Matching Procedure}

 PVI is the isomonodromy deformation equation of  a Fuchsian system of
 differential equations \cite{JMU}: 
\be
   {d\Psi\over d\lambda}=A(\lambda,x,\theta)~\Psi,~~~~~
A(\lambda,x,\theta):=\left[ {A_0(x,\theta)\over \lambda}+{A_x(x,\theta) \over \lambda-x}+{A_1(x,\theta)\over
\lambda-1}\right],~~~\lambda\in{\bf C}.
\label{SYSTEM}
\ee
The  $2\times 2$ matrices  $A_i(x,\theta)$  depend 
 on $x$ in such a way that the monodromy of a fundamental solution $\Psi(\lambda,x)$ 
 does not change for  small deformations of $x$. They also depend on the 
 parameters $\alpha,\beta,\gamma,\delta$ of PVI through more elementary 
parameters $\theta=(\theta_0,\theta_x,\theta_1,\theta_{\infty})$ according to 
the following relations:  
$$
 -A_\infty:=A_0+A_1+A_x = -{\theta_{\infty}\over 2}
 \sigma_3,~~\theta_\infty\neq 0.~~~~~
\hbox{ Eigenvalues}~( A_i) =\pm {1\over 2} \theta_i, ~~~i=0,1,x;
$$
 \begin{equation}
    \alpha= {1\over 2} (\theta_{\infty} -1)^2,
~~~-\beta={1\over 2} \theta_0^2, 
~~~ \gamma={1\over 2} \theta_1^2,
~~~ \left({1\over 2} -\delta \right)={1\over 2} \theta_x^2 
\label{caffe1}
\end{equation}
Here $\sigma_3$ is the Pauli matrix. 
 The equations of monodromy-preserving deformation (Schlesinger equations), can be written in Hamiltonian form and reduce
 to PVI, being the transcendent $y(x)$ solution of
 $A(y(x),x,\theta)_{1,2}=0$. Namely:
\be
y(x)= 
{x~(A_0)_{12} \over x~\left[
(A_0)_{12}+(A_1)_{12}
\right]- (A_1)_{12}},
\label{leadingtermaprile}
\ee
The matrices $A_i(x,\theta)$, $i=0,x,1$, 
 depend on $y(x)$, ${d y(x)\over dx}$ and $\int y(x)$ 
through rational functions, which are given in \cite{JMU}. 
In short, we will  
write $
  A_i=A_i(x)
$.

The product of the monodromy matrices $M_0$, $M_x$, $M_1$ of a fundamental matrix
solution $\Psi$ at $\lambda=0,x,1$ respectively, is equal to the
monodromy at $\lambda=\infty$. The order of the producs depends on the
choice of a basis of loops. As a consequence, the following relation
must hold:  
$$
 \cos(\pi \theta_0) \hbox{tr}(M_1 M_x) + \cos(\pi \theta_1) \hbox{tr}(M_0 M_x) + \cos(\pi \theta_x) \hbox{tr}(M_1 M_0)
$$
$$
= 2\cos(\pi \theta_\infty)+ 4 \cos(\pi \theta_1)\cos(\pi \theta_0)\cos(\pi \theta_x). 
$$

\subsection{Leading Terms of $y(x)$ as a result of Matching}
\label{matchleadingaprile}

 We present the constructive procedure to obtain the leading terms of
a solution $y(x)$, when $x\to 0$. This procedure has been  used for the fifth
Painlev\'e equation by F.V. Andreev and A.V. Kitaev in
\cite{KitaevAndreev}. An analogous scheme 
is used  in \cite{Kitaev}, for a different problem.  In
particular, in \cite{Kitaev}  the non
fuchsian singularity in the $\Psi_{OUT}$-equation (to be introduced
below) appears.  

Since we are considering $x\to 0$, we divide the $\lambda$-plane into
two domains. The ``outside'' domain  is defined for $\lambda$ sufficiently big:
\be
|\lambda|\geq |x|^{\delta_{OUT}},~~~~~\delta_{OUT}>0. 
\label{dominioOUTbasta}
\ee
Therefore, (\ref{SYSTEM}) can be written as: 
\be
{d\Psi\over d \lambda}=
\left[
{A_0+A_x\over \lambda}+{A_x\over \lambda}~\sum_{n=1}^{\infty}\left({x\over \lambda}\right)^n+ {A_1\over \lambda-1}
\right]~\Psi.
\label{SYSTEM1aprile}
\ee
The ``inside'' domain  is defined for $\lambda$ comparable with $x$, namely:
\be
|\lambda|\leq |x|^{\delta_{IN}},~~~~~\delta_{IN}>0. 
\label{dominioINbasta}
\ee
Therefore, $\lambda\to 0$ as $x\to 0$, and we rewrite ({\ref{SYSTEM})
  as: 
\be
{d \Psi \over d\lambda}
=
\left[
{A_0\over \lambda} + {A_x \over \lambda -x} - A_1 \sum_{n=0}^\infty 
\lambda^n
\right]~
\Psi.
\label{SYSTEM0aprile}
\ee

\vskip 0.2 cm 
If the behavior of $A_0(x)$, $A_1(x)$ and $A_x(x)$ is sufficiently
good, we expect that the higher order terms in the series of
(\ref{SYSTEM1aprile}) and (\ref{SYSTEM0aprile}) are small
corrections which can be neglected when $x\to 0$. If this is the case,
(\ref{SYSTEM1aprile}) and (\ref{SYSTEM0aprile}) reduce respectively  to:
\be
{d\Psi_{OUT}\over d \lambda}=
\left[
{A_0+A_x\over \lambda}+{A_x\over \lambda}~\sum_{n=1}^{N_{OUT}}\left({x\over \lambda}\right)^n+ {A_1\over \lambda-1}
\right]~\Psi_{OUT},
\label{nonfuchsianSYSTEMOUT}
\ee
\be
{d \Psi_{IN} \over d\lambda}
=
\left[
{A_0\over \lambda} + {A_x \over \lambda -x} - A_1 \sum_{n=0}^{N_{IN}} 
\lambda^n
\right]~
\Psi_{IN},
\label{nonfuchsianSYSTEMIN}
\ee
where $N_{IN}$, $N_{OUT}$ are suitable integers. 
The simplest  reduction  is to Fuchsian systems:
\be
\label{fuchsianSYSTEMOUT}
{d\Psi_{OUT}\over d \lambda}=
\left[
{A_0+A_x\over \lambda}+ {A_1\over \lambda-1}
\right]~\Psi_{OUT},
\ee
\be
{d \Psi_{IN} \over d\lambda}
=
\left[
{A_0\over \lambda} + {A_x \over \lambda -x}
\right]~
\Psi_{IN}.
\label{fuchsianSYSTEMIN}
\ee

It is a new feature of this paper that we will use reduced
non-fuchsian systems. In the literature, the fuchsian reduction has
been privileged, but we show that in some relevant cases it cannot be
used, being the non-fuchsian reduction necessary. 

\vskip 0.2 cm

 Generally speaking, we can parameterize the elements of
 $A_0+A_x$ and $A_1$ of (\ref{fuchsianSYSTEMOUT}) in terms of
 $\theta_1$, the eigenvalues of $A_0+A_x$ and  the eigenvalues
 $\theta_\infty$ of $A_0+A_x+A_1$. We also need an additional unknown
 function of $x$.   In the same way, we 
can explicitly parameterize the elements of
 $A_0$ and $A_x$  in (\ref{fuchsianSYSTEMIN}) in terms of
 $\theta_0$, $\theta_x$, the   eigenvalues of $A_0+A_x$  and  another
 additional unknown
 function of $x$.  When the reductions   
 (\ref{nonfuchsianSYSTEMOUT}) and (\ref{nonfuchsianSYSTEMIN}) are
 non-fuchsian,  
    particular 
 care must be payed. This will be explained case by case in the paper.
 Our purpose is to find the leading term  of the unknown functions 
 when $x\to 0$, in order to determine  the  critical behavior of 
 $A_0(x)$, $A_1(x)$, $A_x(x)$ and  
 (\ref{leadingtermaprile}).

 The leading term can be  obtained as a result of two facts: 

\noindent
 i) Systems (\ref{nonfuchsianSYSTEMOUT})
 and (\ref{nonfuchsianSYSTEMIN}) are isomonodromic. This imposes 
  constraints on the form of the unknown functions. Typically, one of
 them must be constant. 

\noindent
ii) Two fundamental matrix 
 solutions $\Psi_{OUT}(\lambda,x)$, $\Psi_{IN}(\lambda,x)$ must
match in the region of overlap, provided this is not empty: 
\be\Psi_{OUT}(\lambda,x)
\sim 
\Psi_{IN}(\lambda,x), ~~~~~
|x|^{\delta_{OUT}} \leq |\lambda|\leq |x|^{\delta_{IN}},~~~x\to 0  
\label{overlapINOUT}
\ee  
This relation is to be intended in the sense that the leading terms 
of the local behavior of $\Psi_{OUT}$ and $\Psi_{IN}$ for $x\to
0$ must be equal. 
This  determines a simple relation between the two functions of $x$
appearing in $A_0$, $A_x$, $A_1$, $A_0+A_x$.    (\ref{overlapINOUT})
also implies that $\delta_{IN}\leq\delta_{OUT}$.

\vskip 0.2 cm 
Practically, to fulfill point ii), we will match  a fundamental solution of 
(\ref{nonfuchsianSYSTEMOUT}) for $\lambda\to 0$,  with 
 a fundamental solution of (\ref{nonfuchsianSYSTEMIN}) when $\mu:=\lambda/x\to
\infty$, namely with a solution of: 
\be
{d \Psi_{IN} \over d\mu}
=
\left[
{A_0\over \mu} + {A_x \over \mu -1} - xA_1~ \sum_{n=0}^{N_{IN}} 
x^n\mu^n
\right]~
\Psi_{IN},~~~~~\mu:={\lambda\over x}.
\label{nonfuchsianSYSTEMINmu}
\ee

\vskip 0.2 cm

To summarize,  {\it matching} two fundamental 
 solutions of the  reduced   {\it 
isomonodromic} systems (\ref{nonfuchsianSYSTEMOUT}) and
(\ref{nonfuchsianSYSTEMIN}), we  obtain the leading term(s), for $x\to 0$,  
 of the entries of the matrices of the original system
 (\ref{SYSTEM}). The procedure is algorithmic, no {\it a priori}
 assumption about the behavior being necessary. 

\vskip 0.3 cm
 This method is sometimes called {\it coalescence of singularities},
 because the singularity $\lambda=0$ and $\lambda=x$ coalesce to
 produce  system (\ref{nonfuchsianSYSTEMOUT}), while the singularity
 $\mu={1\over x}$ and $\mu=\infty$ coalesce to produce system
 (\ref{nonfuchsianSYSTEMINmu}). 
Coalescence of singularities was first used by M. Jimbo in 
\cite{Jimbo} to compute the monodromy matrices of (\ref{SYSTEM}) for a
class of solutions of (PVI) with leading term $y(x)\sim ~a~x^{1-\sigma}$,
$0<\Re\sigma<1$.

\subsection{Computation of the Monodromy Data}
\label{MonodromyPasqua}

Let $\Psi$ be a
fundamental matrix solution of (\ref{SYSTEM}), and let $M_0$, $M_x$,
$M_1$, $M_\infty$  be its monodromy matrices at $\lambda=0,x,1,\infty$ 
respectively ($M_\infty$ is the product of $M_0$, $M_x$, $M_1$, the
order depending on the choice of a basis of loops).     
As a consequence of isomonodromicity, 
there exists a fundamental solution $\Psi_{OUT}$ of
(\ref{nonfuchsianSYSTEMOUT}) such that
 $$
M^{OUT}_1=M_1,~~~~~M^{OUT}_\infty=M_\infty,
$$
where $M^{OUT}_1$ and $M^{OUT}_\infty$ are the monodromy matrices of  
$\Psi_{OUT}$ at
$\lambda=1,\infty$. Moreover,  $M^{OUT}_0= M_0M_x$ or $M_xM_0$,
depending on the order of loops. A detailed proof of these facts 
can be found in 
\cite{guz1}. 
 There also exists a fundamental solution $\Psi_{IN}$ of
(\ref{nonfuchsianSYSTEMIN}) such that:
$$
M^{IN}_0=M_0,~~~~~M^{IN}_x=M_x,
$$
where $M^{IN}_0 $ and $M^{IN}_x $  are the monodromy matrices of  
$\Psi_{IN}$ at $\lambda=0,x$. 

\vskip 0.2 cm 
The method of coalescence of singularities is useful when 
the monodromy of the reduced systems
(\ref{nonfuchsianSYSTEMOUT}), (\ref{nonfuchsianSYSTEMIN}) can be
explicitly computed. This is the case  when the reduction
is fuchsian (namely (\ref{fuchsianSYSTEMOUT}),
(\ref{fuchsianSYSTEMIN})), because fuchsian systems with three
singular points are equivalent  to a Gauss hyper-geometric equation
(see Appendix 1). For the non-fuchsian reduction, 
in  general  we can compute the monodromy 
when (\ref{nonfuchsianSYSTEMOUT}),
(\ref{nonfuchsianSYSTEMIN}) are solvable in terms of special or
elementary functions. This will be discussed case by case in the
paper. 

\vskip 0.2 cm 

In order for this procedure to work, not only $\Psi_{OUT}$ and $\Psi_{IN}$ 
 must match with each other, as in subsection
 \ref{matchleadingaprile}, 
but also $\Psi_{OUT}$ 
 must match with a fundamental matrix solution  $\Psi$ of
 (\ref{SYSTEM}) in a domain of the $\lambda$ plane, and   $\Psi_{IN}$ 
 must match with {\it the same}
 $\Psi$ in another domain of the $\lambda$ plane. 
  
\vskip 0.2 cm 
The standard choice of $\Psi$ is as 
follows: 
\be
\Psi(\lambda)= 
\left\{ 
\matrix{
\left[
I+O\left({1\over \lambda}\right)
\right]~\lambda^{-{\theta_\infty\over 2}\sigma_3} \lambda^{R_\infty},&~~~\lambda\to\infty;
\cr
\cr
\psi_0(x) \bigl[I+O(\lambda)\bigr]~\lambda^{{\theta_0\over
    2}\sigma_3}\lambda^{R_0}C_0,&~~~\lambda\to 0;
\cr
\cr
\psi_x(x)\bigl[I+O(\lambda-x)\bigr]~(\lambda-x)^{{\theta_x\over
    2}\sigma_3}(\lambda-x)^{R_x}C_x,&~~~\lambda\to x;
\cr
\cr
\psi_1(x)\bigl[I+O(\lambda-1)\bigr]~(\lambda-1)^{{\theta_1\over
    2}\sigma_3}(\lambda-1)^{R_1}C_1,&~~~\lambda\to 1;
}\right.
\label{PSIlocale}
\ee
Here $\psi_0(x)$,  $\psi_x(x)$, $\psi_1(x)$ are the diagonalizing
matrices of $A_0(x)$, $A_1(x)$, $A_x(x)$ respectively. They are
defined by multiplication to the right by arbitrary diagonal matrices,
possibly depending on $x$. $C_\kappa$, $\kappa=\infty,0,x,1$,
 are invertible {\it connection matrices}, independent of $x$
\cite{JMU}. Each $R_\kappa$,  
$\kappa=\infty,0,x,1$, is also independent of $x$, and:
$$
R_\kappa=0 \hbox{ if } \theta_\kappa\not\in {\bf Z},~~~~~
R_\kappa=\left\{
\matrix{
\pmatrix{0 & *\cr 0 & 0},~~~ \hbox{ if } \theta_\kappa>0 \hbox{
  integer}
 \cr
\cr
\pmatrix{0 & 0\cr * & 0},~~~ \hbox{ if } \theta_\kappa<0 \hbox{
  integer}
}
\right.
$$
If $\theta_i=0$, $i=0,x,1$, then 
 $R_i$ is to be considered the Jordan form $\pmatrix{0 & 1 \cr 0 & 0}$
of $A_i$. If $\theta_\infty=0$, $R_\infty=0$. 
Note that for the loop $\lambda \mapsto \lambda e^{2\pi i}$,
$|\lambda|>\max\{1,|x|\}$, we immediately compute the monodromy at infinity:  
$$
M_\infty=\exp\{-i\pi\theta_\infty\}~\exp\{ 2\pi i R_\infty\}. 
$$ 

\vskip 0.5 cm
Let  
$\Psi_{OUT}$ and $\Psi_{IN}$ be the solutions of 
 (\ref{nonfuchsianSYSTEMOUT}) and (\ref{nonfuchsianSYSTEMIN}) 
matching as in (\ref{overlapINOUT}). We explain how they are matched
with (\ref{PSIlocale}).

\vskip 0.3 cm
\noindent
{\bf (*) Matching $\Psi~\leftrightarrow~ \Psi_{OUT}$:}~~

$\lambda=\infty$ is a fuchsian singularity of
(\ref{nonfuchsianSYSTEMOUT}), with residue
  $-A_\infty/\lambda$. Therefore, we can always find a
  fundamental matrix solution with
  behavior:   
$$
\Psi_{OUT}^{Match}~=\left[
I+O\left({1\over \lambda}\right)
\right]~\lambda^{-{\theta_\infty\over 2}\sigma_3} \lambda^{R_\infty},
~~~\lambda\to\infty.
$$
This solution matches with $\Psi$. 
Also $\lambda=1$ is a fuchsian singularity of
(\ref{nonfuchsianSYSTEMOUT}). Therefore, we have: 
$$
 \Psi_{OUT}^{Match}~
=\psi_1^{OUT}(x)\bigl[I+O(\lambda-1)\bigr]~(\lambda-1)^{{\theta_1\over
    2}\sigma_3}(\lambda-1)^{R_1}C^{OUT}_1,~~~\lambda\to 1;
$$
Here $C^{OUT}_1$ is a suitable connection matrix. $\psi_1^{OUT}(x)$ is
the matrix that diagonalizes the leading 
terms of $A_1(x)$. 
Therefore, $\psi_1(x)\sim \psi_1^{OUT}(x)$ for
$x\to 0$. As a consequence of  isomonodromicity, $R_1$ is the same of 
$\Psi$.  

As a consequence of the matching  $\Psi~\leftrightarrow~
\Psi_{OUT}^{Match}$, the monodromy of $\Psi$ at $\lambda=1$ is:  
$$
M_1={C_1}^{-1}\exp\{i\pi \theta_1\sigma_3\} \exp\{2\pi i R_1\} C_1,
~~\hbox{ with } 
C_1\equiv C^{OUT}_1.
$$
\vskip 0.2 cm 
We finally need an invertible connection matrix  $C_{OUT}$ to connect  
$\Psi_{OUT}^{Match}$ with the solution  $\Psi_{OUT}$ appearing in 
 (\ref{overlapINOUT}). Namely,   $
\Psi_{OUT}^{Match}= \Psi_{OUT} C_{OUT}.
$

\vskip 0.3 cm
\noindent
{\bf (*) Matching $\Psi~\leftrightarrow~ \Psi_{IN}$:}~~ 

\vskip 0.2 cm
As a consequence of the matching $\Psi ~\leftrightarrow~
\Psi_{OUT}^{Match}$, we have to choose  the  IN-solution which
 matches with $\Psi_{OUT}^{Match}$. This is $
\Psi_{IN}^{Match}:=\Psi_{IN} C_{OUT}$.  
 
\vskip 0.18 cm
Now,  $\lambda =0, x$ are fuchsian singularities of
 (\ref{nonfuchsianSYSTEMIN}). Therefore: 
$$
\Psi_{IN}^{Match}=
\left\{
\matrix{
\psi_0^{IN}(x) \bigl[I+O(\lambda)\bigr]~\lambda^{{\theta_0\over
    2}\sigma_3}\lambda^{R_0}C^{IN}_0,&~~~\lambda\to 0;
\cr
\cr
\psi_x^{IN}(x)\bigl[I+O(\lambda-x)\bigr]~(\lambda-x)^{{\theta_x\over
    2}\sigma_3}(\lambda-x)^{R_x}C^{IN}_x,&~~~\lambda\to x;
}
\right.
$$
The above hold for fixed small $x\neq 0$.  
Here $C^{IN}_0$ and  $C^{IN}_x$ are suitable connection matrices. 
$\psi_0^{IN}(x)$ and $ \psi_x(x)^{IN}$ are diagonalizing matrices of 
the leading terms of  
$A_0(x)$ and $A_x(x)$. For $x\to 0$ they match with  $\psi_0(x)$ and 
$ \psi_x(x)$ of $\Psi$ in (\ref{PSIlocale}).  
On the other hand, as a consequence of isomonodromicity, the matrices
$R_0$ and $R_x$ are the same of
$\Psi$. 

By virtue of the matching   $\Psi~\leftrightarrow~
\Psi_{IN}^{Match}$,  
the connection matrices  $C_0$ and $C_x$ coincide with 
  the $x$-independent connection matrices  $C^{IN}_0$,
$C^{IN}_x$  respectively. As a result, we obtain
the monodromy matrices for $\Psi$:
$$
M_0= {C_0}^{-1} \exp\{i\pi\theta_0\sigma_3\} \exp\{ 2\pi i R_0\} C_0,
~~~~~ C_0\equiv C_0^{IN},
$$ 
$$
 M_x= {C_x}^{-1} \exp\{i\pi\theta_x\sigma_3\} \exp\{ 2\pi i R_x\}
C_x,~~~~~ C_x\equiv C_x^{IN}. 
$$

\vskip 0.2 cm 
Our reduction is useful if 
 the connection matrices $C_1^{OUT}$, $C^{IN}_0$, $C^{IN}_x$ can be computed
explicitly. This is  possible  for the fuchsian reduced systems   
(\ref{fuchsianSYSTEMOUT}), (\ref{fuchsianSYSTEMIN}). 
For  non-fuchsian reduced systems, we will discuss the 
computability case by case.

\section{Results}

 In the following, it is understood that $x\to 0$ inside a sector. 
 Namely, $\arg(x)$ is bounded. 

\subsection{Critical Behaviors: Result {\bf R3)}.}

  The  novelty  of this paper is that the matching
 procedure is applied to  non-fuchsian systems 
 (\ref{nonfuchsianSYSTEMOUT}) and
 (\ref{nonfuchsianSYSTEMIN}). As a result, 
we obtain all the solutions that admit a Taylor expansion 
$$
y(x)=b_0+b_1x+b_2x^2+...=\sum_{n=0}^{\infty} b_n x^n,~~~~~x\to 0.
$$
 Precisely, we obtain the representative
 solutions of three equivalence classes, the equivalence relation
 being the birational transformations \cite{Okamoto} 
of  Appendix 3 and formula 
(\ref{nuovasimmetria}). Our result is the following.

\bth
\label{TAYLORpasqua}
 The solutions of (PVI) with Taylor expansion at $x=0$ 
 are divided into four equivalent
 classes (one being that of singular solutions $y=0,1,x$). 
The representatives can be chosen as follows:  

\vskip 0.3 cm
\noindent
1) Singular solution $y=1$.

\vskip 0.3 cm
\noindent 
2) $\theta_\infty\neq 1$, $\theta_1-\theta_\infty\not\in{\bf Z}$ ~{\rm 
 [representative of $\theta_1\pm\theta_\infty\not\in{\bf Z}$]}:
\be
y(x)= {\theta_1-\theta_\infty +1\over 1-\theta_\infty} 
~+
{\theta_1 [(\theta_1-\theta_\infty)(\theta_1-\theta_\infty+2)
+\theta_x^2
-\theta_0^2]
\over 
2 (\theta_\infty-1)(\theta_\infty-\theta_1)(\theta_\infty-\theta_1-2)}
~x
~+
\sum_{n=3}^\infty b_n(\theta_1,\theta_\infty,\theta_0,\theta_x)~x^n.
\label{form1}
\ee
The coefficients are certain rational 
functions of $\theta_0,\theta_\infty,\theta_0,\theta_x$.

\vskip 0.3 cm
\noindent
3) $\theta_1=\theta_\infty\neq 1$, $\theta_0=\pm\theta_x$
~{\rm [representative of $\theta_1\pm\theta_\infty\in{\bf Z}$,
  $\theta_x\pm\theta_0\in{\bf Z}$]}: 
\be
y(x)= {1\over 1-\theta_\infty}+ax ~+\sum_{n=2}^\infty b_n(a;\theta_0,\theta_\infty)x^n.
\label{form2}
\ee
 The coefficients are certain rational 
functions of $\theta_0$, $\theta_\infty$  and a parameter $a\in{\bf C}$.

\vskip 0.3 cm
\noindent
4) $\theta_\infty=1$, $\theta_1=0$~{\rm [representative of
  $\theta_1\pm\theta_\infty\in{\bf Z}$, $\theta_\infty\in{\bf
  Z}\backslash \{0\}$]}: 
\be
y(x)=a+{1-a\over 2} (1+\theta_0^2-\theta_x^2)~x ~+\sum_{n=2}^\infty
b_n(a;\theta_0;\theta_x)x^n. 
\label{form3}
\ee
 The coefficients are certain rational 
functions of $\theta_0$, $\theta_x$ and a parameter $a\in{\bf C}$. 
\eth

 The monodromy data 
associated to the above solutions is given in theorem \ref{thMONODROMY}. 

\vskip 0.2 cm 
 The symmetry $\theta_1\mapsto
 -\theta_1$, which leaves (PVI) invariant, transforms   (\ref{form1})
 into: 
\be
y(x)= {\theta_1+\theta_\infty -1\over \theta_\infty-1} 
~+
{\theta_1 [(\theta_1+\theta_\infty)(\theta_1+\theta_\infty-2)+\theta_x^2
-\theta_0^2]
\over 
2 (1-\theta_\infty)(\theta_\infty+\theta_1)(\theta_\infty+\theta_1-2)}
~x~
+\sum_{n=3}^\infty b_n(-\theta_1,\theta_\infty,\theta_0,\theta_x)~x^n.
\label{riuffa}
\ee
Here $\theta_\infty\neq 1$, $\theta_1+\theta_\infty\not\in{\bf Z}$. The
coefficients $b_n$ are the same of    (\ref{form1}).

\vskip 0.2 cm
The convergence of the Taylor series can be proved by a Briot-Bouquet
like argument. This will not be done here, for reasons of space. The reader can
find the general procedure in \cite{INCE} and an application to the
fifth Painlev\'e equation in \cite{kaneoya}

\vskip 0.3 cm
\noindent
{\bf Comments:}
\vskip 0.2 cm
\noindent 
 {\bf 1) Characterization of solutions  $
y(x)=\sum_{n=0}^\infty b_nx^n$, $b_0\neq 0$}.

\vskip 0.2 cm
{\bf (a)} There always exists one solution  
(\ref{form1}) when $\theta_1-\theta_\infty \not \in {\bf Z}$;
 there always exists one solution (\ref{riuffa}) when
 $\theta_1+\theta_\infty\not \in {\bf Z}$. The coefficients $b_n$
 depend rationally on $\theta_\kappa$, $\kappa=0,x,1,\infty$. 
 {\bf (b)} There is a one-parameter family of solutions equivalent to 
(\ref{form2}), when  $\theta_1\pm\theta_\infty\in {\bf Z}$   and  
$\theta_0\pm\theta_x$ has a particular integer value.  
The coefficients $b_n$ depend rationally on a complex parameter $a$
and $\theta_\infty, \theta_0$. {\bf (c)}  Finally, there is a
one-parameter 
family of solutions 
 equivalent to (\ref{form3}), when $\theta_1\pm\theta_\infty\in{\bf
   Z}$, and 
  $\theta_\infty$ has a  particular integer value; 
the coefficients $b_n$ depend rationally on a complex parameter $a$ and 
 $\theta_0,\theta_x$. 
The singular solutions $y=0,1,x$ are possibly obtained by birational
transformations of (\ref{form1}), (\ref{form2}), (\ref{form3}).

 The
coefficients $b_n$ can  always be computed recursively by direct 
substitution into (PVI).
 We will clarify these facts by some examples in Appendix 4.

\vskip 0.2 cm
\noindent
 {\bf 2) Characterization of solutions  $
y(x)=\sum_{n=1}^\infty b_n x^n$, $b_1\neq 0$}.

\vskip 0.2 cm
These solutions are obtained from those of theorem \ref{TAYLORpasqua} by the 
symmetry. 
\be
\theta_x\mapsto \theta_1,~~~\theta_0\mapsto
\theta_\infty-1,~~~\theta_1\mapsto \theta_x,~~~\theta_\infty\mapsto \theta_0+1;~~~~~
y(x)\mapsto {x\over y(x)}. 
\label{nuovasimmetria}
\ee
The solutions obtained from the singular solution $y=1$ and
(\ref{form1}), (\ref{form2}), (\ref{form3}) are respectively: 
\vskip 0.2 cm
\noindent
1) Singular solution $y(x)=x$. 

\vskip 0.2 cm
\noindent
2) $\theta_0\neq 0$, $\theta_0\pm\theta_x\not\in {\bf Z}$: 
\be 
y(x)= {\theta_0\over \theta_0\pm\theta_x}x
~\pm 
~{\theta_0\theta_x ~\left[
(\theta_0\pm\theta_x)^2+\theta_1^2-\theta_{\infty}^2+2\theta_\infty-2
 \right]\over 2(\theta_0\pm\theta_x)^2~\left[
(\theta_0\pm\theta_x)^2-1
\right]
}x^2
~+\sum_{n=3}^\infty b_n(\theta_0,\theta_x,\theta_1,\theta_\infty)x^n.
\label{taylor1}
\ee

\vskip 0.2 cm
\noindent
3) $\theta_0+\theta_x=1$, $\theta_0\neq 0$,
   $\theta_1=\pm(\theta_\infty-1)$:
\be
y(x)=\theta_0x~+a~x^2~+ \sum_{n=3}^\infty b_n(a;\theta_0,\theta_\infty)x^n.
\label{taylor2}
\ee
%For example, $
% b_3(a;\theta_0,\theta_\infty)={(\theta_0-1)[\theta_0(2\theta_0-1)
%\theta_\infty(\theta_\infty-2)-6a]\over  
%  6}$.

\vskip 0.2 cm
\noindent
4) $\theta_x=\theta_0=0$.
\be
y(x)=ax~+{a(a-1)\over 2} (\theta_1^2-(\theta_\infty-1)^2-1)x^2~+
\sum_{n=3}^\infty b_n(a;\theta_1,\theta_\infty)x^n.
\label{taylor3}
\ee

\vskip 0.2 cm 
{\bf (a)} (PVI) has always  
one or both  solutions  (\ref{taylor1}) when $\theta_0\pm
\theta_x\not \in{\bf Z}$. Also when 
 $\theta_0+\theta_x$ (or $\theta_0-\theta_x$) is
integer, (PVI) has  a solution (\ref{taylor1}) corresponding to
$\theta_0-\theta_x$ not integer (or $\theta_0+\theta_x$ not
integer). {\bf (b)} When 
 $\theta_0+\theta_x$ or $\theta_0-\theta_x$ is integer, (PVI) has a 
 1-parameter family of solutions equivalent (by birational
 transformations) to (\ref{taylor2}); this family 
 exists provided that $\theta_1\pm\theta_\infty$ has a particular
integer value. {\bf (c)} 
When $\theta_0+\theta_x$ or $\theta_0-\theta_x$ is integer and  
$\theta_0$
 has a particular integer value, there is a one parameter family  of
 solutions equivalent to (\ref{taylor3}).

\vskip 0.2 cm
\noindent
{\bf 3)}~ (PVI) has a 
 one-parameter family of 
 solutions of the type:
\be
y(x)=y_0(x) +
y_1(x)~ax^{\omega}~+y_2(x)~\bigl(ax^{\omega}\bigr)^2~+...~=\sum_{N=0}^\infty
y_N(x)~\bigl(ax^{\omega}\bigr)^{N},~~~x\to 0;
\label{1parameter}
\ee
 where  the parameter is $a\in {\bf C}$, and the $y_N(x)$'s are Taylor
 series:  
 $$
y_N(x)=\sum_{k=0}^\infty
b_{k,N}(\theta_1,\theta_\infty,\theta_0,\theta_x)~x^k,~~~~~~~x\to 0. 
$$  
Either  $y_0(x)$ 
 is (\ref{riuffa}) and
  $\omega=\pm(\theta_1+\theta_\infty-1)$, or $y_0(x)$  is (\ref{form1}) and
  $\omega=\pm(\theta_\infty-\theta_1-1)$ . 
The conditions $|\Re\omega|<1$, $\omega\neq 0$ 
hold.  
The coefficients $b_{k,N}(\theta_1,\theta_\infty,\theta_0,\theta_x)$
are certain rational functions that can be recursively  determined
by direct substitution into (PVI). 
 These solutions are the  immages  of solutions  
(\ref{passero}) and (\ref{aquila}) respectively, through the symmetry
(\ref{nuovasimmetria}). 
Solutions (\ref{passero}) and 
(\ref{aquila}) are a sub-case of theorem \ref{thsigmano1}, 
obtained by the matching procedure. 
 
Taylor solutions (\ref{form1}), (\ref{riuffa}) are a special case of
(\ref{1parameter}), when the parameter is zero.  Solutions
(\ref{form2}) and (\ref{form3}) -- and their images by symmetry -- are  
one parameters families  of type (\ref{1parameter}), in
non generic cases when $\omega\in{\bf Z}$. 

Further study of one-parameter solutions,  
including non-generic cases when  $\theta_\nu$ and/or some sum of
two $\theta_\nu$'s  are integer (including logarithmic one parameter
families), will be
presented in another paper devoted to the general classification
problem.

\vskip 0.3 cm
\noindent
{\bf 4)}~
Solutions (\ref{form1})
 and the equivalent solutions (\ref{riuffa}), (\ref{taylor1})  
were also derived  in
\cite{kaneko} by substitution of a Taylor expansion in (PVI). 
The corresponding monodromy  was computed by coalescence of
singularities of a Heun's type (scalar) equation.

\subsection{ Critical Behaviors:  Results {\bf R1), R2)}.}
\vskip 0.2 cm 

\noindent
 We now consider cases when (\ref{SYSTEM}) can be reduced to the 
fuchsian systems (\ref{fuchsianSYSTEMOUT}) and
 (\ref{fuchsianSYSTEMIN}).  Let $\sigma$ be a complex number defined, up
to sign, by:
$$
\hbox{tr}~ (M_0M_x)=2\cos(\pi\sigma),~~~~~| \Re \sigma|\leq 1. 
$$
Actually, $\pm\sigma/2$ are the eigenvalues of $\lim_{x\to 0}
(A_0+A_x)$. The matching procedure yields the following result.

\bth

Let $r\in {\bf C}$ and  $\sigma$ be as above, with the restriction $
|\Re\sigma| <1$. 
 (PVI)  has a family of solutions depending on the two parameters
 $r$, $\sigma$.  
 The  leading terms of the critical behavior for $x\to 0$ 
may be  parametrized as follows:   

\vskip 0.2 cm
\noindent
For $\sigma\neq 0$: 
\be
y(x)\sim ~
\left\{
\matrix{
{1\over r}~{[\sigma^2-(\theta_0+\theta_x)^2][(\theta_0-\theta_x)^2-\sigma^2]
\over 16 \sigma^3}~x^{1-\sigma},~~~~~~~~~~~~ &~~~\hbox{\rm if $\Re \sigma>0$};
\cr
\cr
-{r\over \sigma}~x^{1+\sigma},~~~~~~~~~~~~~~~~~~~~~~~~~~~~~~~~~~~~
 &~~~\hbox{\rm if $\Re \sigma<0$};
\cr
\cr
x\left\{
i{A}~\sin\bigl(i\sigma\ln
x+\phi\bigr)
+{\theta_0^2-\theta_x^2+\sigma^2\over 2\sigma^2}
\right\}
,&~~~\hbox{\rm if $\Re\sigma=0$}.
}
\right.
\label{asyjimb}
\ee
In the above formulae, $r\neq 0$ and 
$$
\phi:= i 
\ln {2r\over \sigma A},~~~~~{A}:=\left[{\theta_0^2\over
    \sigma^2}-\left({\theta_0^2-\theta_x^2+\sigma^2\over
    2\sigma^2}\right)^2\right]^{1\over 2}
$$
\vskip 0.2 cm
\noindent
For special values of $\sigma\neq 0$: 
\be
y(x)\sim{\theta_0\over \theta_0+\theta_x}~x ~\mp~{r\over\theta_0+\theta_x}
~x^{1+\sigma},~~~\sigma=\pm(\theta_0+\theta_x)\neq 0,
\label{passero}
\ee
\be
y(x)\sim {\theta_0\over \theta_0-\theta_x}~x~\mp
~{r\over\theta_0-\theta_x}~x^{1+\sigma},
~~~
\sigma=\pm(\theta_0-\theta_x)\neq 0.
\label{aquila}
\ee
\vskip 0.2 cm
\noindent
For $\sigma = 0$:
\be
y(x)\sim
\left\{
\matrix{
x\left\{
{\theta_x^2-\theta_0^2\over 4} 
\left[
\ln x +{4r+2\theta_0\over \theta_0^2-\theta_x^2}
\right]^2 +{\theta_0^2\over \theta_0^2-\theta_x^2}
\right\}
,&~~~\theta_0\neq\pm\theta_x,
\cr
\cr
x~(r~\pm~\theta_0~\ln x),&~~~\theta_0=\pm\theta_x.
}
\right.
\label{asyjimb0}
\ee
\label{thsigmano1}
\eth
%%%
%%%
%%Si puo` ipotizzare che:
%$$
%y(x)=\sum_{n\geq m\geq 1} x^n(a_{nm}x^{-\sigma}+b_{nm}+c_{nm}x^\sigma)^m
%$$
%
% Credo che sia sbagliata. Probabilmente la seguente e' corretta:
%$$
%y(x)=x\Bigl[
%x^{-\sigma}\sum_{k=0}^\infty \sum_{m=0}^\infty\sum_{n=-m}^m a_{kmn}
%x^{k+m-n\sigma} + \sum_{k=0}^\infty \sum_{m=0}^\infty\sum_{n=-m}^m b_{kmn}
%x^{k+m-n\sigma} +x^{\sigma}\sum_{k=0}^\infty
%\sum_{m=0}^\infty\sum_{n=-m}^m c_{kmn} 
%x^{k+m-n\sigma}  
%\Bigr]
%$$
\noindent
{\bf Comments:}

\vskip 0.2 cm
\noindent
{\bf 1)} $r$ can be computed as a function of the monodromy data. See
  (\ref{bfa}) and comments there. The sign (branch) of the two quare
  roots appearing 
  in $\phi$ and $A$ is the same. $x\to 0$ in a sector of width
  less then $2\pi$.  

\vskip 0.2 cm
\noindent
{\bf 2) Sub-cases of theorem \ref{thsigmano1}}.

\vskip 0.2 cm
i) When $\sigma\neq 0$,  the result of the theorem includes the  
 sub-cases (\ref{passero}) and (\ref{aquila}). 
%\be
%y(x)={\theta_0\over \theta_0+\theta_x}~x ~\mp~{r\over\theta_0+\theta_x}
%~x^{\sigma+1},~~~\sigma=\pm(\theta_0+\theta_x)\neq 0,
%\label{passero}
%\ee
%\be
%y(x)= {\theta_0\over \theta_0-\theta_x}~x~\mp
%~{r\over\theta_0-\theta_x}~x^{1+\sigma},
%~~~
%\sigma=\pm(\theta_0-\theta_x)\neq 0.
%\label{aquila}
%\ee
If $r=0$, $\theta_0\neq 0$, $\theta_0\pm\theta_x\not\in{\bf Z}$,
direct substitution into (PVI) gives the two Taylor expansions
(\ref{taylor1}).

If $r\neq 0$, (\ref{passero}) and (\ref{aquila}) are a 1-parameter
family, with the restriction $|\Re\sigma|<1$.  
The symmetry (\ref{nuovasimmetria}) transforms them into the solutions 
(\ref{1parameter}), the leading terms being respectively:
$$
y(x)~\sim{\theta_\infty+\theta_1-1\over \theta_\infty-1}~\left(1 \pm
{r\over \theta_\infty-1}~x^{\omega}\right),~~~~~
  \omega=\pm(\theta_\infty+\theta_1-1)\neq 0,
$$
$$
y(x)~\sim{\theta_\infty-\theta_1-1\over \theta_\infty-1}~\left(1 \pm
{r\over \theta_\infty-1}~x^{\omega}\right),~~~~~
  \omega=\pm(\theta_\infty-\theta_1-1)\neq 0, 
$$
with the restriction $|\Re\omega|<1$ .

\vskip 0.2 cm
ii) The case $\sigma=0$ includes the sub-case
$y(x)\sim r x$, which occurs for $\theta_0=\theta_x$, $\theta_0=0$.
 By direct substitution in (PVI) we obtain a series: 
$$
y(x)=r~x~+
\sum_{n=3}^\infty b_n(r,\theta_1,\theta_\infty)x^n,
~~~~~
\theta_0=\theta_x=0,~~~r\neq 0,1.
$$
This is again solution (\ref{taylor3}). 
Note that for $r=0,1$ we have the singular solutions $y=0$, $y=1$. 
Note also that the special sub-sub-case $\theta_0=\theta_x=\theta_1=0$  
has applications  in the theory of  {\it semi-simple Frobenius manifolds} of
dimension three \cite{Dub2}  \cite{guz2}. 

\vskip 0.2 cm
\noindent
{\bf 3)} The first two behaviors  (\ref{asyjimb}) were derived 
 in \cite{Jimbo}.  
 The existence of such solutions  was proved
by assuming that the matrices $A_0$, $A_x$, $A_1$ have a certain
critical behavior for $x\to 0$, and proving that such matrices solve
the Schlesinger equations. Then, the monodromy data were computed
by coalescence of singularity. These solutions where further studied
 in  \cite{DM}, 
\cite{guz1}, \cite{guz3}, \cite{Boalch}. We show
that 
 these solutions can
be obtained  without any assumption by the matching procedure,
together with  the solutions (\ref{asyjimb0}) and the first of 
(\ref{asyjimb}), which
 do not appear  in \cite{Jimbo}.

 \vskip 0.2 cm
\noindent
{\bf 4)} 
The class of solutions (\ref{asyjimb}) was enlarged in
\cite{Sh} and 
\cite{guz3}, to the values  $\sigma\in{\bf C}$, $\sigma\not \in
(-\infty,0]\cup[1,+\infty)$.  When $\Re \sigma \geq 1$ or $\Re
  \sigma\leq 0$,   the critical behavior is like the 
first of  (\ref{asyjimb}), and it holds  for
  $x\to 0$ in a spiral-shaped domain in the universal covering of a
  punctured neighborhood of $x=0$, along a  
  paths joining a point $x_0$  to $x=0$. Along special paths which
   approach the movable poles, these solution may have
  behavior $y(x)\sim \sin^{-2}\bigl({i\sigma\over 2}\ln x
  +\varphi(x,r)\bigr)$, 
where
  $\varphi(x,r)$ is a phase depending on the parameter $r$. 
 The transformation 
  $\sigma \mapsto \pm\sigma+2N$, $N\in{\bf Z}$,  
leaves the identity   
tr$(M_0M_x)=
  2\cos(\pi \sigma)$ invariant. Its effect on the solutions is
  studied in \cite{guz3}. 
 As a result, one can reduce to the values $0\leq
  \Re\sigma\leq 1$, $\sigma\neq 0,1$. 
 We cannot enter into more details here. The reader may find a synthetic
description  of
these results in the review paper \cite{guz4}.

\vskip 0.2 cm 
\noindent
{\bf 5)} Solutions with expansion:   
{\small
$$
y(x)= x (A_1 + B_1 \ln x + C_1 \ln^2 x + D_1 \ln^3 x + ...)+
x^2(A_2+B_2\ln x +...)+...,~~~~~x\to 0.  
$$
}
are all included it theorems \ref{TAYLORpasqua} and
\ref{thsigmano1}. Actually, only the following cases are possible:      
{\small
\be
y(x)=
\left\{
\matrix{ 
 {\theta_0\over \theta_0\pm\theta_x} x + O(x^2) ~~\hbox{ [Taylor
     expansion]},
\cr\cr 
x~\left(
{\theta_0^2-B_1^2\over \theta_0^2-\theta_x^2} + B_1\ln x + 
{\theta_x^2-\theta_0^2\over 4} \ln^2 x 
\right)
+x^2(...)  +...,
\cr\cr
x~(A_1\pm \theta_0 \ln x)+x^2(...)+...,~~~\hbox{ and }\theta_0=\pm \theta_x.
}
\right.
\label{fivecases}
\ee
}
$A_1$ and $B_1$ are parameters. We see that the higher orders in
(\ref{asyjimb0})  are $O(x^2\ln^m x)$, for some integer $m>0$. 

\vskip 0.2 cm 
\noindent
{\bf 6)} The symmetry (\ref{nuovasimmetria}) applied to solutions
(\ref{asyjimb0}) gives:
{\small
$$
y(x)\sim 4\left( 
\bigl(\theta_1^2-(\theta_\infty-1)^2\bigr)
\left[
\ln x -{4r +2(\theta_\infty-1)\over \theta_1^2-(\theta_\infty-1)^2}
\right]^2 -{4(\theta_\infty-1)^2\over \theta_1^2-(\theta_\infty-1)^2}
\right)^{-1},
$$
}
namely:
\be
y(x)={4\over \bigl[\theta_1^2-(\theta_\infty-1)^2\bigr]~ \ln^2 x}
~\left[1+{8r+4(\theta_\infty-1)\over 
      \theta_1^2-(\theta_\infty-1)^2}~{1\over  \ln x}
  +O\left({1\over \ln^2 x}\right)\right],
\label{chazygeneral}
\ee
and 
$$
y(x)= {\pm1\over (\theta_\infty-1)~\ln x} 
\left[ 1\mp{r\over(\theta_\infty-1)~ \ln x} 
 +O\left({1\over \ln^2 x}\right)\right],~~~~~~~ \theta_\infty\mp\theta_1=1.
$$ 
The higer orders $O(1/\ln^2 x)$  include powers $x^n(\ln x)^{\pm m}$. 
The so called {\it Chazy solutions}, studied in \cite{M} 
 for the special case
$\theta_0=\theta_x=\theta_1=0$, $\theta_\infty=-1$, have the behavior
 (\ref{chazygeneral}).

\vskip 0.2 cm 
\noindent
{\bf 7)} When this paper was completed, I received a communication 
by the first author of \cite{Bruno}. In \cite{Bruno}  it is proved that 
 (PVI) has 
solutions with expansion at $x=\infty$, or $x=0$, of the form 
 $y= c_r x^r + \sum_s c_s
x^s$, $c_r\in{\bf C}$. The $c_s$'s
 are either complex constants or polynomials
in $\ln x$. $r$ and $s$  are integer or complex. If
$r$ is complex, the restriction  $\Re r\in(0,1)$ holds. 
The method used in \cite{Bruno} is a power geometry
technique. The connection problem and the characterization of the
associated monodromy data are not studied.

\vskip 0.2 cm 
\noindent
{\bf 8)} When this paper was already accepted for publication, I received a 
private communication about a recent work \cite{qlu},  
 on  the asymptotics of the 
{\it real} solutions of (PVI).  
The asymptotic behaviors obtained in \cite{qlu}  are 
 of the type of our theorem  \ref{thsigmano1}, 
namely (\ref{asyjimb}) and the first behavior in (\ref{asyjimb0}). The tool 
used is a method of successive approximations. 
So, the results are local, and the 
 connection problem is not studied. Moreover, some genericity conditions on 
the coefficients of PVI seem to be necessary 
 (so, for example, the second solution in 
(\ref{asyjimb0}) cannot be obtained).

\subsection{Monodromy: Result {\bf R4)}.} 
In this paper, we computed the monodromy for the Taylor-expanded
solutions, which correspond to non-fuchsian reductions of system
(\ref{SYSTEM}).  
Because of the symmetries of (PVI), we can limit ourselves to 
 the monodromy data for the representative solutions 
(\ref{form1}), (\ref{form2}) and (\ref{form3}). 

\bth
\vskip 0.2 cm
a) Let $\theta_\kappa\not \in {\bf Z}$, $\kappa=0,1,x,\infty$. 
A representation for the monodromy matrices of the solution
(\ref{form1}) is:
$$
M_0=C_{0\infty} ~\exp\{i\pi\theta_0\sigma_3\}~
C_{0\infty}^{-1},
$$
$$
M_x=C_{0\infty}~ C_{01}^{-1}~ 
\exp\{i\pi\theta_x\sigma_3\}
 ~ C_{01}~
C_{0\infty}^{-1}.
$$
$$
M_1=\exp\{-i\pi\theta_1\sigma_3\},
~~~M_\infty=\exp\{-i\pi\theta_\infty \sigma_3\}.
$$
The matrices $C_{0\infty}$ and $C_{01}$ are:
\vskip 0.3 cm 
\be
C_{0\infty}:=
\left[
\matrix{
{ 
\Gamma\left(1+{\theta_1\over 2}-{\theta_\infty\over 2}\right)\Gamma(1+\theta_0)
e^{ i{\pi\over 2} \left[
    \theta_0+\theta_x+\theta_\infty-\theta_1\right]}
\over
\Gamma\left({\theta_0\over 2}+{\theta_x\over 2} +{\theta_1\over
  2}-{\theta_\infty\over 2}+1
\right)
\Gamma\left({\theta_0\over 2}-{\theta_x\over 2} +{\theta_1\over
  2}-{\theta_\infty\over 2}+1
\right)
}
&
{ 
\Gamma\left(1+{\theta_1\over 2}-{\theta_\infty\over 2}\right)\Gamma(1-\theta_0)
e^{ i{\pi\over 2} \left[
    \theta_x-\theta_0+\theta_\infty-\theta_1\right]}
\over
\Gamma\left(-{\theta_0\over 2}-{\theta_x\over 2} -{\theta_\infty\over
  2}+{\theta_1\over 2}+1
\right)
\Gamma\left({\theta_x\over 2}-{\theta_0\over 2} +{\theta_1\over
  2}-{\theta_\infty\over 2}+1
\right)
}
\cr
\cr
-{ 
\Gamma\left({\theta_\infty\over 2}-{\theta_1\over 2}-1\right)\Gamma(1+\theta_0)
e^{ i{\pi\over 2} \left[
    \theta_0+\theta_x+\theta_1-\theta_\infty\right]}
\over
\Gamma\left({\theta_0\over 2}+{\theta_x\over 2} +{\theta_\infty\over
  2}-{\theta_1\over 2}
\right)
\Gamma\left({\theta_0\over 2}-{\theta_x\over 2} +{\theta_\infty\over
  2}-{\theta_1\over 2}
\right)
}
&
-{ 
\Gamma\left({\theta_\infty\over 2}-{\theta_1\over 2}-1\right)\Gamma(1-\theta_0)
e^{ i{\pi\over 2} \left[
    \theta_x-\theta_0+\theta_1-\theta_\infty\right]}
\over
\Gamma\left(-{\theta_0\over 2}-{\theta_x\over 2} -{\theta_1\over
  2}+{\theta_\infty\over 2}
\right)
\Gamma\left({\theta_x\over 2}-{\theta_0\over 2} +{\theta_\infty\over
  2}-{\theta_1\over 2}
\right)
}
}
\right],
\label{C0inftypasqua}
\ee
\vskip 0.3 cm
\be
C_{01}:=
\left[
\matrix{
{ 
\Gamma(-\theta_x)\Gamma(1+\theta_0)
\over
\Gamma\left({\theta_0\over 2}-{\theta_x\over 2} +{\theta_1\over
  2}-{\theta_\infty\over 2}+1
\right)
\Gamma\left({\theta_0\over 2}-{\theta_x\over 2} +{\theta_\infty\over
  2}-{\theta_1\over 2}
\right)
}
&
{
\Gamma(-\theta_x)\Gamma(1-\theta_0)
\over
\Gamma\left(-{\theta_0\over 2}-{\theta_x\over 2} -{\theta_\infty\over
  2}+{\theta_1\over 2}+1
\right)
\Gamma\left(-{\theta_0\over 2}-{\theta_x\over 2} -{\theta_1\over
  2}+{\theta_\infty\over 2}
\right)
}
\cr
\cr
{
\Gamma(\theta_x)\Gamma(1+\theta_0)
\over
\Gamma\left({\theta_0\over 2}+{\theta_x\over 2} +{\theta_\infty\over
  2}-{\theta_1\over 2}
\right)
\Gamma\left({\theta_0\over 2}+{\theta_x\over 2} +{\theta_1\over
  2}-{\theta_\infty\over 2}+1
\right)
}
&
{
\Gamma(\theta_x)\Gamma(1-\theta_0)
\over
\Gamma\left({\theta_x\over 2}-{\theta_0\over 2} +{\theta_\infty\over
  2}-{\theta_1\over 2}
\right)
\Gamma\left({\theta_x\over 2}-{\theta_0\over 2} +{\theta_1\over
  2}-{\theta_\infty\over 2}+1
\right)
}
}
\right],
\label{C01pasqua}
\ee
\vskip 0.3 cm

 The subgroup generated by $M_0M_x$ and $M_1$ is
reducible. 
As for the solution (\ref{riuffa}), we just need to change
$\theta_1\mapsto -\theta_1$.

\vskip 0.2 cm

b) It is convenient to  re-parameterize the 
solution (\ref{form2}) by introducing  a parameter
$s$ through the equality:
$$
a={\theta_\infty(2s+\theta_x+1)\over 
2(\theta_\infty-1)}.
$$
Let $\theta_x$, $\theta_\infty\not \in{\bf Z}$. 
Then, a representation for the monodromy group is:
$$
M_0= G~\exp\bigl\{
 i\pi\theta_x\sigma_3\bigr\}~G^{-1} ,~~~~~
M_1=\exp\bigl\{-i\pi\theta_\infty\sigma_3\bigr\}
$$
$$
M_x= G~\exp\bigl\{-
 i\pi\theta_x\sigma_3\bigr\}~G^{-1},~~~~~
M_\infty=\exp\bigl\{-i\pi\theta_\infty\sigma_3\bigr\}
$$
In particular, $M_1=M_\infty$, $M_0M_x=I$. 
We can choose $G$ as follows:
$$
G= \pmatrix{ 1 & 1 \cr {s+\theta_x\over r} & {s\over
    r}}. 
%~~\Rightarrow~~
%G^{-1}= {1\over \theta_x} \pmatrix{ -s & r \cr s+\theta_x & -r},
%~~~r,~\theta_x\neq 0.
$$

 Conversely, we may express $s$ as a funcition of the monodromy data: 
$$
s={\theta_x\bigl[2\cos(\pi(\theta_\infty+\theta_x))-\hbox{\rm
    tr}(M_1M_0)\bigr]\over
 2\bigl[\cos(\pi(\theta_\infty-\theta_x))-\cos(\pi(\theta_\infty+\theta_x))\bigr]}.
$$

\vskip 0.2 cm
c) We re-parameterize solution (\ref{form3}) introducing a new
parameter $s$ defined by $
a=:(1-s)^{-1} $. ~Let $\theta_0$, $\theta_x\not \in {\bf Z}$. Then, 
a monodromy representation for the solutions (\ref{form3}) is:
$$
M_0=  \bigl(C_{\infty0}\bigr)^{-1}~
\exp\{i\pi\theta_0\sigma_3\} ~C_{\infty 0},~~~~~
M_{\infty}= \pmatrix{ -1 & 0 \cr 2\pi i \left(1-s\right) & -1 }
$$
$$
M_x=\bigl(C_{\infty0}\bigr)^{-1}\bigl(C_{01}\bigr)^{-1}~
\exp\{i\pi\theta_x\sigma_3 \}~C_{01}C_{\infty0},~~~~~
M_1= \pmatrix{ 1 & 0 \cr 2\pi i s & 1}.
$$
where $C_{\infty0}$ and $C_{01}$ are (\ref{Cinfty0PASQUA}) and
(\ref{C01PASQUA}) given below.  
Conversely, we may express $s$ as a function of the monodromy data:
$$
s= {\hbox{\rm tr}(M_1M_0)-2\cos(\pi \theta_0) \over 4\pi
  \sin(\pi\theta_0)}{(C_{\infty0})_{21}\over (C_{\infty0})_{22}}.
$$

The matrices  $C_{\infty0}$ and $C_{01}$ are:

\be
C_{\infty0}=~2~\pmatrix{
0 
&
{\Gamma(-\theta_0)~e^{-i\pi\{{\theta_0\over 2}+{\theta_x\over
      2}+{3\over 2}\}}
\over 
\Gamma\Bigl(-{\theta_0\over 2}-{\theta_x\over 2} +{3\over 2}\Bigr)
\Gamma\Bigl(-{\theta_0\over 2}+{\theta_x\over 2} +{3\over 2}\Bigr)
}
\cr
\cr
-{\Gamma\Bigl(-{\theta_0\over 2}-{\theta_x\over 2} -{1\over 2}\Bigr)
\Gamma\Bigl(-{\theta_0\over 2}+{\theta_x\over 2} -{1\over 2}\Bigr)
\over 
\Gamma(1-\theta_0) 
e^{-i\pi\{
{\theta_0\over 2}-{\theta_x\over 2}-{3\over 2}
\}}}
&
{\Gamma(\theta_0)~e^{-i\pi\{-{\theta_0\over 2}+{\theta_x\over
      2}+{3\over 2}\}}
\over 
\Gamma\Bigl({\theta_0\over 2}-{\theta_x\over 2} +{3\over 2}\Bigr)
\Gamma\Bigl({\theta_0\over 2}+{\theta_x\over 2} +{3\over 2}\Bigr)
}
}.
\label{Cinfty0PASQUA}
\ee

\vskip 0.3 cm

\be
C_{01}=
\pmatrix{
{\Gamma(-\theta_x)\Gamma(1+\theta_0)
\over 
\Gamma\Bigl({\theta_0\over 2}-{\theta_x\over 2}+{3\over 2}\Bigr)
\Gamma\Bigl({\theta_0\over 2}-{\theta_x\over 2}-{1\over 2}\Bigr)
}
&
{\Gamma(-\theta_x)\Gamma(1-\theta_0)
\over 
\Gamma\Bigl(-{\theta_0\over 2}-{\theta_x\over 2}+{3\over 2}\Bigr)
\Gamma\Bigl(-{\theta_0\over 2}-{\theta_x\over 2}-{1\over 2}\Bigr)
}
\cr
\cr
{\Gamma(\theta_x)\Gamma(1+\theta_0)
\over 
\Gamma\Bigl({\theta_0\over 2}+{\theta_x\over 2}+{3\over 2}\Bigr)
\Gamma\Bigl({\theta_0\over 2}+{\theta_x\over 2}+{1\over 2}\Bigr)
}
&
{\Gamma(\theta_x)\Gamma(1-\theta_0)
\over 
\Gamma\Bigl(-{\theta_0\over 2}+{\theta_x\over 2}+{3\over 2}\Bigr)
\Gamma\Bigl(-{\theta_0\over 2}+{\theta_x\over 2}-{1\over 2}\Bigr)
}
}
\label{C01PASQUA}
\ee

\label{thMONODROMY}
\eth

The conditions  $\theta_\kappa\not\in{\bf Z}$  can be
eliminated, and  the computations can be repeated without conceptual
changes, but with different results.

\vskip 0.2 cm 

In the above theorem, the subgroups
generated by $M_0M_x$ and $M_1$ are {\it reducible}. This characterizes
the monodromy associated to solutions which have a Taylor series at
$x=0$. The same characterization  at $x=1$ involves the subgroup
generated by $M_1M_x$ and $M_0$. At $x=\infty$, it involves  the
subgroup generated by $M_0M_1$ and $M_x$.
\footnote{In the appendix of 
\cite{guz3}, the reader may find explanantions about how to obtain results at
 $x=1,\infty$ from the results  at $x=0$}
In another  paper, we will consider again this characterization, 
 together with the general problem of 
classification. 

\vskip 0.2 cm
 Let us  define again $\sigma$ by tr$(M_0M_x)=2\cos \pi\sigma$. Then, in case
 a), $\sigma= \pm(\theta_1-\theta_\infty)$ [and
 $\pm(\theta_1+\theta_\infty)$ for the change $\theta_1\mapsto
 -\theta_1$]. In case b),  tr$(M_0M_x)=2$ and $\sigma=0$. In case c),    
 tr$(M_0M_x)=-2$, $\sigma=\pm 1$. 
The matching procedure is 
effective to produce solutions corresponding to monodromy data for
which  the connection problem is {\it so far} not well studied, such
 as the case 
tr$(M_iM_j)=-2$ (see \cite{guz4}).\footnote{
 Here I remark that the
 formula (1.30), page 1293, of my paper \cite{guz3} is
 wrong. The correct one is tr$(M_iM_j)\not\in(-\infty,-2]$. In 
\cite{guz3} the connection problem is solved for 
tr$(M_iM_j)\neq \pm 2$. The case 
tr$(M_iM_j)=2$  yields  (\ref{asyjimb0}).  For the special choice of the
 parameters $\theta_0=\theta_x=\theta_1=0$, 
it was studied in \cite{DM}
 and \cite{guz1} (no logarithmic terms
 appear in such a special case). The result (\ref{asyjimb0}) 
for the general (PVI), 
 corresponding to tr$(M_0M_x)=2$, appears in the present paper for the
 first time.}
 
\vskip 0.2 cm 
\noindent
{\it Note:} Also the 1-parameter solutions (\ref{passero})
(\ref{aquila})  and the second
solution in (\ref{asyjimb0}) are 
characterized by a reducible subgroup generated by 
$M_0$, $M_x$. 

\vskip 0.3 cm 
\noindent
{\bf Comments.}

\vskip 0.2 cm 
\noindent
{\bf 1)} The  monodromy group for the solutions (\ref{taylor1})  
was derived  also in
\cite{kaneko}, by confluence of singularities of scalar equations
(including a Heun's type  equation).  The result is equivalent to that
in point a) 
of the above theorem.

\vskip 0.2 cm 
\noindent
{\bf 2)} 
 The computation of the monodromy group
 of the fuchsian systems (\ref{fuchsianSYSTEMOUT})
 and (\ref{fuchsianSYSTEMIN}) is quite clear \cite{Jimbo} \cite{DM}
 \cite{guz3} \cite{Boalch}. It allows
 to express the parameter 
 $r$ of   (\ref{asyjimb}), (\ref{passero}), (\ref{aquila}) 
  and  (\ref{asyjimb0})  as a function of the monodromy
 data. Though the computation for  (\ref{passero}), (\ref{aquila}) 
  and  (\ref{asyjimb0}) does not appear in the literature, the
 procedure is clear (see section
 \ref{senzaconti}), so we do not repeat it.  
We just report the result for  (\ref{asyjimb}), which can be
 found in  \cite{Jimbo} \cite{guz3} \cite{Boalch}:
\be
r=
 {(\theta_0-\theta_x+\sigma)(\theta_0+\theta_x-\sigma)(\theta_\infty+
\theta_1-\sigma)\over
 4\sigma(
 \theta_\infty+\theta_1+\sigma)} ~ {1\over {\bf F}}, 
\label{bfa}
\ee
where
$$
{\bf F}:=  {
\Gamma(1+\sigma)^2\Gamma \left({1\over 2}(\theta_0+\theta_x-\sigma)+1
\right) \Gamma\left(  
{1\over 2} ( \theta_x-\theta_0-\sigma)+1
\right)
\over 
\Gamma(1-\sigma)^2 \Gamma \left({1\over2}(\theta_0+\theta_x+\sigma)+1
\right) \Gamma\left(   
{1\over 2} ( \theta_x-\theta_0+\sigma)+1
\right)
}~\times
$$
$$
\times  {
\Gamma \left({1\over 2}(\theta_\infty+\theta_1-\sigma)+1 \right)
 \Gamma\left(  
{1\over 2} ( \theta_1-\theta_\infty-\sigma)+1
\right)
\over
\Gamma \left({1\over 2}(\theta_\infty+\theta_1+\sigma)+1 \right)
 \Gamma\left( 
{1\over 2} ( \theta_1-\theta_\infty+\sigma)+1
\right)
}~ {V\over U},
$$
and:  
$$
   U:=\left[
{i \over 2}\sin(\pi\sigma)\hbox{tr}(M_1M_x)- \cos(\pi\theta_x) \cos(\pi \theta_\infty) - \cos(\pi\theta_0) \cos(\pi \theta_1) \right]e^{i\pi\sigma}~+$$
$$+
{i\over 2} \sin(\pi\sigma)\hbox{tr}(M_0M_1)+ \cos(\pi\theta_x) \cos(\pi \theta_1) + \cos(\pi\theta_\infty) \cos(\pi \theta_0)
$$
$$ 
V:=4 \sin{\pi\over 2} (\theta_0+\theta_x-\sigma) \sin{\pi\over 2} (\theta_0 - \theta_x+\sigma)~ \sin{\pi\over 2} (\theta_\infty+\theta_1-\sigma)
 \sin{\pi\over 2} (\theta_\infty-\theta_1+\sigma).
$$
The above formula was computed with the assumption that
$\sigma\pm(\theta_0+\theta_x)$, $\sigma\pm(\theta_0-\theta_x)$, 
$\sigma\pm(\theta_1+\theta_\infty)$,
$\sigma\pm(\theta_1-\theta_\infty)$ are not even integers. 
\footnote{ In \cite{guz3} there is a missprint in
  formula (A.30), which must be re-calculated.  In \cite{Jimbo}, in
  formula (1.8) at the bottom of page 1141,  the last 
sign is $\pm\sigma$
  instead of $\mp \sigma$.
}

\vskip 0.2 cm
\noindent
{\bf 3) Reducible Monodromy.} 
The monodromy groups in theorem \ref{thMONODROMY} 
 are not reducible, but they
have a reducible subgroup. If the entire group itself is completely
reducible, all the Painleve\'e transcendents are known. 
 Solutions of (PVI) corresponding to a reducible monodromy were found
in \cite{hitchin}). We summarize the results:  

\bpr 
All the solutions of (PVI) corresponding to a reducible monodromy group 
 are  equivalent by birational canonical transformations 
 to the following one-parameter family of solutions, with 
$\theta_\infty+\theta_1+\theta_0+\theta_x=0$:
\be 
y(x)=
{\theta_1+\theta_\infty-1+x(1+\theta_x)\over \theta_\infty-1}-
{1\over \theta_\infty-1}~{x~(1-x)\over u(x;a)}~{du(x;a)\over dx},
\label{sabishii}
\ee
where $u(x;a)=u_1(x)+a u_2(x)$; $a \in {\bf C}$,  $u_1(x)$ and
$u_2(x)$ 
are linear independent solutions of the hyper-geometric equation: 
$$
x(1-x) {d^2 u\over dx^2} +
\left\{
[2-(\theta_\infty+\theta_1)]-(4-\theta_\infty+\theta_x)x
\right\}{du\over dx} 
-(2-\theta_\infty)(1+\theta_x)u=0
$$
\epr

\vskip 0.3 cm

The monodromy matrices are 
$$
M_0=\pmatrix{{\theta_0\over 2} & * \cr 0 & -{\theta_0\over 2}}
,~~~M_x=\pmatrix{{\theta_x\over 2} & * \cr 0 & -{\theta_x\over 2}},~~~
M_1=\pmatrix{{\theta_1\over 2} & * \cr 0 & -{\theta_1\over 2}}.
$$
The parameter $a$ does not appear in the monodromy.

\vskip 0.3 cm
\noindent
{\it Remark:}  The {\it rational solutions} of 
(PVI) are a special case of the above proposition. They were studied in 
\cite{mazzoccoratio}. 
 Up to canonical birational transformations, 
they are realized for $\theta_\infty+\theta_1+\theta_0+\theta_x=0$ and:
$$
\theta_0=1:~~~~~y(x)= {\theta_\infty+\theta_1\over \theta_\infty}~ {x-1\over
x(1+\theta_1)-(\theta_1+\theta_\infty)};
$$
\vskip 0.2 cm
$$
\theta_0=-2:~~~~~y(x)={\bigl(2-(\theta_\infty+\theta_1)+\theta_1~x\bigr)^2-2+
\theta_\infty+\theta_1-\theta_1~x^2
\over (1-\theta_\infty)\bigl(2-(\theta_\infty+\theta_1)+\theta_1~x\bigr)}.
$$

\vskip 0.3 cm
The computation of the expansion at $x=0$ of (\ref{sabishii}) is just
a consequence of the expansions of $u_1(x)$ and $u_2(x)$. The reader
can find by himself a behavior  $y\sim x(r(a)\pm \theta_x\ln(x))$ for 
$\theta_1+\theta_\infty=\theta_0+\theta_x=0$, namely a sub-case of 
the second solution in  (\ref{asyjimb0}). 
For $\theta_1+\theta_\infty\not\in {\bf Z}$, we
 find behaviors of the type  (\ref{1parameter}) (and (\ref{form1}),
 (\ref{taylor1}) for 
 $a=0$).

\vskip 1 cm
\centerline{\Large \bf  PART II --  Derivation  Results:  Fuchsian Reduction}

\section{Fuchsian Case}
Let $x\to 0$. The reduction to the fuchsian systems (\ref{fuchsianSYSTEMOUT}) is
possible if  in the domain (\ref{dominioOUTbasta}) we have: 
\be
 |(A_0+A_x)_{ij}| \gg \left|(A_x)_{ij}~{x\over \lambda}\right|,~~~\hbox{
   namely: }~
|(A_0+A_x)_{ij}| \gg \left|(A_x)_{ij}~x^{1-\delta_{OUT}}\right|.
\label{condition1}
\ee
Let us denote with $\hat{A}_i$ the leading term of the matrix $A_i$,
$i=0,x,1$. We can  substitute 
(\ref{fuchsianSYSTEMOUT}) with: 
\be
{d\Psi_{OUT}\over d\lambda} = \left[{\hat{A}_0+\hat{A}_x\over
    \lambda}+ 
{\hat{A}_1\over \lambda-1}\right]~\Psi_{OUT}
\label{system1}
\ee
 We suppose that $\theta_\infty\neq 0$. This is not a loss
in generality, because 
$\theta_\infty=0$ is equivalent to $\theta_\infty=2$.
% LEMMA

\ble
   If the approximation 
(\ref{fuchsianSYSTEMOUT}) is possible, then $ 
  \hat{A}_0+\hat{A}_x $ has eigenvalues $\pm {\sigma\over 2}\in {\bf
  C}$  independent of $x$, defined (up to sign and addition 
of an integer) by $
\hbox{\rm tr}(M_xM_0)= 2\cos(\pi \sigma)
$. Let $r_1\in{\bf C}$, $r_1\neq 0$. 
For $\theta_\infty\neq 0$, the leading terms are:
\be
\hat{A_1}= 
\pmatrix{ {\sigma^2-\theta_\infty^2-\theta_1^2\over 4 \theta_\infty} & 
-r_1
\cr
{[\sigma^2-(\theta_1-\theta_\infty)^2][\sigma^2-(\theta_1+\theta_\infty)^2]\over
  16 \theta_\infty^2 }~{1\over r_1} 
&
- 
 {\sigma^2-\theta_\infty^2-\theta_1^2\over 4 \theta_\infty} 
},
\label{hatA1}
\ee
and
\be
\hat{A_0}+\hat{A_x}=
\pmatrix{
{\theta_1^2-\sigma^2-\theta_\infty^2\over 4 \theta_\infty}
&
r_1
\cr
-{[\sigma^2-(\theta_1-\theta_\infty)^2][\sigma^2-(\theta_1+\theta_\infty)^2]\over
  16 \theta_\infty^2 }~{1\over r_1}
& 
-{\theta_1^2-\sigma^2-\theta_\infty^2\over 4 \theta_\infty}
}.
\label{hatA0Ax}
\ee
\label{elicopter1}
\ele

\vskip 0.2 cm
\noindent
{\it Proof:}  Observe that $
\hbox{tr}(\hat{A}_0+\hat{A}_x)=~\hbox{tr}(A_0+A_x)=0$, 
thus, for any $x$,  
$\hat{A}_0+\hat{A}_x$ has eigenvalues of opposite sign, that we
denote $\pm \tilde{\sigma}(x)/2$.    
Then, we recall that 
 $x$ is a monodromy preserving deformation, therefore the monodromy 
matrices 
of (\ref{system1}) are  
independent of $x$.  At $\lambda=0,1,\infty$ they are:
$$
   M^{OUT}_0=\left\{\matrix{M_xM_0
                                     \cr 
                            M_0M_x}
\right.
, ~~~M^{OUT}_1=M_1,~~~M^{OUT}_{\infty}=M_{\infty}.
$$
Thus,  det$(M^{OUT}_0)=1$, because
det$(M_x)$=det$(M_0)=1$. Therefore, there exists a constant matrix $D$
and a complex constant number $\sigma$  
such that:  
$$ 
D^{-1}~M^{OUT}_0~D~=\left\{ \matrix{
 \hbox{diag}(\exp\{-i\pi \sigma\}, \exp\{i\pi \sigma\})
,\cr
\cr
 \pmatrix{\pm 1 & * \cr 0 & \pm 1} ,\hbox{ or } \pmatrix{\pm 1 & 0 \cr * & \pm 1},~~~\sigma\in{\bf Z}
}
\right. 
$$ 
We conclude that $\tilde{\sigma}(x)\equiv \sigma$. We also have 
 $\hbox{tr}(M_0^{OUT})= 2\cos(\pi\sigma)$.

\vskip 0.2 cm 
 Now consider the gauge:  
\be
\Phi_1 := \lambda^{-{\sigma\over 2}} (\lambda-1)^{-{\theta_1\over 2}} 
~\Psi_{OUT}.~~~~~
{d\Phi_1\over d\lambda} = \left[{\hat{A}_0+\hat{A}_x-{\sigma\over 2}\over \lambda}+ {\hat{A}_1-{\theta_1\over 2}\over \lambda-1}\right]~\Phi_1
\label{systemPhi1}
\ee
%$$
%\hbox{Eigenvalues }(\hat{A}_0+\hat{A}_x-{\sigma\over 2})=0,~-\sigma,
%$$
%$$
%\hbox{Eigenvalues }(\hat{A}_1-{\theta_1\over 2})=0,~-\theta_1,
%$$
%$$
%\hat{A}_0+\hat{A}_x+\hat{A}_1-{\theta_1\over 2}-{\sigma\over 2}
%=\pmatrix{-{\theta_\infty\over 2} -{\theta_1\over 2}-{\sigma\over 2} &
%0\cr
%0&   {\theta_\infty\over 2} -{\theta_1\over 2}-{\sigma\over 2}}.
%$$
We can identify $\hat{A}_0+\hat{A}_x-{\sigma\over 2}$ and
  $\hat{A}_1-{\theta_1\over 2}$ with $B_0$ and $B_1$  of  Proposition
 \ref{matrices} in Appendix 1, case  (\ref{1}), with $
a=  {\theta_\infty\over 2} +{\theta_1\over 2}+{\sigma\over 2}$, 
        $ b= - {\theta_\infty\over 2} +{\theta_1\over 2}+{\sigma\over
  2}$, $
         c= \sigma$. 
 \qed

\vskip 0.3 cm
\noindent
{\it Remark:} 
$r_1$ may be a function of $x$. If the monodromy of system (\ref{system1}) 
depends on $r_1$, then $r_1$ is a constant independent of $x$. This is
the case here, but we do not prove it for reasons of space. See the
references in Part I.

\vskip 0.3 cm
 Lemma \ref{elicopter1} (and Lemma  \ref{genericfuchs1} which 
 follows)
 includes 
 all cases (\ref{1})---(\ref{5}) for system (\ref{systemPhi1}). 
Cases  
(\ref{2})---(\ref{5}) are obtained substituting 
$\sigma=-(\theta_\infty+\theta_1)$, $\theta_\infty-\theta_1$,
 $\theta_\infty+\theta_1$, $\theta_1-\theta_\infty$ respectively. 
 For all the computations which  follow, 
 involving system (\ref{system1}) or (\ref{systemPhi1}), we note that 
the hypothesis $\theta_\infty\neq 0$ excludes  cases (\ref{6}),  (\ref{7}) and the Jordan cases (\ref{8})--(\ref{10}).

\vskip 0.5 cm

%%%%%%%%%%%%%%%%%%%%%%%%%%%%%%%%%%%%%%%%%%%%%%%%%%%%%%%

The reduction to the fuchsian system (\ref{fuchsianSYSTEMIN}) is
possible for $x\to 0$ in the domain (\ref{dominioINbasta}) if: 
\be
\left|{(A_0)_{ij}\over \lambda}+{(A_x)_{ij}\over \lambda-x}
\right|
\gg
\left|(A_1)_{ij}\right|,
~~~\hbox{ namely:}~
\left|
{(A_0+A_x)_{ij}\over x^{\delta_{IN}}} 
\right| \gg 
\left|(A_1)_{ij}\right|.
\label{condition0}
\ee
We can rewrite   (\ref{fuchsianSYSTEMIN}) using just   the leading terms of
the matrices: 
\be
{d\Psi_{IN}\over d\lambda}=\left[
{\hat{A_0}\over \lambda}+{\hat{A_x}\over \lambda-x}\right] \Psi_{IN},
\label{systemapprox0}
\ee
Then, we re-scale $\lambda$ and consider the following system:  
$$
{d \Psi_{IN} \over d\mu} = \left( 
         {\hat{A}_0\over \mu} +{\hat{A}_x \over \mu-1} 
\right) \Psi_{IN},~~~~~\mu:={\lambda\over x}
$$
We know that there exists  a    
matrix $K_0(x)$ such that: 
$$
  {K_0}^{-1}(x)~(\hat{A_0}+\hat{A}_x)~K_0(x)= \pmatrix{{\sigma\over 2}& 0 \cr 
                                                0 & -{\sigma\over 2}}
,~~\hbox{ or }~\pmatrix{0 & 1 \cr 0 & 0 }.
$$
Let $
\hat{\hat{A_i}} := {K_0}^{-1} \hat{A}_i K_0$, $i=0,x$. By a gauge
transformation, we get the system: 
\be
\Psi_{IN}=: K_0(x)~\Psi_0,
~~~~~~~
{d\Psi_0\over d \mu}=\left[{\hat{\hat{A_0}}\over
    \mu}+{\hat{\hat{A_x}}\over \mu-1} 
 \right] \Psi_0,
\label{system0}
\ee

\ble
Let $r\in{\bf C}$, $r\neq 0$. 
If $\sigma\neq 0$, we have:  
\be
\hat{\hat{A_0}}
=
\pmatrix{{\theta_0^2-\theta_x^2+\sigma^2\over 4\sigma}
&
r
\cr
-{[\sigma^2-(\theta_0-\theta_x)^2][\sigma^2-(\theta_0+\theta_x)^2]
\over 16 \sigma^2 }~{1\over r_0}
&
-{\theta_0^2-\theta_x^2+\sigma^2\over 4\sigma}
},
\label{hathatA0}
\ee
\be
\hat{\hat{A_x}}= \pmatrix{ {\sigma^2+\theta_x^2-\theta_0^2\over 4 \sigma} 
&
-r
\cr
{[\sigma^2-(\theta_0-\theta_x)^2][\sigma^2-(\theta_0+\theta_x)^2]\over
  16 \sigma^2 }~{1\over r}
&
- {\sigma^2+\theta_x^2-\theta_0^2\over 4 \sigma} 
}.
\label{hathatAx}
\ee
\label{elicopter0}
\ele

\vskip 0.3 cm
\noindent
{\it Proof:} 
We do a gauge transformation:  
\be
\Phi_0:= \mu^{-{\theta_0\over 2}} (\mu-1)^{-{\theta_x\over 2}} ~\Psi_0,~~~~~
{d \Phi_0\over d\mu}= \left[ {\hat{\hat{A_0}}-{\theta_0\over 2}\over 
\mu} + {\hat{\hat{A_x}}-{\theta_x\over 2}\over \mu-1} \right] ~\Phi_0. 
\label{systemPhi0}
\ee
We 
identify $\hat{\hat{A_0}}-{\theta_0\over 2}$, 
$\hat{\hat{A_x}}-{\theta_x\over 2}$ 
with $B_0$ and $B_1$ in the Appendix 1, Proposition
\ref{matrices}, case (\ref{1}),  with $
  a={\theta_0\over 2} + {\theta_x\over 2} -{\sigma\over 2}$,  
$b= {\theta_0\over 2}+{\theta_x\over 2} + {\sigma\over 2}$, 
$c=\theta_0$. 
\qed

\vskip 0.3 cm

\noindent
{\it Remark:} If the monodromy of the system ({\ref{system0}) 
depends on $r$, then $r$ is a constant independent of $x$. This is the
case here.

\vskip 0.3 cm

The  Lemma \ref{elicopter0}  (and lemma \ref{genericfuchs0}
 which follows) includes also the cases 
 (\ref{2})--(\ref{5}) for the system
 (\ref{systemPhi0}). 
 These cases correspond respectively to the values
 $\sigma=\theta_0+\theta_x$, $-\theta_0-\theta_x$,
 $\theta_x-\theta_0$, $\theta_0-\theta_x$, with $\theta_0\neq \pm
 \theta_x$.

%%%%%%%%%%%%%%%%%%%%%%%%%%%%%%%%%%%%%%%%

\subsection{Matching for $\sigma\not \in {\bf Z}$ and proof of
 (\ref{asyjimb})} 

We match $\Psi_{OUT}$ and $\Psi_{IN}$ in the intersection of the
``outside'' and ``inside'' domains, namely the  
region $|x|^{\delta_{OUT}}\leq |\lambda|\leq |x|^{\delta_{IN}}$, $x\to
0$.  
As a consequence, 
we obtain the leading term of $y(x)$. 

\ble
 If $\sigma\not\in{\bf Z}$ and $\theta_\infty\neq 0$, system
 (\ref{system1}) has a fundamental matrix solution 
 $\Psi_{OUT}(\lambda)$ with the following behavior at $\lambda=0$: 
$$\Psi_{OUT}= \sum_{n=0}^\infty G_n\lambda^n\pmatrix{
\lambda^{{\sigma\over 2}} & 0 \cr
0 & \lambda^{-{\sigma\over 2}}
},~~~~~
G_0=
\pmatrix{1 & 1 
             \cr
       {(\theta_\infty+\sigma)^2 - \theta_1^2 \over 4 \theta_\infty r_1} 
&
  {(\theta_\infty-\sigma)^2-\theta_1^2\over 4 \theta_\infty r_1}
}.
$$ 
$G_n$ are matrices which depend
rationally on $\theta_\infty$, $\theta_1$, $\sigma$, $r_1$. 
The series is convergent for $|\lambda|<1$. 
\label{genericfuchs1}
\ele
\noindent
{\it Proof:} It is an immediate consequence of the standard theory of 
linear systems of fuchsian differential equations. 
\qed

\vskip 0.3 cm

\ble
If $\sigma \not \in {\bf Z}$,  system 
(\ref{system0}) has a fundamental matrix solution  with the following
behavior at $\mu=\infty$: 
$$
\Psi_0(\mu)=
\left[I+
\sum_{n=1}^{\infty} K_n \mu^{-n}\right]
~\pmatrix{ \mu^{\sigma\over 2}
&
0
\cr
0
&
\mu^{-{\sigma\over 2} }
},
$$
where $I$ is the identity matrix, $K_n$ are matrices which depend
rationally on $\theta_0$, $\theta_x$, $\sigma$, $r$. The series is
convergent for $|\mu|>1$. 
\label{genericfuchs0}
\ele

\noindent
{\it Proof:} It is a consequence of the standard theory of systems of
fuchsian  equations.  
\qed

%%%%%%%%%%%%%%%%%%%%%%%%%%%%%%%%%%%%%%%%%%%%%%%%

\vskip 0.3 cm
The matching relation $
 \Psi_1(\lambda)\sim 
 K_0(x) \Psi_0\left({\lambda/ 
 x}\right)$, $|x|^{\delta_{OUT}}\leq |\lambda|\leq
 |x|^{\delta_{IN}}$, $x\to 
0$, is: 
$$
 G_0 \pmatrix{ \lambda^{{\sigma\over 2}} & 0 
\cr
0 &\lambda^{-{\sigma\over 2}}}~
\sim
~
K_0(x)  \pmatrix{ \lambda^{{\sigma\over 2}} & 0 
\cr
0 &\lambda^{-{\sigma\over 2}}}~ \pmatrix{ x^{-{\sigma\over 2}} & 0 
\cr
0 &x^{{\sigma\over 2}}}.  
$$
 This gives the result:
$$
K_0(x)\sim \pmatrix{1 & 1 
             \cr
       {(\theta_\infty+\sigma)^2 - \theta_1^2 \over 4 \theta_\infty r_1} 
&
  {(\theta_\infty-\sigma)^2-\theta_1^2\over 4 \theta_\infty r_1}
}
 \pmatrix{x^{{\sigma\over 2}} & 0 
\cr
0 & x^{-{\sigma\over 2}}
}.
$$
%%%%%%%%%%%%%%%%%%%%%%%%%%%%%%%
%%%%%%%%%%%%%%%%
%%%%%%%%%%%       CANCELED IN THE PREPRINT
%%%%%%%%%%%%%%%%
%%%%%%%%%%%%%%%%%%%%%%%%%%%%%%%
%and
%$$
%{K_0(x)}^{-1}=
%\pmatrix{x^{-{\sigma\over 2}} & 0
%\cr
%0 & x^{{\sigma\over 2}}
%}
%~\pmatrix{  {\theta_1^2-(\theta_\infty-\sigma)^2\over 4\theta_\infty \sigma}
%&
%{r_1\over \sigma} 
%\cr
%{(\theta_\infty+\sigma)^2-\theta_1^2\over 4 \theta_\infty \sigma} 
%&
%-{r_1\over \sigma}
%}.
%$$
%%%%%%%%%%%%%%%%%%%%%%%%%%%%%%%%
We  compute the matrices $
\hat{A_0}(x)=K_0(x) \hat{\hat{A_0}}{K_0(x)}^{-1}$, 
$\hat{A_x}(x)=K_0(x) \hat{\hat{A_x}} {K_0(x)}^{-1}
$
making use of  Lemma \ref{elicopter0}. We obtain: 
$$
\hat{A_0}(x)= G_0 
\pmatrix{{\theta_0^2-\theta_x^2+\sigma^2\over 4\sigma}
&
r ~x^{\sigma}
\cr
-{(\sigma+\theta_x-\theta_0)(\sigma+\theta_x+\theta_0)(\sigma-\theta_x+\theta_0)(\sigma-\theta_x-\theta_0)\over 16 \sigma^2 r}~x^{-\sigma}
&
-{\theta_0^2-\theta_x^2+\sigma^2\over 4\sigma}
}~{G_0}^{-1}
$$
$$
\hat{A_x}(x)= G_0\pmatrix{ {\sigma^2+\theta_x^2-\theta_0^2\over 4 
\sigma} &
-r~x^{\sigma} 
\cr
{(\sigma+\theta_x-\theta_0)(\sigma+\theta_x+\theta_0)(\sigma-\theta_x+\theta_0)(\sigma-\theta_x-\theta_0)\over 16 \sigma^2 r}~x^{-\sigma}
&
- {\sigma^2+\theta_x^2-\theta_0^2\over 4 \sigma} 
}
{G_0}^{-1}.
$$
This result shows that the matrix elements of 
 $\hat{A_0}$ and $\hat{A_x}$ diverge as $|x|^{-|\Re \sigma|}$ when 
 $x\to 0$ inside a sector (i.e. for $|$arg$(x)|$  bounded). 
In particular, we find $
(\hat{A_1})_{12}=-r_1$ and
{\small 
$$
(\hat{A_0})_{12}= {r_1\over r}~{
[\sigma^2-(\theta_0+\theta_x)^2][(\theta_0-\theta_x)^2-\sigma^2]
\over 16 \sigma^3 } ~x^{-\sigma}
~+~r_1~{\theta_0^2-\theta_x^2+\sigma^2\over 2 \sigma^2}  ~-~{r r_1 \over \sigma} ~x^{\sigma}.
$$
}
The above are enough to compute the 
 leading term(s) of $y(x)$ from the
 formula:
\be
y(x)= {x (A_0)_{12}\over x[(A_0)_{12}+(A_1)_{12}]-(A_1)_{12}}= 
 -{x(A_0)_{12}\over (A_1)_{12}}~\left[ 1~ - ~x~\left(
1+ {(A_0)_{12}\over (A_1)_{12}}
\right)\right]^{-1} 
\label{leadingterm}
\ee
Thus:
{\small
\be
y(x)\sim -x ~ {(\hat{A_0})_{12}\over (\hat{A_1})_{12}}
= 
\left[
{1\over r}~{[\sigma^2-(\theta_0+\theta_x)^2][(\theta_0-\theta_x)^2-\sigma^2]
\over 16 \sigma^3}~x^{1-\sigma}~+{\theta_0^2-\theta_x^2+\sigma^2\over 
2\sigma^2}~x~-~{r\over \sigma}~x^{1+\sigma}
\right].
\label{genericlead}
\ee
}
We have ignored $1/[1-x(1+(\hat{A_0})_{12}/(\hat{A_1})_{12})]$ 
because condition (\ref{condition1}) is  equivalent to:  $
\left| x^{1-\delta_1\pm\sigma}\right| \to 0$ for $x\to 0$, 
 which implies that $
\left| x~{(\hat{A_0})_{12}/ (\hat{A_1})_{12}} \right|\sim
\left|x^{1\pm \sigma}  \right|\to 0$. Therefore,
  $1/[1-x(1+(\hat{A_0})_{12}/(\hat{A_1})_{12})]=1/(1+O(x))=1+O(x)$.

If $\Re \sigma\neq 0$, 
the leading term of (\ref{genericlead}) is certainly correct, but some
higher
 order corrections may be bigger than the next two terms of 
 (\ref{genericlead}).     
 If $\Re\sigma=0$, the three terms of
  (\ref{genericlead}) are of the same order, and their combination gives the 
trigonometric expression in theorem  \ref{thsigmano1}.

%%%%%%%%%%%%%%%%%%%%%%%%
\subsection{Range of $\sigma$}

 Conditions (\ref{condition1}) and (\ref{condition0}) must be verified. Let 
$C$ denote a non zero constant. We suppose that $x\to 0$ inside a
 sector with center on $x=0$.  Then:
$$
\hbox{Condition (\ref{condition0}) is }
|x|^{-\delta_{IN}}\gg C~~ \Longleftrightarrow~~ \delta_{IN}>0.
$$
$$
 \hbox{Condition (\ref{condition1}) is }
C \gg |x|^{-|\Re\sigma|+1-\delta_{OUT}}~~\Longleftrightarrow~~  |\Re \sigma|<1-\delta_{OUT}.
$$
 The last condition  implies that $|\Re \sigma|<1$. We also conclude
 that 
 $0<\delta_{IN}\leq \delta_{OUT}<1$. 

%%%%%%%%%%%%%%%%%%%%%%%%%
%%%%%%%%%%%%%%
%%%%%%%%%%      CANCELED IN THE PREPRINT
%%%%%%%%%%%%%%
%%%%%%%%%%%%%%%%%%%%%%%%%
% The above discussion holds 
%when $x\to 0$ inside a sector with center on $x=0$. We can 
%be more general, allowing $x$ to converge to $x=0$ along any
% path in the covering of a neighborhood 
%of $x=0$ (with $x=0$ removed). Therefore $|x^{\sigma}|=\exp\{ \Re
% \sigma \log|x| - \Im \sigma \arg(x)\}|$ and the condition
% (\ref{condition1}) reads: 
%$$
 %\hbox{Condition (\ref{condition1}) }
%\Longleftrightarrow 
%\left\{
%\matrix{
%\exp\{(\Re\sigma+1-\delta_{OUT})\log|x|  - \Im\sigma \arg (x)\} \ll C
%\cr
%\cr
%\exp\{(-\Re\sigma+1-\delta_{OUT})\log|x|  + \Im\sigma \arg (x)\} \ll C
%}\right.
%$$
%This is verified when the exponents are [...].
%%%%%%%%%%%%%%%%%%%%%%%%%%%%%%%%%%%%%%%%%%%%%%%%%%%%%

%%%%%%%%%%%%%%%%%%%

\subsection{Leading term for $\sigma=\pm(\theta_0+\theta_x),
  \pm(\theta_0-\theta_x)~\neq~0$. Proof of (\ref{passero}) and
  (\ref{aquila})} 

Formula (\ref{genericlead}) holds  for any  $\sigma\neq 0$ such that
 $|\Re \sigma|<1$. However, we cannot naively substitute the value of
 $\sigma=\pm(\theta_0+\theta_x),~ \pm(\theta_0-\theta_x)$, for which
 the coefficient of $x^{1-\sigma}$ vanishes. This 
 is  
 because
 only the leading term is certainly correct, and it may be the term in
 $x^{1-\sigma}$. 
  Therefore, 
here we briefly give the explicit derivation of  (\ref{passero}) and
 (\ref{aquila})},  
 using cases (\ref{2})--(\ref{5}) for system (\ref{systemPhi0}).  
% We also comment on cases (\ref{6}) and (\ref{7}).

\vskip 0.2 cm
\noindent
{\it Case (\ref{2}), $a=0$:} This is the case $
\sigma=\theta_0+\theta_x\neq 0$. The matching procedure does not
change.   From (\ref{2}) we compute:
$$
\hat{A_0}=G_0~\pmatrix{ {\theta_0\over 2} & r~x^{\sigma} \cr
                          0 & -{\theta_0\over 2}
}~{G_0}^{-1},
~~~
\hat{A_x}=
G_0~\pmatrix{ {\theta_x\over 2} & -r~x^{\sigma} \cr
                          0 & -{\theta_x\over 2}
}~{G_0}^{-1}.
$$
This implies that $
(\hat{A_0})_{12}= r_1 \left(
{\theta_0\over \theta_0+\theta_x}-{r\over \theta_0+\theta_x}~x^{\sigma}
\right),$
while $(\hat{A_1})_{12}=-r_1$ as in the generic case. Therefore, 
$$
y(x)\sim {\theta_0\over \theta_0+\theta_x}~x~-{r\over
  \theta_0+\theta_x}~x^{\sigma+1}. 
$$
It is interesting to note that for $\Re \sigma>0$ we have $\hat{y}(x)
\sim  {\theta_0/(\theta_0+\theta_x)}~x$. Such a behavior is what one
would naively expect from the generic behavior (\ref{genericlead})
when $\sigma=0$.   

\vskip 0.2 cm 
 Case (\ref{3}), $b=0$, is  $\sigma=-\theta_0-\theta_x\neq 0$. 
Case (\ref{4}), $a=c$, is $\sigma=\theta_x-\theta_0$.  Case (\ref{5}),
 $b=c$, is $\sigma=\theta_0-\theta_x$. Proceeding as above, we find 
 (\ref{passero}) and (\ref{aquila}). 

%%%%%%%%%%%%%%%%%%%%%%%%%
%%%%%%%%%%%%%%%%
%%%%%%%%%%%%%%    CANCELED IN THE PREPRINT
%%%%%%%%%%%%%%%%
%%%%%%%%%%%%%%%%%%%%%%%%%
%
%
%\vskip 0.2 cm
%\noindent
%{\it Case (\ref{3}), $b=0$}: In this case 
% $$\sigma=-\theta_0-\theta_x\neq 0.
%$$
% From (\ref{3}) we compute:
%$$
%\hat{\hat{A_0}}=\pmatrix{-{\theta_0\over 2} & r \cr 
%                          0 & {\theta_0\over 2}},
%~~~~~
%\hat{\hat{A_x}}=\pmatrix{-{\theta_x\over 2} & -r \cr 
%                          0 & {\theta_x\over 2}}. 
%$$
%In the same way of the case (\ref{2}) we get:
%$$
%\hat{y}(x)= {\theta_0\over \theta_0+\theta_x} ~x +~{r\over
% \theta_0+\theta_x}~x^{1+\sigma}. 
%$$
%
%
%\vskip 0.2 cm
%\noindent
%{\it Case (\ref{4}), $a=c$:} In this case 
%$$
%\sigma=\theta_x-\theta_0\neq 0.
%$$
% We compute:
%$$
%\hat{\hat{A_0}}=\pmatrix{-{\theta_0\over 2} & r \cr 
%                          0 & {\theta_0\over 2}},
%~~~~~
%\hat{\hat{A_x}}=\pmatrix{{\theta_x\over 2} & -r \cr 
%                          0 & -{\theta_x\over 2}}. 
%$$
%In the same way of the case (\ref{2}) we get:
%$$
%\hat{y}(x)= {\theta_0\over \theta_0-\theta_x} ~x +~{r\over
%\theta_0-\theta_x}~x^{1+\sigma}. 
%$$
%
%\vskip 0.2 cm
%\noindent
%{\it Case (\ref{5}), $b=c$:} In this case  
%$$
%\sigma=\theta_0-\theta_x\neq 0.
%$$
%We find:
%$$
%\hat{\hat{A_0}}=\pmatrix{{\theta_0\over 2} & r \cr 
%                         0 & -{\theta_0\over 2}},
%~~~~~
%\hat{\hat{A_x}}=\pmatrix{-{\theta_x\over 2} & -r \cr 
%                          0 & {\theta_x\over 2}}. 
%$$
%In the same way of the case (\ref{2}) we get:
%$$
%\hat{y}(x)= {\theta_0\over \theta_0-\theta_x} ~x -~{r\over
%\theta_0-\theta_x}~x^{1+\sigma}. 
%$$

\vskip 0.3 cm
\noindent
{\it Remark:} If we substitute $y=b_1x+b_2x^2+b_3x^3+...$ into (PVI) we find 
all the coefficients $b_n$ by identifying equal powers of $x$. The
result is (\ref{taylor1}).   
We need to assume that $\theta_0\pm\theta_x$ is not 
integer or zero.

\subsection{Matching for $\sigma=0$. Proof of (\ref{asyjimb0})}

\subsection{ Case $\theta_0\pm \theta_x\neq 0$}
\label{secsigma0}
\ble
Let $r_1\in{\bf C}$, $r_1\neq 0$. 
The matrices of system (\ref{system1}) are:  
$$
\hat{A_1}= 
\pmatrix{ -{{\theta_\infty}^2+{\theta_1}^2\over 4 \theta_\infty} & 
-r_1
\cr
{[\theta_1^2-\theta_\infty^2]^2\over 16 \theta_\infty^2 r_1} 
&
 {{\theta_\infty}^2+{\theta_1}^2\over 4 \theta_\infty} 
},~~~
\hat{A_0}+\hat{A_x}=
\pmatrix{
{\theta_1^2-\theta_\infty^2\over 4 \theta_\infty}
&
r_1
\cr
-{[\theta_\infty^2-\theta_1^2]^2\over 16 \theta_\infty^2 r_1} 
& 
{\theta_\infty^2-\theta_1^2\over 4 \theta_\infty}
}, ~~~~\forall r_1\neq 0.
$$
A fundamental matrix solution can be chosen with the following
behavior at $\lambda=0$: 
$$
\Psi_{OUT}(\lambda)=[G_0 + O(\lambda)]~\pmatrix{1 & \log \lambda \cr
0 & 1},~~~~~~~
G_0= \pmatrix{ 1 & 0 \cr
                     {{\theta_\infty}^2-{\theta_1}^2\over 4 \theta_\infty~r_1}
&
{1\over r_1}
}.
$$
\ele

\noindent
{\it Proof:}  The system (\ref{systemPhi1}) is: 
 $$
{d\Phi_1\over d\lambda} = \left[{\hat{A}_0+\hat{A}_x\over \lambda}+ {\hat{A}_1-{\theta_1\over 2}\over \lambda-1}\right]~\Phi_1, 
$$
%$$
%\hbox{Eigenvalues }(\hat{A}_0+\hat{A}_x)=0,~0
%.~~~~~\hbox{Eigenvalues }(\hat{A}_1-{\theta_1\over 2})=0,~-\theta_1.
%$$
%$$
%\hat{A}_0+\hat{A}_x+\hat{A}_1-{\theta_1\over 2}
%=\pmatrix{-{\theta_\infty\over 2} -{\theta_1\over 2} &
%0\cr
%0&   {\theta_\infty\over 2} -{\theta_1\over 2}}.
%$$
 We identify 
 $\hat{A}_0+\hat{A}_x$ and $\hat{A}_1-{\theta_1\over 2}$ with $B_0$
 and $B_1$ of   proposition
 \ref{matrices} in Appendix 1, diagonalizable case (\ref{1})--(\ref{5}) 
(we recall that (\ref{6})--(\ref{10}) never occur when $\theta_\infty\neq 0$) 
with $
a=  {\theta_\infty\over 2} +{\theta_1\over 2}$, 
      $   b= - {\theta_\infty\over 2} +{\theta_1\over 2}$, 
          $c= 0$. 

The behavior of a fundamental solution 
is a standard result in the theory of Fuchsian
systems. The matrix $G_0$ is defined by $
{G_0}^{-1} \left(\hat{A_0}+\hat{A_x}\right)G_0=
\pmatrix{ 0 & 1 \cr 0 & 0}
$. 
\qed

\ble
Let $r\in{\bf C}$. The matrices of system (\ref{system0}) are:
\be
\hat{\hat{A_0}}= \pmatrix{ r +{\theta_0\over 2} &
{4~r~(r+\theta_0) \over \theta_x^2-\theta_0^2} \cr
{\theta_0^2-\theta_x^2\over 4} & -r-{\theta_0\over 2}
},~~~~~
\hat{\hat{A_x}}=\pmatrix{
-r-{\theta_0\over 2}  & 1- {4~r~(r+\theta_0) \over \theta_x^2-\theta_0^2}
\cr
{\theta_x^2-\theta_0^2\over 4} 
&
r+{\theta_0\over 2}.
}.
\label{saravero?0}
\ee
There exist a fundamental solution of (\ref{system0}) with the
following behavior at $\mu=\infty$:
$$
\Psi_0(\mu)=
\left[I+O\left({1\over \mu}\right)\right]~
\pmatrix{ 1 & \log \mu \cr 0 &1 }, ~~~\mu\to\infty.
$$
\ele

\noindent
{\it Proof:} 
To compute $\hat{\hat{A_0}}$ and $\hat{\hat{A_x}}$ for the generic
case,  we consider the  case (\ref{8}) in Proposition \ref{matrices},
applied to the system(\ref{systemPhi0}). The parameters are $
a={\theta_0\over 2}+{\theta_x\over 2}$, $c= \theta_0$. In particular, 
 \be
\hat{\hat{A_0}}-{\theta_0\over 2} +\hat{\hat{A_x}}-{\theta_x\over 2}
=
\pmatrix{-{\theta_0+\theta_x\over 2} & 1 \cr 0 & -{\theta_0+\theta_x\over 2} }
\label{matrinf0}
\ee
Here the values of the  
parameters satisfy the conditions  
  $a\neq 0$ and $ a\neq c$,  namely 
 $\theta_0\pm\theta_x\neq 0$. From the  matrices (\ref{8}), we obtain
$\hat{\hat{A_0}}=B_0+\theta_0/2$ and 
$\hat{\hat{A_x}}=B_1+\theta_x/2$. 
Keeping into account (\ref{matrinf0}),  by the standard theory
of fuchsian systems we have:  
$$
\Phi_0(\mu)=\left[I+O\left({1\over \mu}\right)\right]~\mu^{-{\theta_0+\theta_x\over 2}}~
\pmatrix{ 1 & \log \mu \cr 0 &1 }, ~~~\mu\to\infty.
$$
This proves the behavior of $\Psi_0(\mu)$. 
\qed

\vskip 0.3 cm
The matching condition  
$
\Psi_{OUT}(\lambda)\sim
K_0(x) ~\Psi_0\left(\lambda/ x\right)$ becomes: 
{\small
$$
K_0(x)~ \pmatrix{1 & \log\left({\lambda\over x}\right) \cr
0 & 1 }~\sim~
G_0~\pmatrix{1 & \log \lambda 
 \cr 
0 & 1 }~~~
\Longrightarrow
~~~
K_0(x)\sim \pmatrix{1 & 0 \cr
{\theta_\infty^2-\theta_1^2\over 4~\theta_\infty~r_1} &
{1\over r_1}
}~\pmatrix{ 1 & \log x \cr 0 & 1 }.
$$
}
 From the above result, together with (\ref{saravero?0}), we  compute $
\hat{A_0}=K_0\hat{\hat{A_0}}{K_0}^{-1}$, 
$\hat{A_1}=K_0\hat{\hat{A_1}}{K_0}^{-1}$. For example, 
{\small 
$$
\hat{A_0}
= G_0~
\pmatrix{
r+{\theta_0\over 2} + {\theta_0^2-\theta_x^2\over 4}\log x
~~&
{\theta_x^2-\theta_0^2\over 4}\log^2 x~ -2\left(
r+{\theta_0\over 2}
\right) \log x ~+{4~r(r+\theta_0)\over \theta_x^2 -\theta_0^2}
\cr
\cr
{\theta_0^2-\theta_x^2 \over 4} &
{\theta_x^2-\theta_0^2\over 4}\log x - \left(
r+{\theta_0\over 2}
\right)
}
~{G_0}^{-1}.
$$
}
A similar expression holds for $\hat{A_x}$. The reader can verify that
  the matching conditions (\ref{condition1}), (\ref{condition0})
   are satisfied. 

The leading terms of $y(x)$ are obtained from (\ref{leadingterm}) with matrix 
entries $
(\hat{A}_1)_{12}=-r_1$ and:  
$$
(\hat{A_0})_{12}=r_1~\left[ {\theta_x^2-\theta_0^2\over 4} \log^2
  x-2\left(r+{\theta_0\over 2}\right)\log x +{4~r(r+\theta_0)\over
    \theta_x^2-\theta_0^2}\right]. 
$$
The result is: 
\be
y(x)\sim x ~\left[
{\theta_x^2-\theta_0^2\over 4} \log^2 x-2\left(r+{\theta_0\over
  2}\right)\log x +{4~r(r+\theta_0)\over \theta_x^2-\theta_0^2} 
\right]
.
\label{STELLA}
\ee

%%%%%%%%%%%%%%%%%%%%%%%%%%%%%%%%%%%%%%%%%%%%%%%%%%

\subsection{Case  $\theta_0\pm\theta_x=0$}
\label{specialcasesec}

We consider here the cases (\ref{9}), (\ref{10}) 
of Proposition \ref{matrices} applied to the system (\ref{systemPhi0}). 

\vskip 0.3 cm
\noindent
{\it Case (\ref{9})}  is the case $
\sigma=0$,  $\theta_0=-\theta_x$, 
 with $a=0$, $c=\theta_0$ in the system (\ref{systemPhi0}). From Proposition 
\ref{matrices} we immediately have:
$$
\hat{\hat{A_0}}=\pmatrix{{\theta_0\over 2} & r \cr
                            0 & -{\theta_0\over 2} \cr
},~~~\hat{\hat{A_x}}=\pmatrix{{\theta_x\over 2} & 1-r \cr
                            0 & -{\theta_x\over 2} \cr
}. 
$$
The behavior of $\Psi_0$ and $\Psi_{OUT}$, and the matching are the
same of  subsection \ref{secsigma0}.  We
obtain the same $K_0(x)$. Therefore: 
$$
(\hat{A_0})_{12}= r_1~(r-\theta_0~\ln x),~~~(\hat{A_1})_{12}=-r_1.
$$
This gives the leading terms: 
\be
y(x)\sim x(r-\theta_0~\ln x)=x(r+\theta_x~\ln x).
\label{CECILIA}
\ee

In the same way, we treat the other cases. 
{\it Case (\ref{9})} with $a=c$, is the case $
\sigma=0$, $\theta_0=\theta_x$. As above, we find $
y(x)\sim x(r-\theta_0~\ln x)=x(r-\theta_x~\ln x)
$. {\it Case (\ref{10})} with $a=0$, is the case $
\sigma=0$, $\theta_0=-\theta_x$. We find $
y(x)\sim  x(r+\theta_0~\ln x)=x(r-\theta_x~\ln x)
$. {\it Case (\ref{10})} with $a=c$,  is the case 
$  
\sigma=0$, $\theta_0=\theta_x$. We find  $
y(x)\sim  x(r+\theta_0~\ln x)=x(r+\theta_x~\ln x)
$. 

%%%%%%%%%%%%%%%%%%%%%%%%
%%%%%%%%%%%%%
%%%%%%%%%%      CANCELED FROM THE PREPRINT
%%%%%%%%%%%%%
%%%%%%%%%%%%%%%%%%%%%%%%
%\vskip 0.3 cm
%\noindent
%{\it Case (\ref{9}) with $a=c$:} This is the case:
%$$
%\sigma=0,~~~\theta_0=\theta_x,
%$$
%and $a=c=\theta_0$ in system (\ref{systemPhi0}). 
%$$
%\hat{\hat{A_0}}=\pmatrix{{\theta_0\over 2} & r_0 \cr
%                            0 & -{\theta_0\over 2} \cr
%},~~~\hat{\hat{A_x}}=\pmatrix{-{\theta_x\over 2} & 1-r_0 \cr
%                            0 & {\theta_x\over 2} \cr
%}. 
%$$
%Proceeding as  above,
%$$
%\hat{y}(x)=  x(r_0-\theta_0~\ln x)=x(r_0-\theta_x~\ln x).
%$$
%
%
%\vskip 0.3 cm
%\noindent
%{\it Case (\ref{10}) with $a=0$:} This is the case 
%$$
%\sigma=0,~~~~\theta_0=-\theta_x
%$$
%and $a=0$, $c=\theta_0$  in system (\ref{systemPhi0}).  Therefore:
%$$
%\hat{\hat{A_0}}=\pmatrix{-{\theta_0\over 2} & r_0 \cr
%                            0 & {\theta_0\over 2} \cr
%},~~~\hat{\hat{A_x}}=\pmatrix{-{\theta_x\over 2} & 1-r_0 \cr
%                            0 & {\theta_x\over 2} \cr
%}. 
%$$
%The leading term is:
%$$
%\hat{y}(x)=  x(r_0+\theta_0~\ln x)=x(r_0-\theta_x~\ln x).
%$$
%
%
%\vskip 0.3 cm
%\noindent
%{\it Case (\ref{10}) with $a=c$:} This is the case 
%$$
%\sigma=0,~~~~\theta_0=\theta_x
%$$
%and $a=c=\theta_0$  in system (\ref{systemPhi0}).  Therefore:
%$$
%\hat{\hat{A_0}}=\pmatrix{-{\theta_0\over 2} & r_0 \cr
%                            0 & {\theta_0\over 2} \cr
%},~~~\hat{\hat{A_x}}=\pmatrix{{\theta_x\over 2} & 1-r_0 \cr
%                            0 & -{\theta_x\over 2} \cr
%}. 
%$$
%The leading term is:
%$$
%\hat{y}(x)=  x(r_0+\theta_0~\ln x)=x(r_0+\theta_x~\ln x).
%$$

\vskip 0.3 cm 
Both (\ref{STELLA}) and (\ref{CECILIA}) contain more than one term,
and in principle only the leading one is certainly correct. To
prove that they are all correct, we observe that  (\ref{STELLA}) and
(\ref{CECILIA}) 
can be obtained also by direct substitution of {\small $
y(x) = x (A_1 + B_1 \ln x + C_1 \ln^2 x + D_1 \ln^3 x + ...)+
x^2(A_2+B_2\ln x +...)+... $} 
into (PVI). We can recursively 
determine the coefficients by identifying 
 the same powers of $x$ and $\ln x$. As a result we obtain only the  
 five cases (\ref{fivecases}), which include  (\ref{STELLA}) and 
 (\ref{CECILIA}).

%%%%%%%%%%%%%%%%%%%%%%%%%%%%%%%%%%%%%%%%%%%%%%%%%%%%%

%%%%%%%%%%%%%%%%%%%%%%%%%%%%%%%%%%%%%%%%%%%%%%%%%

\subsection{No Naive Matching for $\sigma=1$}
\label{naive} 

The  condition $|\Re \sigma| <1 $ suggests that the    
matching above does not work 
in the case $\sigma=1$ (and $\sigma=-1$, being equivalent).
Let 
us convince ourselves of this fact by repeating the procedure above. 
A fundamental matrix solution for (\ref{system1}) at $\lambda=0$ is non-generic:
\be
\Psi_{OUT}(\lambda)= \left( G_0~+ ~O(\lambda) \right) ~\pmatrix{\lambda^{1\over 2} & 0 
\cr
0 & \lambda^{-{1\over 2}}
}
~
\pmatrix{ 1 & \log \lambda 
\cr
                                           0 & 1 }.
\label{philog1}
\ee
where:
$$
G_0=\pmatrix{ 1 & {4\over \theta_1^2-(\theta_\infty -1)^2} \cr
\cr
{(\theta_\infty+1)^2-\theta_1^2 \over 4 \theta_\infty ~r_1} 
&
-{1\over \theta_\infty ~r_1}
},~~~~~ \forall r_1\neq 0. 
$$
A fundamental matrix solution of (\ref{system0}) at $\mu=\infty$ is non-generic:
\be
\Psi_0(\mu)= \left(I+O\left({1\over \mu}\right)\right)~
\pmatrix{ \mu^{1\over 2} & 0 \cr
0 & \mu^{-{1\over 2}}}
~
\pmatrix{1 & 0 \cr
R~\log \mu & 1 },
\label{philog0}
\ee
where:
$$
R:=(\hat{\hat{A}_x})_{21}={[(\theta_0+\theta_x)^2-1][(\theta_0-\theta_x)^2-1]
\over 16 ~r},~~~~~ r\neq 0.
$$
The matching relation:
{\small
$$
K_0(x)~ \pmatrix{ \left({\lambda \over x}\right)^{1\over 2} 
& 0 \cr
R~ \left({\lambda \over x}\right)^{-{1\over 2}}\log \left({\lambda \over x}\right) & \left({\lambda \over x}\right)^{-{1\over 2}}
}
~\sim ~
\pmatrix{ 1 & {4\over \theta_1^2-(\theta_\infty -1)^2} 
\cr
\cr
{(\theta_\infty+1)^2-\theta_1^2 \over 4 \theta_\infty ~r_1} 
&
-{1\over \theta_\infty ~r_1}
} ~ 
\pmatrix{ \lambda^{1\over 2} &  \lambda^{1\over 2}\log \lambda
\cr
         0 & \lambda^{-{1\over 2} }
},
$$
}
\noindent
 shows that we cannot eliminate $\lambda$ to obtain $K_0(x)$. 

One  case $\sigma=1$ is studied in Part III, making use of a 
non-fuchsian reductions of the system 
(\ref{SYSTEM}).

\subsection{Monodromy Data}
\label{senzaconti}
Systems (\ref{system1})  
(\ref{system0}) are equivalent to Gauss
hyper-geometric equations, as it is explained in Appendix 1 (make use of
the systems (\ref{systemPhi1}) and (\ref{systemPhi0})
respectively). Therefore, the monodromy can be computed in a standard
way, using the connection formulae for the hyper-geometric functions.

 We
obtain in this way the monodromy of $\Psi_{OUT}$ and
$\Psi_{IN}$.  As it is
explained in section \ref{MonodromyPasqua}, 
it may be necessary to do a transformation
$\Psi_{OUT}\mapsto \Psi_{OUT}^{Match}:= \Psi_{OUT}C_{OUT}$, 
in order to match the ``out'' and ``in'' solutions with
a solution $\Psi$ of (\ref{SYSTEM}). In this way, the monodromy
matrices 
$M_0$, $M_x$, $M_1$ of
$\Psi$ can be obtained. They depend on $r$. We then compute  the traces of
$M_iM_j$ and extract $r$, which is thus obtained as a function of the
monodromy data. 

We do not repeat the computations here. One example is the computation
of (\ref{bfa}) in  \cite{Jimbo} and   \cite{DM} \cite{guz1}
\cite{guz3} \cite{Boalch}.

\vskip 1 cm 
\centerline{\Large \bf PART III -- Derivations for   
 Non-fuchsian Reduction}

\section{ Case $\theta_1=\pm\theta_\infty$, $\theta_x=\pm \theta_0$, 
$\lim_{x\to 0}(A_x+A_0)=0$. Solution (\ref{form2})}  
\label{irregular0}

 We begin by observing that for $\theta_x=\pm \theta_0$, system 
(\ref{systemPhi0}) may fall in  cases  (\ref{6}) and (\ref{7}). 
 If it is so, then 
$\hat{\hat{A_0}}+\hat{\hat{A_x}}=0$, and therefore  $
\hat{A_0}+\hat{A_x}=0$. More precisely, we start from  the following hypotheses:
$$
\lim_{x\to 0} \bigl(A_0(x)+A_x(x)\bigr)=0,
$$
$$
A:=\lim_{x\to 0} A_x(x) = \hbox{ a constant matrix with eigenvalues
}\pm{\theta_x\over 2}
$$
The first hypothesis means that we can write (the trace is zero): 
$$
A_0+A_x= \pmatrix{ a(x) & b(x)~r \cr
                               c(x)~{1\over r} & -a(x) },~~~\lim_{x\to 0}a(x)=
\lim_{x\to 0}b(x)=\lim_{x\to 0}c(x)=0,
$$
The second hypotheses implies that the general form of  $A$ is:
$$
A=\pmatrix{ s+{\theta_x\over 2} & -r 
\cr
                   (s+\theta_x)~s\over r &  -s-{\theta_x\over 2}
},~~~r,s\in{\bf C},~~r\neq 0.
$$  
We also write:
$$
A_x(x)-A=:\Delta_x(x),~~~A_0+A=:\Delta_0(x),~~~~~\Delta_0+\Delta_x\equiv A_0+A_x.
$$
$\Delta_x(x)$ and $\Delta_0(x)$ are vanishing. 
We suppose that the slowest vanishing behavior be of order 
$x^{\sigma_0}$, for some $\sigma_0>0$. Namely:
$$
a(x),~b(x),~c(x),~
(\Delta_x)_{ij}(x),~(\Delta_0)_{ij}(x)~=O(x^{\sigma_0}), ~~~\sigma_0>0.
$$
Finally, we have:
$$
A_1(x)= -{\theta_{\infty}\over 2}~\sigma_3~-(A_0+A_x)
~\longrightarrow~ -{\theta_{\infty}\over 2}~\sigma_3,~~~x\to 0
.
$$
\subsection{Coalescence of Singularities}

{\bf 1) THE SYSTEM for $\Psi_{OUT}$.} 
\vskip 0.2 cm
We consider system (\ref{nonfuchsianSYSTEMOUT}), in the domain
$ |\lambda|\geq|x|^{\delta_{OUT}}$. 
Let us determine the conditions to neglect a term $x^n
A_x/\lambda^{n+1}$ -- and all the terms following it -- 
with respect to $(A_0+A_x)/\lambda$, when $x\to 0$, 
$\lambda\sim x^{\delta}$, $\delta\leq \delta_{OUT}$.
$$
\hbox{ We can neglect }~ {x^n A_x\over \lambda^{n+1}} 
\Longleftrightarrow
\left|
{x^n\over \lambda^n } A_x 
\right| \ll \bigl| (A_0+A_x)_{ij} \bigr|,~~~\forall i,j\in\{1,2\}.
$$
 Since $\lim_{x\to 0} (A_x)_{ij}$ are non-zero constants, 
the above condition is: $|x|^{n-n\delta}\ll |x|^{\sigma_0}$, namely:
$\delta<1-\sigma_0/n$. We state this result as a lemma. 
 
\ble
 Let $N_{OUT}\geq 2$ be an integer.  We can approximate 
(\ref{nonfuchsianSYSTEMOUT}) with:
$$
{d\Psi_{OUT}\over d \lambda}=
\left[
{(A_0+A_x)\over \lambda}+{A_x\over \lambda}~\sum_{n=1}^{N_{OUT}-1}\left({x\over \lambda}\right)^n+ {A_1\over \lambda-1}
\right]~\Psi_{OUT}.
$$
 if and only if: 
\be
\delta_{OUT}<1-{\sigma_0\over N_{OUT}}.
\label{condout.april}
\ee
\ele

Suppose that all term $x^n A_x / \lambda^{n+1}$, with $n\geq
  N_{OUT}$,  have been neglected. We can also make the substitution $
A_1\mapsto -{\theta_{\infty}\over 2} \sigma_3$ 
in ${A_1\over \lambda-1}$, if and only if the error term 
$-{A_0+A_x\over \lambda-1}$,
 is smaller than   
${ x^{N_{OUT}-1} A_x \over \lambda^{N_{OUT}} }$. Namely, if and only if:
$$
\bigl|(A_0+A_x)_{ij}  \bigr|\ll \left|
{ x^{N_{OUT}-1} (A_x)_{ij} \over \lambda^{N_{OUT}} }\right|. 
$$
This is $|x|^{\sigma_0}\ll |x|^{N_{OUT}-1-N_{OUT}\delta_{OUT}}$,
namely $
\delta_{OUT} > 1-{1+\sigma_0\over N_{OUT}}$. 

\vskip 0.2 cm 
We can also do the substitution $A_x\mapsto A$, provided that
$\delta_{OUT}<1$. This is because we can 
neglecting terms $ {x^n\Delta_x\over \lambda ^{n+1}}$ with respect
to ${A_0+A_x\over \lambda}$, where both $\Delta_x$ and $A_0+A_x$ are
$O(x^{\sigma_0})$,  $\lambda\sim x^{\delta}$, $\delta<1$.  
We summarize the result in the following lemma.

\ble 
 Let $N_{OUT}\geq 2$ be  integer.  We can approximate
 (\ref{nonfuchsianSYSTEMOUT}) with:
$$
{d\Psi_{OUT}\over d \lambda}=
\left[
{A_0+A_x\over \lambda}+{A\over
  \lambda}~\sum_{n=1}^{N_{OUT}-1}\left({x\over
  \lambda}\right)^n-{\theta_\infty\over 2} {\sigma_3\over \lambda-1}
\right]~\Psi_{OUT}.
$$
 if and only if: 
$$
                     1-{1+\sigma_0\over N_{OUT}}
                     <\delta_{OUT}<1-{\sigma_0\over N_{OUT}}.
$$
In particular, this means that $\delta_{OUT}<1$. 
\ele

\vskip 0.2 cm
\noindent
{\bf Example:}
\vskip 0.2 cm

If $\sigma_0=1$ and $N_{IN}=2$ we have:
$$
{d\Psi_{OUT}\over d \lambda}=
\left[{x~A\over \lambda^2}+
{A_0+A_x\over \lambda}
  -{\theta_\infty\over 2} {\sigma_3\over \lambda-1}
\right]~\Psi_{OUT},
~~~~~~~
                    0 <\delta_{OUT}<{1\over 2}.
$$

\vskip 0.2 cm
If $\sigma_0=1$ and $N_{IN}=3$ we have:
$$
{d\Psi_{OUT}\over d \lambda}=
\left[{x^2 A\over \lambda^3}+{x~A\over \lambda^2}+
{A_0+A_x\over \lambda}
  -{\theta_\infty\over 2} {\sigma_3\over \lambda-1}
\right]~\Psi_{OUT},
~~~~~~~
                    {1\over 3} <\delta_{OUT}<{2\over 3}.
$$

\vskip 0.3 cm
\noindent
{\bf 2) THE SYSTEM for $\Psi_{IN}$.}

\vskip 0.2 cm

 We consider system (\ref{nonfuchsianSYSTEMIN}) in the domain
 $|\lambda|\leq |x|^{\delta_{IN}}$. 
We investigate the condition necessary and sufficient to neglect a
term 
$\lambda^n A_1$ (and all its next terms) with respect to ${A_0\over
  \lambda}+{A_x\over \lambda-x}$. It is convenient to write:
$$
{A_0\over
  \lambda}+{A_x\over \lambda-x}= {A_0+A_x\over \lambda-x} - {xA_0\over 
\lambda(\lambda-x)}.
$$
Suppose that $\lambda\sim x^{\delta}$, $\delta\geq\delta_{IN}$.
$$
\hbox{ We neglect }~~ A_1\lambda^n ~~
\Longleftrightarrow 
\left\{
\matrix{
\left|
{xA_0\over 
\lambda(\lambda-x)}
\right|
\gg 
\bigl|
A_1\lambda^n
\bigr|,~\hbox{ namely: }~ |x|^{1-2\delta}\gg
|x|^{n\delta}~~\Leftrightarrow~~\delta>{1\over n+2};
\cr
\cr
\left|
 {A_0+A_x\over \lambda-x}
\right|
\gg
\bigl|
A_1\lambda^n
\bigr|
,~\hbox{ namely: }~ |x|^{\sigma_0-\delta}\gg
|x|^{n\delta}~~ \Leftrightarrow~~ \delta>{\sigma_0\over n+1}.
}
\right. 
$$
Thus, we have the condition $\delta>\max \left\{ {\sigma_0\over
  n+1},~{1\over n+2}\right\}$. We have proven   the following:

\ble
Let $N_{IN}\geq 1$ be an integer. We approximate
(\ref{nonfuchsianSYSTEMIN}) with:
 $$
{d \Psi_{IN} \over d\lambda}
=
\left[
{A_0\over \lambda} + {A_x \over \lambda -x} - A_1 \sum_{n=0}^{N_{IN}-1} 
\lambda^n
\right]~
\Psi_{IN}.
$$
if and only if:
$$
\delta_{IN}>\max \left\{ {\sigma_0\over N_{IN}+1},~{1\over N_{IN}+2}\right\}.
$$
\ele

We further make the substitution $A_1\mapsto
-{\theta_\infty\over 2}\sigma_3$.  This is possible if and only if  
 two conditions are true: 
1) $\left|
{xA_0\over \lambda(\lambda-x)}
\right|$ and  $ \left|
{A_0+A_x\over \lambda-x}
\right|  $ are dominant w.r.t. the term $A_0+A_x$ appearing in  $A_1=-{\theta_\infty\over
  2}\sigma_3-(A_0+A_x)$. ~~2) $\bigl|\lambda^{N_{IN}-1}A_1\bigr|$, i.e. $\bigl|\lambda^{N_{IN}-1}\sigma_3\bigr|$, is 
 dominant w.r.t. the term $A_0+A_x$ in $A_1$.  Esplicitely, the conditions are: 
$$ 
\left|
{A_0+A_x\over \lambda-x}
\right|\gg \bigl|
A_0+A_x
\bigr|
~\hbox{ (this is always true)},
$$
$$
\left|
{xA_0\over \lambda(\lambda-x)}
\right|
\gg 
\bigl|
A_0+A_x
\bigr|
~~
\Longleftrightarrow
~~
|x|^{1-2\delta}>x^{\sigma_0},\hbox{ namely: } \delta>{1-\sigma_0\over
  2},
$$
and
$$
\bigl| \lambda^{N_{IN}-1}\sigma_3\bigr|
\gg 
\bigl|
A_0+A_x
\bigr|~~\Longleftrightarrow~~ |x|^{(N_{IN}-1)\delta}>|x|^{\sigma_0}, \hbox{
  namely: }\delta<{\sigma_0\over N_{IN}-1}. 
$$
We have:

\ble
Let $N_{IN}\geq 1$ be an integer. We approximate
(\ref{nonfuchsianSYSTEMIN}) with:
 $$
{d \Psi_{IN} \over d\lambda}
=
\left[
{A_0\over \lambda} + {A_x \over \lambda -x} +{\theta_\infty\over
  2}\sigma_3
 \sum_{n=0}^{N_{IN}-1} 
\lambda^n
\right]~
\Psi_{IN}
$$
$$
=
\left[
{A_0+A_x\over \lambda-x} -{xA_0\over \lambda(\lambda-x)} +{\theta_\infty\over
  2}\sigma_3
 \sum_{n=0}^{N_{IN}-1} 
\lambda^n
\right]~
\Psi_{IN},
$$
if and only if:
$$
\max \left\{{1-\sigma_0\over 2},  {\sigma_0\over N_{IN}+1},~{1\over
  N_{IN}+2}\right\}<\delta_{IN}<{\sigma_0\over N_{IN}-1}.
$$
\ele

\vskip 0.3 cm
As a final simplification, we substitute 
$A_0=-A+\Delta_0\mapsto -A$. 
 This is possible if and only if:
 $$
\left|{x~\Delta_0\over
  \lambda(\lambda-x)}\right|
\ll \left|{A_0+A_x\over \lambda -x}\right|~~\Longleftrightarrow ~~
|x|^{1+\sigma_0-2\delta}<|x|^{\sigma_0-\delta}, \hbox{ namely }
\delta<1, 
$$  
and:
$$
\left|{x~\Delta_0\over
  \lambda(\lambda-x)}\right|
\ll \bigl|
\lambda^{N_{IN}-1} {\theta_\infty\over 2}~\sigma_3
\bigr|
~~
~~\Longleftrightarrow ~~ |x|^{1+\sigma_0-2\delta}\ll x^{(N_{IN}-1)\delta}
, \hbox{ namely: } \delta<{\sigma_0+1\over N_{IN}+1}. 
$$
We have proven the following:

\ble
Let $N_{IN}\geq 1$ be an integer. We can approximate  
(\ref{nonfuchsianSYSTEMIN}) 
 with:
 $$
{d \Psi_{IN} \over d\lambda}
=
\left[
{A_0+A_x\over \lambda-x} +{xA\over \lambda(\lambda-x)} +{\theta_\infty\over
  2}\sigma_3
 \sum_{n=0}^{N_{IN}-1} 
\lambda^n
\right]~
\Psi_{IN},
$$
if and only if:
$$
\max \left\{{1-\sigma_0\over 2},  {\sigma_0\over N_{IN}+1},~{1\over
  N_{IN}+2}\right\}<\delta_{IN}<\min\left\{{\sigma_0\over N_{IN}-1}, {\sigma_0+1\over N_{IN}+1}\right\}.
$$
\ele

\vskip 0.2 cm
\noindent
{\bf Examples:}
\vskip 0.2 cm  
If $\sigma_0=1$ and $N_{IN}=1$, we have:
 $$
{d \Psi_{IN} \over d\lambda}
=
\left[
{A_0+A_x\over \lambda-x} +{xA\over \lambda(\lambda-x)} +{\theta_\infty\over
  2}\sigma_3
\right]~
\Psi_{IN},
~~~~~~~
 {1 \over 2} <\delta_{IN} <1.
$$
If we keep $-A_0$ instead of  $A$, with no change in the condition on
$\delta_{IN}$, we can also rewrite:
$$
{d \Psi_{IN} \over d\lambda}
=
\left[
{A_0\over \lambda} +{A_x\over \lambda-x} +{\theta_\infty\over
  2}\sigma_3
\right]~
\Psi_{IN},
$$

\vskip 0.2 cm
If $\sigma_0=1$ and $N_{IN}=2$, we have:
 $$
{d \Psi_{IN} \over d\lambda}
=
\left[
{A_0+A_x\over \lambda-x} +{xA\over \lambda(\lambda-x)} +{\theta_\infty\over
  2}\sigma_3(1+\lambda) \right]~
\Psi_{IN},
~~~~~~~
{1\over 3} <\delta_{IN}<{2\over 3}.
$$
Equivalently, we can write:
$$
{d \Psi_{IN} \over d\lambda}
=
\left[
{A_0\over \lambda} +{A_x\over \lambda-x} +{\theta_\infty\over
  2}\sigma_3(1+\lambda) \right]~
\Psi_{IN}.
$$

\subsection{Matching}

We do the matching in the overlapping region
 $|x|^{\delta_{OUT}}\leq|\lambda|\leq
|x|^{\delta_{IN}}$. This imposes: $\delta_{IN}\leq \delta_{OUT}$. 
In order for the overlapping region not to be empty, we 
 must choose  suitable reductions of (\ref{nonfuchsianSYSTEMOUT}) and
(\ref{nonfuchsianSYSTEMIN}). 
If we expect $\sigma_0$ to be close to 1, we try to match solutions
$\Psi_{OUT}$ and $\Psi_{IN}$ satisfying one of the following sets of
systems:

\vskip 0.2 cm
\noindent
{\it First choice:}
$$
{d\Psi_{OUT}\over d \lambda}= \left[
{xA \over \lambda^2} +{A_0+A_x\over \lambda} -{\theta_\infty\over 2}{\sigma_3\over \lambda-1} 
\right] ~\Psi_{OUT},
$$
$$
{d \Psi_{IN} \over d\lambda}
=
\left[
{A_0\over \lambda} +{A_x\over \lambda-x} +{\theta_\infty\over
  2}\sigma_3(1+\lambda) \right]~
\Psi_{IN}.
$$
The condition to be satisfied for $\sigma_0\cong 1$ is: ${\sigma_0\over
  3}<\delta_{IN}\leq\delta_{OUT} <1-{\sigma_0\over 2}$. For $\sigma_0=1$,
this is:
$$
{1\over 3}<\delta_{IN}\leq\delta_{OUT}<{1\over 2}.
$$
\vskip 0.2 cm 
\noindent
{\it Second choice:}
$$
{d\Psi_{OUT}\over d \lambda}= \left[{x^2A\over \lambda^3}+
{xA \over \lambda^2} +{A_0+A_x\over \lambda} -{\theta_\infty\over 2}{\sigma_3\over \lambda-1} 
\right] ~\Psi_{OUT},
$$
$$
{d \Psi_{IN} \over d\lambda}
=
\left[
{A_0\over \lambda} +{A_x\over\lambda-x} +{\theta_\infty\over
  2}\sigma_3 \right]~
\Psi_{IN}.
$$
 For $\sigma_0\cong 1$, the condition to be satisfied is:
 ${\sigma_0\over 2}<\delta_{IN}\leq \delta_{OUT} < 1-{\sigma_0\over
 3}$. For $\sigma_0=1$,  this is:
$$
{1\over 2}<\delta_{IN}\leq\delta_{OUT}<{2\over 3}.
$$

\vskip 0.2 cm
In both cases, the overlapping regions are not empty. The matching
procedure will determine the leading terms (order $x^{\sigma_0}$) 
 of  the unknown matrix
elements  
 $a(x)$, $b(x)$, $c(x)$ of $A_0+A_x$.

\subsection{Matching for the  First Choice: ${1\over 3}<
  \delta_{IN}\leq \delta_{OUT}<{1\over 2}$.} 

We rewrite the systems in a more convenient form:
$$
\nu:={1\over \lambda},~~~\mu:={\lambda\over x};
$$
\be
{d\Psi_{OUT}\over d \nu}= \left[
-xA -{A_0+A_x\over \nu} -{\theta_\infty\over 2}{\sigma_3\over \nu(\nu-1)} 
\right] ~\Psi_{OUT}~~~~~~~~~~~~~~~~~~~~~
\label{tildeoutapril}
\ee
\vskip 0.2 cm
$$
{d \Psi_{IN} \over d\mu}
=
\left[x^2~{\theta_\infty\over 2}\sigma_3~\mu~+x~{\theta_\infty\over
    2}\sigma_3     
+{A_0\over \mu} +{A_x\over \mu-1} \right]~
\Psi_{IN}~~~~~~~~~~~~~~~~~~~~~
$$
$$
=
\left[x^2~{\theta_\infty\over 2}\sigma_3~\mu~+x~{\theta_\infty\over 2}\sigma_3   +{A_0+A_x\over \mu} -{A_0\over \mu(\mu-1)} \right]~
\Psi_{IN}.
$$
$$
\hbox{ Then we substitute $A_0\mapsto -A$ in the last term}.
$$
In the matching region $|x|^{-\delta_{IN}}\leq |\nu|\leq
|x|^{-\delta_{OUT}}$, $|x|^{\delta_{OUT}-1}\leq\mu\leq
|x|^{\delta_{IN}-1}$, we have 
$\nu\to \infty$, $
\mu\to\infty$. The point at infinity is a non-fuchsian singularity. 
%%%%%%%%%%%%
\footnote{ System (\ref{tildeoutapril}) can be also written with 
${\theta_\infty\over 2}\sigma_3\mapsto -A_1$:
$$
{d\Psi_{OUT}\over d \nu}= \left[
-xA -{A_0+A_x\over \nu} +{A_1\over \nu(\nu-1)} 
\right] ~\Psi_{OUT}
$$
$$
\left[
-xA -{A_0+A_x+A_1\over \nu} +{A_1\over \nu-1} 
\right] ~\Psi_{OUT}
$$
After diagonalization, we get:
$$
{d\tilde{\Psi}_{OUT}\over d \nu}= \left[
x{\theta_x\over 2} \sigma_3 +{\theta_\infty\over 2}{G^{-1}\sigma_3~G\over \nu} +{G^{-1}A_1~G\over \nu-1} 
\right] ~\tilde{\Psi}_{OUT}.
$$
This form is that of a system of isomonodromy deformation for the
fifth Painlev\'e equation. 
}
%%%%%%%%%%%%%%%%%
 In order to find the local behavior at this point, it is convenient
 to put the leading term in diagonal form. 
 Let $G$ be the invertible matrix such that:
$$
 G^{-1}~ A G = -{\theta_x\over 2}\sigma_3,~~~\hbox{ for example: }~G=\pmatrix{1 & {r\over
 s} \cr {s+\theta_x\over r} & 1},
$$
and put:
$$
\Psi_{OUT}=:G~\tilde{\Psi}_{OUT}.
$$
Then:
\be
{d\tilde{\Psi}_{OUT}\over d\nu}=
\left[
x~{\theta_x\over 2}\sigma_3 -{G^{-1}(A_0+A_x)G\over \nu}
-{\theta_\infty\over 2}~G^{-1}\sigma_3 G~\left({1\over \nu^2}+{1\over \nu^3}+...\right)
\right]
~\tilde{\Psi}_{OUT},~~~\nu\to\infty;
\label{irroutApril}
\ee
\be
{d \Psi_{IN} \over d\mu}
=
\left[x^2~{\theta_\infty\over 2}\sigma_3~\mu~+x~{\theta_\infty\over
    2}\sigma_3   +{A_0+A_x\over \mu} +A ~\left({1\over \mu^2}+{1\over \mu^3}+...
\right) \right]~
\Psi_{IN},~~~\mu\to \infty.
\label{irrinApril}
\ee

In order to write the local behavior of $\tilde{\Psi}_{OUT}$ and
  $\Psi_{IN}$ at infinity, we observe that the systems
  (\ref{irroutApril}) and (\ref{irrinApril}) respectively 
have the following forms:
\be
{dY_1\over dz}= \left[\Omega+{D_1\over z}+{D_2\over z^2}+{D_3\over z^3}+...\right]~Y_1,
\label{IRR1}
\ee
\be
{dY_2\over dz}= \left[x^2\Lambda~z+x\Lambda~z+{E_1\over z}+{E_2\over z^2}+{E_3\over z^3}+...\right]~Y_2,
\label{IRR2}
\ee
where $\Omega$ and $\Lambda$ are diagonal matrices with {\it distinct
  eigenvalues}. In our case:
$$
\Omega= x~ {\theta_x\over 2} \sigma_3,~~~~~\Lambda={\theta_\infty\over
  2}\sigma_3.
$$
The eigenvalues are distinct iff $\theta_x\neq 0$, $\theta_\infty\neq
0$. 
\vskip 0.2 cm 
The theory for such systems is developed in \cite{BJL1} (see also
\cite{BJL2}).  
 For any sector of angular width $\pi+\epsilon$,
$\epsilon>0$ sufficiently small, there exists a unique solution of (\ref{IRR1}) with  asymptotic
expansion:
$$
Y_1(z)\sim\left[I+{G_1\over z} +{G_2\over z^2}+...\right] ~\exp\bigl\{
\Omega~z 
\bigr\} ~z^{\Omega_1},~~~z\to \infty.
$$
$$
\Omega_1=\hbox{ diagonal part of } D_1.
$$
 For any sector of angular width ${\pi\over 2} +\epsilon$,
$\epsilon>0$ sufficiently small, there exist a unique solution of (\ref{IRR2}) with asymptotic
expansion: 
$$
Y_2(z)\sim\left[I+{K_1\over z} +{K_2\over z^2}+...\right] ~\exp\left\{
{x^2\over 2}\Lambda~z^2 ~+ x\Lambda~z
\right\} ~z^{\Lambda_1},~~~z\to \infty.
$$
$$
\Lambda_1=\hbox{ diagonal part of } E_1.
$$

We can always find two solutions $Y_1(z)$ and $Y_2(z)$ as above, such
that the sectors where the asymptotic expansions hold are
overlapping. We refer the reader to  \cite{BJL1} for the general
description of irregular system with a Stokes phenomenon, and to the
Appendix 2 for the computation of the matrices $G_i$, $K_i$, $i=1,2,...$

\vskip 0.3 cm

The systems (\ref{irroutApril}), (\ref{irrinApril}) are {\it
  isomonodromic}. This imposes that $\Omega_1$ and $\Lambda_1$ must be
  independent of $x$. They are:
$$
\Omega_1=\hbox{ diagonal of } \bigl(-G^{-1}(A_0+A_x)~G\bigr),
$$
$$
\Lambda_1=\hbox{ diagonal of } (A_0+A_x)~=\pmatrix{ a(x) & 0 \cr 0 & -a(x)} .
$$
We compute:
$$
G^{-1}(A_0+A_x)~G= {1\over \theta_x}~
\pmatrix{
-(2s+\theta_x)a-s(s+\theta_x)b+c & \left(-2a-sb+{c\over s}\right)~r
\cr
\cr
{s\over r} ~\bigl(
2(s+\theta_x)a+(s+\theta_x)^2b-c
\bigr)
&
(2s+\theta_x)a+s(s+\theta_x)b-c
}
$$
Since $a,b,c$ vanish, the condition of isomonodromicity implies that: 
$$
\Omega_1=0,~~~~~\Lambda_1=0
$$
This means that {\it the leading terms of} $a(x)$, $b(x)$, $c(x)$
satisfy the conditions:
$$
a(x)=0,~~~~~c(x)=s(s+\theta_x)b(x).
$$
The above conditions mean that if $b(x),~ c(x) = O(x^{\sigma_0})$, then 
$a(x)$ is of higher order, i.e. it vanishes faster than
$x^{\sigma_0}$. Note that with this choice of $a$, $b$, $c$ we get:
$$
G^{-1}(A_0+A_x)~G=\pmatrix{0 & br 
\cr
\cr
{s(s+\theta_x)\over r}~b & 0
}
$$ 

\vskip 0.3 cm
We are ready to write the behavior of $\Psi_{OUT}$:
$$
\Psi_{OUT}= G~\left[
I+{G_1\over \nu}+{G_2\over \nu^2}+...
\right]~\exp\left\{
x {\theta_x\over 2}\sigma_3\nu
\right\},~~~~~\nu\to\infty
$$
$$
=
$$
$$
=G~\bigl[
I+G_1\lambda+G_2\lambda^2+...
\bigr]~\exp\left\{
{\theta_x\over 2} \sigma_3~{x\over \lambda} 
\right\},~~~~~\lambda\to 0.
$$
We use the formulae of  Appendix 2 to determine $G_1$:  
$$
(G_1)_{ij}= 2~{\bigl[G^{-1}(A_0+A_x)~G\bigr]_{ij}\over 
x ~\theta_x \bigl( (\sigma_3)_{ii}-(\sigma_3)_{jj} \bigr)
},~~~~~i\neq j.
$$
\vskip 0.2 cm
$$
(G_1)_{ii}= {\theta_\infty\over 2}\bigl(G^{-1}\sigma_3~ G \bigr)_{ij}
+
2~{
\bigl[G^{-1}(A_0+A_x)~G\bigr]_{ij}\bigl[G^{-1}(A_0+A_x)~G\bigr]_{ji}
\over 
x ~\theta_x \bigl( (\sigma_3)_{jj}-(\sigma_3)_{ii} \bigr)
}.
$$
In the second term of the last formula $j=2$ if $i=1$,  $j=1$ if
$i=2$. We compute:
\be
G^{-1}\sigma_3~ G
={1\over \theta_x}~ 
 \pmatrix{
-(2s+\theta_x) & -2 r 
\cr
\cr
{2s(s+\theta_x)\over r} & 2s+\theta_x
}. 
\label{SIGMA3}
\ee
Therefore:
$$
(G_1)_{12}={r\over x~\theta_x}~b,~~~~~(G_1)_{21}=
-{s(s+\theta_x)\over x~\theta_x~r}~b,
$$
\be
(G_1)_{11}= -{\theta_\infty\over 2\theta_x}(2s+\theta_x) ~+
  {s(s+\theta_x)\over x~\theta_x}~b^2,~~~~~(G_1)_{22}=-(G_1)_{11}
\label{G1OUT}
\ee

\vskip 0.3 cm
\noindent 
On the other hand, the local behavior of $\Psi_{IN}$ is:
$$
\Psi_{IN}= \left[I + {K_1\over \mu} +{K_2\over \mu^2}+...\right]~\exp\left\{
x^2{\theta_\infty\over 4}\sigma_3~{\mu^2}~+x{\theta_\infty\over 2}\sigma_3~\mu
\right\},~~~~~\mu\to\infty
$$
$$
 =\left[I + K_1~{x\over \mu} +K_2~{x^2\over \lambda^2}+...\right]~
~\exp\left\{
{\theta_\infty\over 4}\sigma_3~{\lambda^2}~+{\theta_\infty\over 2}\sigma_3~\lambda
\right\},~~~~~\lambda\to 0.
$$
We determine $K_1$ from the formulas of Appendix 2. 
$$
K_1= \hbox{ diagonal part of } (-A)~= \pmatrix{-\left(s+{\theta_x\over 2}\right) & 0
  \cr\cr
0 & s+{\theta_x\over 2}
}
$$ 

\vskip 0.2 cm
The matching conditition:  
$$
\Psi_{OUT}(\lambda,x)\sim\Psi_{IN}(\lambda, x),~~~~x\to 0,~~~|x|^{\delta_{OUT}}\leq |\lambda|\leq |x|^{\delta_{IN}}.
$$
 is
restricted to the overlapping sector where both expansions of
$\Psi_{OUT}$ and $\Psi_{IN}$ hold. Noting that $\Psi_{OUT}\sim G$ and
 $\Psi_{IN}\sim I$, we choose the 
new solution
$\Psi_{OUT}\mapsto\Psi_{OUT}~ G^{-1}$. Then, we expand the exponents:
$$
\Psi_{OUT}= ~\bigl[
I+GG_1G^{-1}\lambda+GG_2G^{-1}\lambda^2+...
\bigr]~\left[
I+{\theta_x\over 2} G\sigma_3G^{-1}{x\over \lambda}+{{\theta_x}^2\over
  8}~{x^2\over \lambda^2}~+...
\right]. 
$$
%%%%%%%%%%%%%%%%%%%%%%%
%%%%%%%%%%%%%%
%%%%%%%          CANCELED IN THE PREPRINT 
%%%%%%%%%%%%%%
%%%%%%%%%%%%%%%%%%%%%%%
%$$
%=I +{\theta_x\over 2} G~\sigma_3~G^{-1}{x\over \lambda}+{{\theta_x}^2\over
%  8}~{x^2\over \lambda^2}\right)~+...
%$$
%$$
%+G~G_1~G^{-1}\lambda~+G~G_1~G^{-1}\left[
%{\theta_x\over 2}G\sigma_3G^{-1}~x~+{{\theta_x}^2\over 8}~{x^2\over
%  \lambda}+
%{{\theta_x}^3\over 48}G\sigma_3G^{-1}~{x^3\over\lambda^2}+...
%\right]
%$$
%$$
%+G~G_2~G^{-1}\lambda^2~ +G~G_2~G^{-1}\left[
%{\theta_x\over 2}G\sigma_3G^{-1}~x\lambda~+{{\theta_x}^2\over 8}~x^2+
%{{\theta_x}^3\over 48}G\sigma_3G^{-1}~{x^3\over\lambda}+...
%\right]~+...~.
%$$
The point here is quite delicate. We consider the relation of
dominance among terms -- and write the leading terms of the expansion
--
  as they are  
 in case the $G_n(x)$'s
are not divergent when $x\to 0$.  
Keeping into account that ${\theta_x\over 2}
G\sigma_3G^{-1}=-A$,  the dominant terms are: 
$$
\Psi_{OUT}(\lambda,x)=
I+GG_1G^{-1}\lambda~-A~{x\over\lambda}~+O\left(
\lambda^2,{x^2\over \lambda^2}, x
\right)
$$
It is important to note that $\lambda$ is dominant w.r.t ${x\over
  \lambda}$, because $\delta_{OUT}<{1\over 2}$; namely, $\lambda\sim
  x^{\delta}$ vanishes slower than ${x\over
  \lambda}\sim x^{1-\delta}$, as $x\to 0$.

\vskip 0.2 cm 
We expand the exponent in $\Psi_{IN}$, and keep only the first
dominant terms
(in the spirit of the observation on the dominance relations made above):
$$
\Psi_{IN}(\lambda,x)=I + {\theta_\infty\over 2} \sigma_3~\lambda~+K_1~{x\over
  \lambda}
+O\left(
\lambda^2,{x^2\over \lambda^2}, x
\right)
$$

$\Psi_{OUT}$ and $\Psi_{IN}$ match in the first term $I$. We impose
the matching of the second term, namely the term in $\lambda$:
$$
G~G_1(x)~G^{-1} ~\sim~ {\theta_\infty\over 2}\sigma_3,~~~~~x\to 0.
$$
Namely:
\be
G_1(x)~\sim~ {\theta_\infty\over 2} ~G^{-1}\sigma_3~G,~~~~~x\to 0
\label{matchcond1choice}
\ee
From the explicit form of $G_1$ and $
G^{-1}\sigma_3~G$ given above, we conclude that the matching
 is satisfied if and only if:
$$
b(x)\sim -x~\theta_\infty,~~~~~\hbox{ and } \sigma_0=1.  
$$
The error in $b(x)$ is of higher order w.r.t. $x$.  
The determination of the leading behavior of $A_0+A_x$ is complete, because $
c(x)\sim s(s+\theta_x)~b(x)$, $x\to 0$. Namely:
$$
c(x)\sim  - x~ \bigl\{s(s+\theta_x) ~\theta_\infty\bigr\} ,~~~~~a(x)=o(x).
$$

\vskip 0.2 cm 
With such a choice of $b(x)$, one can verify that the terms in
$\Psi_{OUT}$ and $\Psi_{IN}$ which follow the second (i.e. which
follow  the term in
$\lambda$) are actually of higher order in $x$. Nevertheless,
$\Psi_{OUT}$ and $\Psi_{IN}$ match only in the first and second term,
being already the off-diagonal entries of the third term not matching
(i.e. $-A$ and $K_1=$~diagonal part of $-A$ respectively.)

\subsection{Matching for the Second Choice: ~${1\over 2} <
  \delta_{IN}\leq \delta_{OUT} <{2\over 3}$.}

We rewrite the systems in the convenient form:
$$\nu:={1\over \lambda},~~~\mu:={\lambda\over x};$$
$$
{d\Psi_{OUT}\over d \nu}= \left[-x^2A~\nu~
-xA -{A_0+A_x\over \nu} -{\theta_\infty\over 2}{\sigma_3\over \nu(\nu-1)} 
\right] ~\Psi_{OUT},
$$
\vskip 0.2 cm
\be
{d \Psi_{IN} \over d\mu}
=
\left[x~{\theta_\infty\over 2}\sigma_3   
+{A_0\over \mu} +{A_x\over \mu-1} \right]~
\Psi_{IN}
=
\left[x~{\theta_\infty\over 2}\sigma_3 +{A_0+A_x\over \mu} -{A_0\over \mu(\mu-1)} \right]~
\Psi_{IN}. 
\label{tildeinapril}
\ee
$$
\hbox{ Then, substitute $A_0\mapsto -A$ in the last term
}
$$
We rewrite the systems at infinity 
%%%%%%%%%%%%%%%%%%
\footnote{ The system (\ref{tildeinapril}) 
is in the form of a system of isomonodromy deformation for the fifth
Painlev\'e equation. 
}
%%%%%%%%%%%%%%%%%%%%%%%%%%%%%%%%%%%%%%
:
$$
{d\tilde{\Psi}_{OUT}\over d\nu}=
\left[x^2~{\theta_x\over 2}\sigma_3~\nu+
x~{\theta_x\over 2}\sigma_3 -{G^{-1}(A_0+A_x)G\over \nu}
-{\theta_\infty\over 2}~G^{-1}\sigma_3 G~\left({1\over \nu^2}+{1\over \nu^3}+...\right)
\right]
~\tilde{\Psi}_{OUT},
$$
\be
{d \Psi_{IN} \over d\mu}
=
\left[x~{\theta_\infty\over
    2}\sigma_3  +{A_0+A_x\over \mu} +A ~\left({1\over \mu^2}+{1\over \mu^3}+...
\right) \right]~
\Psi_{IN},~~~~~\nu,~\mu\to \infty.
\label{starIN}
\ee
This time the system of $\tilde{\Psi}_{OUT}$ is in the form
(\ref{IRR2}), while the system of $\Psi_{IN}$ is in the form
(\ref{IRR1}), where: 
$$
\Omega=x~{\theta_\infty\over 2} \sigma_3,~~~~~\Omega_1
= \hbox{ diagonal part of } (A_0+A_x),
$$
$$
\Lambda={\theta_x\over
  2} \sigma_3,~~~~~\Lambda_1=  \hbox{ diagonal
  part of } \bigl(-G^{-1}(A_0+A_x) ~G\bigr).
$$
We impose that $\Omega_1$ and $\Lambda_1$ do not depend on
$x$, and we get the conditions $a=0$, $c= s(s+\theta_x)~b$.  
Then, we choose the following solutions: 
$$
\Psi_{OUT}=G
\left[
I+{ K_1\over \nu} ~+...
\right]
\exp
\left\{
x^2~{\theta_x\over 4} \sigma_3~\nu^2~+x~{\theta_x\over 2}
\sigma_3~\nu
\right\} ~G^{-1}~~~~~~~~~~~~~~
$$
$$~~~~~~~~~~
= I+{\theta_x\over 2} G\sigma_3 G^{-1} ~{x\over \lambda} +
GK_1G^{-1} \lambda +O\left({x^2\over
  \lambda^2},x,\lambda^2\right),~~~
\nu\to\infty.
$$
\vskip 0.2 cm 
$$
\Psi_{IN}
= \left[
I+{G_1\over \mu} +...
\right]~
\exp\left\{
x~{\theta_\infty\over 2} \sigma_3 \mu
\right\}= 
I+G_1{x\over \lambda} +{\theta_\infty\over 2} \sigma_3 \lambda +O\left({x^2\over \lambda^2},x,\lambda^2\right),~~~\mu\to\infty.
$$
The relation of dominance among terms are considered as if the
$G_n$'s and $K_n$'s do not diverge as $x\to 0$.
The matching conditition:  
$$
\Psi_{OUT}(\lambda,x)\sim\Psi_{IN}(\lambda,x),~~~~x\to 0,~~~|x|^{\delta_{OUT}}\leq |\lambda|\leq |x|^{\delta_{IN}}.
$$
is restricted to the overlapping sector where both expansions of
$\Psi_{OUT}$ and $\Psi_{IN}$ hold.  
 We note that ${x\over\lambda}$ vanishes slower than $\lambda$, because
$\delta_{IN}>{1\over 2}$ (namely, $  {x/\lambda}\sim x^{1-\delta}$, 
$\lambda\sim x^{\delta}$, $\delta>1/2$). 
 $\Psi_{IN}$ and $\Psi_{OUT}$ automatically match in the first term
$I$. We impose the matching of the second leading term, i.e. the term in
${x\over \lambda}$: 
\be
G_1(x)\sim{\theta_x\over 2} ~G~\sigma_3~G^{-1}\equiv -A,~~~~~x\to 0. 
\label{MMMA}
\ee
 As for
$G_1$, the formulae in Appendix 2 give:
$$
(G_1)_{ij} = -2~{(A_0+A_x)_{ij}\over x~ \theta_{\infty}\bigl[
(\sigma_3)_{ii}-(\sigma_3)_{jj}
\bigr] }, ~~~~~i\neq j.
$$
$$
(G_1)_{ii}= -(A)_{ii}
+
2~{(A_0+A_x)_{ij}~(A_0+A_x)_{ji}
\over 
x ~\theta_\infty \bigl[ (\sigma_3)_{jj}-(\sigma_3)_{ii} \bigr]
}.
$$
In the last formula, $j=2$ if $i=1$, $j=1$ if $i=2$. Explicitly:
$$
(G_1)_{12}= -{r\over x~\theta_\infty}~b,~~~~~(G_1)_{21}=
{s(s+\theta_x)\over x~\theta_\infty~r}~b,
$$
\be
(G_1)_{11}= -(A)_{11} -{s(s+\theta_x) \over x~\theta_\infty}~b^2
,~~~
(G_1)_{22}= -(A)_{22}+{s(s+\theta_x) \over x~\theta_\infty}~b^2.
\label{G1IN}
\ee
Therefore $
\hbox{(\ref{MMMA})}~\Longleftrightarrow~b(x)\sim -x~\theta_\infty$,
and $\sigma_0=1$. 
As it must be, we get the same result of the matching for the first
choice.

\subsection{Critical Matching:~~ ${1\over
    2}-\epsilon<\delta_{IN}\leq \delta_{OUT} <{1\over 2}+\epsilon$.}
\label{criticalmatch}

In between the first and the second choice -- which hold respectively for
${1\over 3}<\delta_{IN}\leq \delta_{OUT}<{1\over 2} $ and ${1\over
  2}<\delta_{IN}\leq \delta_{OUT}<{2\over 3} $ -- we can also consider
the following approximations of system (\ref{SYSTEM}): 
$$
{d\Psi_{OUT}\over d \lambda}= \left[
{xA \over \lambda^2} +{A_0+A_x\over \lambda} -{\theta_\infty\over 2}{\sigma_3\over \lambda-1} 
\right] ~\Psi_{OUT},
$$
$$
{d \Psi_{IN} \over d\lambda}
=
\left[
{A_0\over \lambda} +{A_x\over \lambda-x} +{\theta_\infty\over
  2}\sigma_3 \right]~
\Psi_{IN}.
$$
Rigorously speaking, the two systems 
 cannot be considered simultaneously when
$\sigma_0=1$.  But we can consider $\sigma_0=1$ as a ``limit'' value -- or ``critical''
value -- for the matching of the two above systems in the region
specified by  ${1\over
  2}-\epsilon<
\delta_{IN}\leq \delta_{OUT}<{1\over 2}+\epsilon $, where $\epsilon>0$
 is sufficiently small. 
 We write again $
\Psi_{OUT}=:G~\tilde{\Psi}_{OUT} 
$. 
Then:
$$
{d\tilde{\Psi}_{OUT}\over d\nu}=
\left[
x~{\theta_x\over 2}\sigma_3 -{G^{-1}(A_0+A_x)G\over \nu}
-{\theta_\infty\over 2}~G^{-1}\sigma_3 G~\left({1\over \nu^2}+{1\over \nu^3}+...\right)
\right]
~\tilde{\Psi}_{OUT},~~~\nu\to\infty;
$$
$$
{d \Psi_{IN} \over d\mu}
=
\left[x~{\theta_\infty\over
    2}\sigma_3   +{A_0+A_x\over \mu} +A ~\left({1\over \mu^2}+{1\over \mu^3}+...
\right) \right]~
\Psi_{IN},~~~\mu\to \infty.
$$

When we  impose isomonodromicity conditions,  the diagonal parts of
$A_0+A_x$ and $G^{-1}(A_0+A_x)G$ must be independent of $x$. This
gives again $a=0$, $c=s(s+\theta_x)~b$. 
Then, we choose the fundamental solutions:
$$
\Psi_{OUT}= G\left[
I+{G_1^{OUT}\over \nu} + ...
\right]~\exp\left\{
x~{\theta_x\over 2}\sigma_3~\nu
\right\}~G^{-1}=~~~~~~~~~~~~~~~~~~~~~~~~~ 
$$
$$~~~~~~~~~~~~~~~~~~~~~~~~~~~~~~~~~~~~~
=~I~+ GG_1^{OUT}G^{-1}\lambda+{\theta_x\over 2}
G\sigma_3G^{-1}~{x\over \lambda} +O\left(\lambda^2,x,{x^2\over
  \lambda^2}
\right)
$$
\vskip 0.2 cm
$$
\Psi_{IN}=\left[
I+G_1^{IN}~{x\over \lambda}+...
\right]~
\exp\left\{
x~{\theta_\infty\over 2}\mu 
\right\}
=~
I~+{\theta_\infty\over 2} \sigma_3~\lambda~+G_1^{IN}~{x\over \lambda}
+O\left(\lambda^2,x,{x^2\over
  \lambda^2}\right)
$$

We match them  for $\lambda\sim x^\delta$ and  
 $\delta\in(1/2-\epsilon,1/2+\epsilon)$,
 in the overlapping sector where the above expansions
hold. Here both $\lambda\sim x^\delta$ and 
 $x/\lambda\sim
x^{1-\delta}$ are the dominant terms. The matching conditions are $
G G_1^{OUT}G^{-1}\sim{\theta_\infty\over 2}\sigma_3$ and $G_1^{IN}\sim
{\theta_x\over 2} G \sigma_3 G^{-1}$. 
Namely:
\be
G_1^{OUT}(x)\sim {\theta_\infty\over 2}
G^{-1}\sigma_3~G,~~~~~G_1^{IN}(x)\sim -A,~~~~~~~x\to 0
\label{MAMACH}
\ee
The matrix $ {\theta_\infty\over 2} G^{-1}\sigma_3 G$ can be derived
from (\ref{SIGMA3}). The matrix $G_1^{OUT}$ is (\ref{G1OUT}), the
matrix $G_1^{IN}$ is (\ref{G1IN}). Condition (\ref{MAMACH}) is
inclusive of both (\ref{matchcond1choice}) and
(\ref{MMMA}). Therefore, $
\hbox{(\ref{MAMACH})} ~\Longleftrightarrow ~
b(x)\sim -x~\theta_\infty$, $~x\to 0$. 
This is again the expected result.

\subsection{Higher Order Terms} 
\label{higherorder}

The final result obtained above is: 
\be
A_0+A_x= \pmatrix{ 0 & -r\theta_\infty ~x
                     \cr
\cr
                   -{(s+\theta_x)s\theta_\infty\over r}~x & 0 
}
+
o(x),~~~~~A_1=-{\theta_\infty\over 2}\sigma_3-(A_0+A_x),
\label{temporaryA0Ax}
\ee
\be
A_x=
\pmatrix{ s+{\theta_x\over 2} & -r 
\cr
                   (s+\theta_x)~s\over r &  -s-{\theta_0\over 2}
}
+o(1),~~~A_0=-\pmatrix{ s+{\theta_x\over 2} & -r 
\cr
                   (s+\theta_x)~s\over r &  -s-{\theta_0\over 2}
} +o(1)
\label{temporaryA1}
\ee
Let us substitute the  above results into (\ref{leadingtermaprile}). 
We obtain the first term with no error:
$$
y(x)\sim {1\over 1-\theta_\infty},~~~~~ x\to 0.
$$
Here, $r$ and $s$ do not appear. Nevertheless, if 
 we substitute in (PVI)  the series $y= {1\over
   1-\theta_\infty}+\sum_{n=1}^\infty b_n x^n$, 
we can compute recursively all the 
terms, for  $\theta_0=\pm\theta_x$ and $\theta_1=\pm\theta_\infty$. We find a
series:
 \be
   y(x)= {1\over 1-\theta_\infty}+a~x +\sum_{n=0}^\infty 
b_n(a;\theta_\infty,\theta_0)~x^n,~~~x\to 0,
\label{temporaryy}
\ee
where $a$ is an arbitrary parameter. This parameter is actually a
function of $s$, as we prove now. 
The convergence of the Taylor expansion can be proved by a
Briot-Bouquet like argument. 
This will not be done here. The reader can find a similar proof in
\cite{kaneko} and the general procedure in \cite{INCE}.

\subsubsection{Determination of $a=a(s)$}

The system (\ref{SYSTEM}) is
isomonodromic. This determines the structure $A_x$, $A_0$ $A_1$ as
can be found in \cite{JMU}, Appendix C, formulae (C.47), (C.49),
(C.51), (C.52), (C.55). 
If we substitute (\ref{temporaryy}) in the formulae, we get a Taylor
expansion for the matrix elements, in terms of the parameter $a$. 
The leading terms have exactly the
structure of  
(\ref{temporaryA0Ax}) and (\ref{temporaryA1}). We can identify the
leading terms to express $a$  as a function 
of $s$ and $r$. The computations are quite long, so we give
the result. When we write  the leading terms as a function of $a$ and
impose that they coincide with  (\ref{temporaryA0Ax}) and
(\ref{temporaryA1}), we find:  
$$
a=  {\theta_\infty(2s+\theta_x+1)\over 
2(\theta_\infty-1)}~\in{\bf C}
$$
The higher order terms are  Taylor expansions. 
%%%%%%%%%%%%%
%%%%%%%%
%%%%%%      CANCELED IN THE PREPRINT
%%%%%%%%
%%%%%%%%%%%%%  
%$$
%A_0+A_x= \pmatrix{ 0 & -r\theta_\infty ~x
%                     \cr\cr
%                   -{(s+\theta_x)s\theta_\infty\over r}~x & 0 
%} +\sum_{n=2}^\infty \Delta_n x^n, 
%$$
%$$
%A_1=-{\theta_\infty\over 2}\sigma_3 -(A_0+A_x),~~~~~
%A_x= \pmatrix{ s+{\theta_x\over 2} & -r 
%\cr
%                   (s+\theta_x)~s\over r &  -s-{\theta_0\over 2}
%}+\sum_{n=1}^{\infty}{\cal A}_nx^n,
%$$
%where 
% ${\cal A}_n$ and  $\Delta_n$ are matrices. 
%%%%%%%%%%%%%%%%
 Explicitly, the first terms are: 
\be
A_1=-{\theta_\infty\over 2}\sigma_3- (A_0+A_x)~=
\label{AS1}
\ee
$$=
 \pmatrix{
-{\theta_\infty\over 2}+(s+\theta_x)s\theta_\infty~x^2 
& 
r\theta_\infty
\left\{
x - {(\theta_\infty+1)(2s+\theta_x-1)\over 2}~x^2
\right\}
\cr
\cr
{s\theta_\infty\over r}\left\{
(\theta_x+s)~x
-{(\theta_\infty-1)(\theta_x+s)(2s+\theta_x+1)s\over 2}~x^2
\right\} 
&
{\theta_\infty\over 2}-  (\theta_x+s)s\theta_\infty~ x^2
}~+O(x^3)
$$
\vskip 0.3 cm
\be
A_x=\pmatrix{
\left(s+{\theta_x\over 2}\right) - 2(s+\theta_x)s\theta_\infty ~x
&
-r \left\{1-(2s+\theta_x-1)\theta_\infty~x\right\}
\cr
\cr
{s(s+\theta_x)\over r}\left\{1-(2s+\theta_x+1)\theta_\infty~x
\right\}
&
-\left(s+{\theta_x\over 2}\right) + 2(s+\theta_x)s\theta_{\infty}~ x
}~+O(x^2),
\label{ASX}
\ee
\vskip 0.3 cm
\be
A_0=\pmatrix{
-\left(s+{\theta_x\over 2}\right) + 2(s+\theta_x)s\theta_{\infty}~ x
&
r \left\{1-(2s+\theta_x)\theta_\infty~x\right\}
\cr
\cr
-{s(s+\theta_x)\over r}\left\{1-(2s+\theta_x)\theta_\infty~x
\right\}
&
\left(s+{\theta_x\over 2}\right) - 2(s+\theta_x)s\theta_\infty ~x
}~+O(x^2).
\label{AS0}
\ee
%%%%%
The above expansions  are enough to obtain first two leading terms of 
 (\ref{leadingtermaprile}):
$$
y(x)= {1\over 1-\theta_\infty} +{\theta_\infty(2s+\theta_x+1)\over 
2(\theta_\infty-1)}~x ~+O(x^2),
$$
Note that $r$ simplifies. This is the solution (\ref{form2}). 

%%%%%%%%%%%%%%%%%%%%%%%%%%%%%%%%
%%%%%%%%%%%%%%%%%%%
%%%%%%%%%%%           CANCELED IN THE PREPRINT
%%%%%%%%%%%%%%%%%%%
%%%%%%%%%%%%%%%%%%%%%%%%%%%%%%%%
%{\it Let $\theta_\infty$ and $\theta_x$ be given,  
% ${\theta_1}^2={\theta_\infty}^2$, ${\theta_0}^2={\theta_x}^2\neq 0$,
%  and 
%$\lim_{x\to 0}(A_0+A_x)=0$. 
% (PVI) has  a 1-parameter family of solutions with behavior: 
%\be
%   y(x)= {1\over 1-\theta_\infty}+a~x +\sum_{n=2}^\infty 
%c_n(a;~\theta_\infty,\theta_x)~x^n,~~~x\to 0,
%\label{sperosiavero}
%\ee
%where $c_n(a;~\theta_\infty,\theta_x)$ are certain 
%rational function of  $\theta_\infty$, $\theta_x$,  the parameter is:
%$$
%a= {\theta_\infty(2s+\theta_x+1)\over 
%2(\theta_\infty-1)}~\in{\bf C},$$
%and the matrices  $A_1(x)$, $A_x(x)$, $A_0(x)$ are (\ref{AS1}),
%(\ref{ASX}),  (\ref{AS0}),  
%}   
%
%\vskip 0.2 cm
%\noindent 
%Note that for $a=0$, $c_n(0;~\theta_\infty,\theta_x)\neq 0$. For
%example:
%$$
%c_2(a;~\theta_\infty,\theta_x)=
%{ \theta_\infty  (2 \theta_\infty-1) ({\theta_x}^2-1)\over
%  6~(\theta_\infty -1)}+
%\theta_\infty~ a+(1-\theta_\infty)~a^2. 
%$$
%

\subsection{Monodromy Data}

We assume that the matching has been completed as above, and in
particular $\sigma_0=1$. 
Thus, the system (\ref{SYSTEM}) can be approximated by:
\be
{d\Psi_{OUT}\over d\lambda} = \left[
{A_0+A_x\over \lambda}+{xA_x\over \lambda^2} -{\theta_\infty\over 2}~
 {\sigma_3\over \lambda-1}
\right]~\Psi_{OUT},
~~~\hbox{ for } |\lambda|\geq |x|^{\delta},~~~\delta<{1\over
  2}; 
\label{primadibumOUT}
\ee
or
\be
{d\Psi_{IN}\over d\lambda}=
\left[
{A_0\over \lambda}+{A_x\over \lambda} +{\theta_\infty\over 2}\sigma_3
\right]~\Psi_{IN},~~~\hbox{ for }~|\lambda|\leq |x|^{\delta},~~~\delta>{1\over
  2}; 
\label{primadibumIN}
\ee
The first two leading terms  are
:
\be
\Psi_{OUT}= I + {\theta_\infty\over 2} \sigma_3 ~\lambda
+O\left({x\over \lambda}\right),~~~\lambda\sim x^{\delta}\to 0,
~~~\delta<{1\over 2};
\label{finalmatchOUT}
\ee
\be
\Psi_{IN}= I -A~{x\over \lambda} + O(\lambda), ~~~\lambda\sim
x^{\delta}\to 0,
~~~\delta>{1\over 2}.
\label{finalmatchIN}
\ee
The above solutions  match in the first three terms in the 
``critical'' region $\lambda\sim x^{\delta}$, $\delta\simeq
   {1\over 2}$ (the region is restricted to a sector). Namely:
\be
   \Psi_{OUT}\sim\Psi_{IN} \sim  I + {\theta_\infty\over 2} \sigma_3 ~\lambda
-A~{x\over \lambda} ~+O\left(x,{x^2\over
   \lambda^2},\lambda^2\right),~~~\lambda\sim x^{\delta},~~\delta\simeq
   {1\over 2}
\label{girigirimatch}
\ee

\vskip 0.3 cm 
\noindent
Now, for $\delta<{1\over 2}$, we have: 
$$
{A_0+A_x\over \lambda}\sim x^{1-\delta}\to
0,~~~~~ 
{xA_x\over\lambda^2}\sim x^{1-2\delta}\to 0.
$$ 
Thus, (\ref{primadibumOUT}) can be further reduced to:
\be
{d\hat{\Psi}_{OUT}\over d\lambda}= -{\theta_\infty\over 2}
~{\sigma_3\over \lambda-1} ~\hat{\Psi}_{OUT}
\label{pesticidaOUT}
\ee

In system (\ref{primadibumIN}) we rewrite:
$$
{A_0\over \lambda}+{A_x\over \lambda-x}= {A_0+A_x\over
  \lambda-x}-{xA_0\over \lambda(\lambda-x)}
$$
Then, we observe that,  for  $\delta>{1\over 2}$, we have the behaviors:
$$
{A_0+A_x\over
  \lambda-x}\sim x^{1-\delta}\to 0,
~~~
  {xA_0\over \lambda(\lambda-x)} \sim -{x A\over \lambda(\lambda-x)}
  \sim
x^{1-2\delta}\to \infty,~~~~~\delta>{1\over 2}
$$
Thus, as $x\to 0$, the system (\ref{primadibumIN}) can be further
reduced to:
\be
{d\hat{\Psi}_{IN}\over d\lambda} = {x A\over \lambda(\lambda-x)}
~\hat{\Psi}_{IN}
.
\label{pesticidaIN}
\ee

\vskip 0.3 cm
System (\ref{pesticidaOUT}) has the same monodromy of (\ref{SYSTEM})
at $\lambda=1,\infty$. System (\ref{pesticidaIN}) has the same
monodromy of 
(\ref{SYSTEM})
at $\lambda=0,x$.

\vskip 0.5 cm 
\noindent
{\bf MATCHING $\Psi\leftrightarrow \hat{\Psi}_{OUT}$:} 

\vskip 0.2 cm 
We choose $\hat{\Psi}_{OUT}$ such that it matches with
$\Psi$ at
$\lambda=\infty$. The behavior of this last is: 
$$
\Psi(\lambda,x)=\left[I+O\left({1\over \lambda}\right)\right]
\lambda^{-{\theta_\infty\over 2}\sigma_3},~~~~~\lambda\to \infty,
$$
for $\theta_\infty\not\in{\bf Z}$. 
The  solution of  (\ref{pesticidaOUT}) with the same behavior is:
\be
\hat{\Psi}_{OUT}(\lambda):= (\lambda-1)^{-{\theta_\infty\over 2}\sigma_3} 
\label{PSIHATOUT}
\ee
As a consequence, the monodromy of $\Psi$ at $\lambda=1,\infty$
coincides with that  of $\hat{\Psi}_{OUT}$. To compute it, we consider the
loops  $\lambda-1\mapsto
(\lambda-1)e^{2\pi i }$ and $\lambda\mapsto \lambda e^{2\pi i}$. The
corresponding monodromy is: 
$$
M_1= M_{\infty}=~\exp\bigl\{
i\pi\theta_\infty\sigma_3
\bigr\}.
$$ 

\vskip 0.5 cm
\noindent
{\bf MATCHING $\hat{\Psi}_{OUT}\leftrightarrow \Psi_{OUT}$:} 

\vskip 0.3 cm 

We match $\hat{\Psi}_{OUT}$ with $\Psi_{OUT}$ 
for $x/\lambda\to 0$.   
Let us choose the branch $(\lambda-1)=(1-\lambda)~e^{i\pi}$, $(1-\lambda)>0$ for
$0<\lambda<1$. Solution (\ref{PSIHATOUT}) has expansion:
$$
\hat{\Psi}_{OUT}(\lambda)=
e^{-i{\pi\over 2}\theta_\infty \sigma_3} \left[I+{\theta_\infty\over
    2}\sigma_3\lambda+O(\lambda^2)\right], ~~~\lambda\to 0
$$
Therefore, for $\lambda\to 0$,
 $\hat{\Psi}_{OUT}$ matches with $\Psi_{OUT}~ e^{-i{\pi\over
    2}\theta_\infty\sigma_3}$, where   $\Psi_{OUT}$ is 
(\ref{finalmatchOUT}).

\vskip 0.5 cm
\noindent
{\bf MATCHING $\Psi_{OUT}\leftrightarrow \Psi_{IN}$:} ~~This is
(\ref{girigirimatch}).

\vskip 0.5 cm
\noindent
{\bf MATCHING $\Psi_{IN}\leftrightarrow \Psi$:} ~~
 The above matchings   imply that 
$\Psi_{IN}~ e^{-i{\pi\over
    2} \theta_\infty \sigma_3}$ matches with $\Psi$  
(where $\Psi_{IN}$ is (\ref{finalmatchIN})).

\vskip 0.5 cm
\noindent
{\bf MATCHING $\Psi_{IN}\leftrightarrow \hat{\Psi}_{IN}$:} 

\vskip 0.3 cm

 In order to determine the
monodromy of system (\ref{SYSTEM}) at $\lambda=0,x$, we need to find a
  fundamental solution $\hat{\Psi}_{IN}$  of (\ref{pesticidaIN}), that   
matches with  $\Psi_{IN}~ e^{-i{\pi\over
    2}\sigma_3}$
 for $x/\lambda\to 0$, where $\Psi_{IN}$ is (\ref{finalmatchIN}). A
 fundamental solution satisfying these requirements is:
$$
\hat{\Psi}_{IN}(\lambda,x):=~G\left[\left(1-{x\over
  \lambda}\right)^{-{\theta_x\over 2}\sigma_3}\right]G^{-1}~e^{-i{\pi \over
      2}\theta_\infty\sigma_3}.
$$
Actually, this has the behavior:
$$
\hat{\Psi}_{IN}(\lambda,x)= \left[I +{\theta_x\over 2}
  G\sigma_3G^{-1}~{x\over \lambda}~+O\left({x^2\over
    \lambda^2}\right)
\right]~e^{-i{\pi\over 2}\theta_\infty \sigma_3},~~~~~~~~{x\over \lambda}\to 0,
$$
where $
{\theta_x\over 2}
G\sigma_3G^{-1}\equiv -A$. 
The first two terms match with $\Psi_{IN} ~e^{-i{\pi\over
    2}\theta_\infty \sigma_3}$, as required. 
\vskip 0.2 cm
As a consequence of the matching, the monodromy matrices  of $\Psi$ and $\hat{\Psi}_{IN}$ at $\lambda=0,x$
coincide. To compute them, we write the local behaviors (for $x\neq 0$
fixed) of $\hat{\Psi}_{IN}$:
$$
\hat{\Psi}_{IN}(\lambda,x)=
\left\{
\matrix{
\Bigr[Gx^{{\theta_x\over 2}\sigma_3}\Bigl]~\bigl[1+O(\lambda-x) \bigr]
~(\lambda-x)^{-{\theta_x\over 2}\sigma_3} ~\left[G^{-1} e^{-i{\pi
      \over 2} \theta_\infty\sigma_3}
\right],~~~\lambda\to x,
\cr
\cr
\Bigr[G(-x)^{-{\theta_x\over 2}\sigma_3}\Bigr]~\bigl[1+O(\lambda) \bigr]~
\lambda^{{\theta_x\over 2}\sigma_3}~ \left[G^{-1} e^{-i{\pi
      \over 2} \theta_\infty\sigma_3}
\right],~~~\lambda\to 0.
}
\right.
$$
It is not necessary to specify the branch of $
(\pm x)^{-{\theta_x\over 2}\sigma_3}$. We resctict to 
the case $\theta_0$ and $\theta_x\not\in{\bf Z}$, so that the
matrices $R_0$, $R_x$ in  (\ref{PSIlocale}) 
are zero and the matching of $\Psi$  with the above
behaviors of $\hat{\Psi}_{IN}$ is realized. 

We can compute the monodromy
matrices for $\lambda\mapsto \lambda e^{2\pi i}$ and
$(\lambda-x)\mapsto (\lambda-x)e^{2\pi i}$ respectively:
$$
M_0= e^{i{\pi\over 2}\theta_\infty\sigma_3}G~\exp\bigl\{
 i\pi\theta_x\sigma_3\bigr\}~G^{-1} e^{-i{\pi\over
     2}\theta_\infty\sigma_3},
$$
$$
M_x= e^{i{\pi\over 2}\theta_\infty\sigma_3}G~\exp\bigl\{-
 i\pi\theta_x\sigma_3\bigr\}~G^{-1} e^{-i{\pi\over
     2}\theta_\infty\sigma_3},
$$
Note that $M_1=M_{\infty}$ is invariant for $M_1\mapsto e^{-i{\pi\over
     2}\theta_\infty\sigma_3}M_1e^{i{\pi\over
    2}\theta_\infty\sigma_3}$.
With this fact in mind, we obtain the result of theorem
\ref{thMONODROMY}, point b). In particular, computing the trace of
$M_0M_1$,
 we get:
$$
 s={\theta_x\bigl[2\cos(\pi(\theta_\infty+\theta_x))-\hbox{\rm
    tr}(M_1M_0)\bigr]\over
 2\bigl[\cos(\pi(\theta_\infty-\theta_x))-\cos(\pi(\theta_\infty+\theta_x))\bigr]}.
$$

%%%%%%%%%%%%%%%%%%%%%%%%%%%%%%%%%%%%%%%%%%%%%%

\section{Case
  $\sigma=\pm(\theta_1-\theta_\infty),~\pm(\theta_1+\theta_\infty)$.~
Solution (\ref{form1}).}
\label{semireducible}
 This case shows new features, 
namely $r$ and $r_1$  may be functions of $x$. 
For $\sigma=\pm(\theta_1-\theta_\infty),~\pm(\theta_1+\theta_\infty)$
the matrices (\ref{hatA1})  (\ref{hatA0Ax}) become: 
\be
  \sigma=\pm(\theta_1-\theta_\infty):
~~~
\hat{A}_1=\pmatrix{-{\theta_1\over 2} & -r_1 
\cr 
              0 & {\theta_1\over 2}},
~~
\hat{A}_0+\hat{A}_x= \pmatrix{{\theta_1-\theta_\infty\over 2} & r_1 
\cr
                                 0 &- {\theta_1-\theta_\infty\over 2}
}.
\label{hatA1A0plusAx}
\ee
$$
  \sigma=\pm(\theta_1+\theta_\infty):
~~~~~
\hat{A}_1=\pmatrix{{\theta_1\over 2} & -r_1 
\cr 
              0 & -{\theta_1\over 2}},
~~~
\hat{A}_0+\hat{A}_x= \pmatrix{-{\theta_1+\theta_\infty\over 2} & r_1 
\cr
                                 0 & {\theta_1+\theta_\infty\over 2}
}.
$$
The transpose matrices may be considered (namely, re-define 
$r_1 \mapsto$$~$$ -{[\sigma^2-(\theta_1-\theta_\infty)^2][\sigma^2-(\theta_1+\theta_\infty)^2]\over 16 \theta_\infty^2 r_1}$). 

\vskip 0.2 cm 
\noindent
The matrices $\hat{\hat{A}_0}$ and  $\hat{\hat{A}_x}$ are again 
(\ref{hathatA0}) and (\ref{hathatAx}). 
For definiteness, we will consider the upper triangular matrices
$\hat{A_1}$ and $\hat{A_0}+\hat{A}_x$, and the choice $
\sigma=\theta_1-\theta_\infty$.  
We distinguish three cases:

\vskip 0.3 cm

{\bf (I) $r_1$ is a non zero constant:} In this case we just repeat
the general matching procedure and find $y(x)$ as in
(\ref{genericlead}). $r\neq 0$ is constant.

\vskip 0.3 cm
{\bf (II) $r_1=0$ and $r$ constant.} This is a very easy case, because 
$
\hat{A_1}=-{\theta_1\over 2}\sigma_3$, 
$\hat{A_0}+\hat{A_x}={\theta_1-\theta_\infty\over 2} \sigma_3$. 
Therefore a fundamental solution of (\ref{system1}) is:
$$
\Psi_{OUT}(\lambda)= (\lambda-1)^{-{\theta_1\over
    2}\sigma_3}\lambda^{{\theta_1-\theta_\infty\over
    2}\sigma_3}~(-1)^{-{\theta_1\over 2}\sigma_3} 
 = \left[ I + {\theta_1\over 2} \sigma_3 + O(\lambda^2)\right] \lambda^{{\theta_1-\theta_\infty\over 2}\sigma_3},~~~\lambda\to 0.
$$
The solution of (\ref{system0}), with $r$ constant, 
 is the same of the general case. 
The matching is possible as in the 
general case, with $K_0= \pmatrix{x^{\theta_1-
\theta_\infty\over 2} & 0 \cr 0 &x^{\theta_\infty-
\theta_1\over 2}} $. Therefore: 
$$
 \hat{A_0}=  
\pmatrix{{\theta_0^2-\theta_x^2+\sigma^2\over 4\sigma}
&
r~x^{\theta_1-\theta_\infty}
\cr
-{[\sigma^2-(\theta_0-\theta_x)^2][\sigma^2-(\theta_0+\theta_x)^2]
\over 16 \sigma^2 r}~x^{\theta_\infty-\theta_1}
&
-{\theta_0^2-\theta_x^2+\sigma^2\over 4\sigma}
}.
$$
We obtain: $
y(x)\sim {x (\hat{A}_0)_{12}\over x[(\hat{A}_0)_{12}+(\hat{A}_1)_{12}]
-(\hat{A}_1)_{12}} ={x~ r~x^{\theta_1-\theta_\infty}\over
  x[~x^{\theta_1-\theta_\infty}+0]-0}=1$. 
This is the {\it singular solution} $y(x)\equiv 1$.
 
The case $r_1=r_1(x)\to 0 $ may give a non-singular
$y(x)$, provided
that also $r=r(x)\to 0$. This will be proved below.

\vskip 0.3 cm
{\bf (III) Case $r_1=r_1(x)$.}
\vskip 0.2 cm
 A priori, we don't know if it is possible to have an $x$-dependent $r_1$. This is not actually possible if we approximate 
(\ref{SYSTEM}) with the  fuchsian system (\ref{system1})
 (i.e. (\ref{fuchsianSYSTEMOUT})).   
Actually,  (\ref{system1}) is a system with reducible monodromy (upper triangular):
$$
{d\Psi_{OUT}\over\ d \lambda}= \left\{
 \pmatrix{{\theta_1-\theta_\infty\over 2} & r_1 
\cr
                                 0 &- {\theta_1-\theta_\infty\over 2}}{1\over\lambda}+
\pmatrix{-{\theta_1\over 2} & -r_1 
\cr 
              0 & {\theta_1\over 2}}{1\over \lambda-1}
\right\}\Psi_{OUT}.
$$
There is a fundamental solution (obtained by variation of parameters):
$$
\Psi_{OUT}=\left[\matrix{
(1-\lambda)^{-{\theta_1\over 2}} &  0 
\cr
\cr
0 &  (1-\lambda)^{\theta_1\over 2}
}
\right]
\left[
\matrix{
1 & {r_1  F(1-\theta_1,\theta_\infty-\theta_1,\theta_\infty-\theta_1+1;\lambda)  \over \theta_\infty-\theta_1}
\cr
\cr
0 & 1 
}
\right]
\left[
\matrix{
\lambda^{\theta_1-\theta_\infty\over 2} & 0 \cr \cr
0 & \lambda^{-{\theta_1-\theta_\infty\over 2}}
}
\right].
$$
Here, $F(...)$ denotes a Gauss hyper-geometric function. 
This solution 
has a diagonal monondromy matrix at $\lambda=0$ and an upper triangular 
monodromy matrix at $\lambda=1$, with $r_1$ in the $(1,2)$ matrix element. Therefore, $r_1$ must be independent of the monodromy preserving deformation $x$. We are sent back to the case {\bf (I)} of constant $r_1$.  

\vskip 0.3 cm
 
The only possibility for $r_1$ to depend on $x$ is that  the matrices 
$\hat{A}_0$, $\hat{A}_x$ have a behavior, for $x\to 0$, such that the
system (\ref{SYSTEM}) is approximated by a system
(\ref{nonfuchsianSYSTEMOUT})
 with singularity of the second kind at $\lambda=0$. Namely, system (\ref{SYSTEM})
 must be approximated {\it at least} by: 
\be
{d\Psi_{OUT}\over d\lambda}= 
\left\{
{xA_x \over \lambda^2}+
{A_0+A_x\over \lambda}+
{A_1\over \lambda-1}
\right\}
\Psi_{OUT}
\label{systemxAx}
\ee

\vskip 0.3 cm
\noindent 
{\bf Hypothesis:} {\it  We consider the case:}
$$
  x A_x= \pmatrix{0 & \rho(x) \cr 0 & 0} + \hbox{ higher orders } ~~~x\to 0.
$$

\noindent
In the above hypothesis,  (\ref{systemxAx}) is:
\be
{d\Psi_{OUT}\over d\lambda} = \left\{  \pmatrix{0 & \rho(x) \cr 0 & 0}{1\over \lambda^2}+
\pmatrix{{\theta_1-\theta_\infty\over 2}  & r_1 \cr 0 & {\theta_1-\theta_\infty\over 2} }{1\over \lambda}+ \pmatrix{-{\theta_1\over 2} & -r_1 \cr 0 &{\theta_1\over 2} } {1\over \lambda-1}
\right\}\Psi_{OUT}
\label{irregolaresystem1}
\ee
Let us write a $\Psi_{OUT}$ as a vector
$\pmatrix{\psi_1\cr\psi_2}$. The system becomes: 
{\small
$$
{d\psi_1\over d\lambda}=\left({\theta_1-\theta_\infty\over 2\lambda}-{\theta_1\over 2(\lambda-1)}\right)\psi_1+\left(
{r_1\over \lambda} +{\rho\over \lambda^2}-{r_1\over \lambda-1}
\right) \psi_2,~~~
{d\psi_2\over d \lambda}= \left(
{\theta_\infty-\theta_1\over 2 \lambda}+ {\theta_1\over 2(\lambda-1)}
\right)\psi_2.
$$
}
This system is solvable by variation of parameters.   
Let $C_1$ and $C_2$ be integration constants. The general solution is:  
\be
\psi_2(\lambda)= C_2~ \lambda^{\theta_\infty-\theta_1\over 2} (\lambda-1)^{\theta_1\over 2},
\label{psi2}
\ee
$$
\psi_1(\lambda)= C_1~\lambda^{\theta_1-\theta_\infty\over 2}(\lambda-1)^{-{\theta_1\over 2}}~+
C_2~ e^{i\pi \theta_1}\left[
{\rho\over \theta_\infty-\theta_1-1}~{ F(1-\theta_1,\theta_\infty-\theta_1-1,\theta_\infty-\theta_1:\lambda)\over \lambda}~+
\right.
$$
\be
\left.+~{r_1-\rho\over \theta_\infty-\theta_1}~F(1-\theta_1,\theta_\infty-\theta_1,\theta_\infty-\theta_1+1;\lambda)
\right]\lambda^{\theta_\infty-\theta_1\over 2} (\lambda-1)^{-{\theta_1\over 2}}.
\label{psi1}
\ee
Here $F(...)$ denotes a Gauss Hyper-geometric equation. 
The choice of the branch is such that for $0<\lambda<1$, we have $0<1-\lambda= e^{-i\pi 
}(\lambda-1)$. 

\vskip 0.2 cm
In order to write the  local behavior for $\lambda \to 
0$, we expand the hypergemoetric functions and 
$(\lambda-1)^{\pm \theta_1/2} =\bigl( e^{i\pi}(1-\lambda)\bigr)^{\pm \theta_1/2}$:
$$
\psi_1=C_1~ e^{-i{\pi\over 2} \theta_1}\lambda^{\theta_1-\theta_\infty\over 2}
\left(1+{\theta_1\over 2}\lambda+O(\lambda^2)\right) +
C_2~ e^{i{\pi\over 2} \theta_1}\lambda^{\theta_\infty-\theta_1\over 2}
\left[
-{\rho\over 1+\theta_1-\theta_\infty}~ {1\over\lambda}~ + 
\right.
$$
$$
~~~~~~~~~~~~~~~~~+\left.
{\rho ~\theta_1(\theta_1-\theta_\infty+2)-2r_1~(\theta_1-\theta_\infty+1)\over 2(\theta_1-\theta_\infty)(\theta_1-\theta_\infty+1)}+O(\lambda)
\right].
$$
\vskip 0.2 cm 
$$
\psi_2=C_2~ e^{i{\pi\over 2}\theta_1}\left(1-{\theta_1\over 2}\lambda+O(\lambda^2)\right)~\lambda^{\theta_\infty-\theta_1\over 2}.
~~~~~~~~~~~~~~~~~~~~~~~~~~~~~~~~~~~~~~~~~~~$$
Therefore, we can take as a fundamental solution the following matrix (choose 
$C_1=e^{i{\pi\over 2}\theta_1}$, $C_2=e^{-i{\pi\over 2}\theta_1}$): 
$$
\Psi_{OUT}=\left[\matrix{
\left(1+{\theta_1\over 2}\lambda +...
\right)\lambda^{\theta_1-\theta_\infty\over 2} 
&
\left[-{\rho\over 1+\theta_1-\theta_\infty} {1\over\lambda}+{\rho ~\theta_1(\theta_1-\theta_\infty+2)-2r_1~(\theta_1-\theta_\infty+1)\over 2(\theta_1-\theta_\infty)(\theta_1-\theta_\infty+1)}+...\right]\lambda^{\theta_\infty-\theta_1\over 2}
\cr
\cr
0 
&
\left(1-{\theta_1\over 2} \lambda+...\right) \lambda^{\theta_\infty-\theta_1\over 2}
}\right]
$$
\vskip 0.3 cm
$$
=
\left\{
\pmatrix{0 & -{\rho\over 1+\theta_1-\theta_\infty}
\cr\cr
0 & 0 
}
{1\over \lambda} + 
\pmatrix{ 1 & {\rho ~\theta_1(\theta_1-\theta_\infty+2)-2r_1~(\theta_1-\theta_\infty+1)\over 2(\theta_1-\theta_\infty)(\theta_1-\theta_\infty+1)}
\cr
\cr
0
 & 
1 }
+O(\lambda)
\right\}
~
\lambda^{
{\theta_1-\theta_\infty\over 2}~\sigma_3}
$$.

\subsection{Matching}

The above solution must be matched with the solution of system (\ref{system0}), with $\sigma=\theta_1-\theta_\infty$. 
$$
\Psi_0\left({\lambda\over x}\right)=\left[
I+K_1{x\over \lambda} +O\left({x^2\over \lambda^2}\right)
\right] \pmatrix{
\left({\lambda\over x}\right)^{\theta_1-\theta_\infty\over 2} & 0 
\cr
0 &  \left({\lambda\over x}\right)^{\theta_\infty-\theta_1\over 2}}.
$$
From a standard computation we find: 
$$
K_1= 
\pmatrix{ {\theta_0^2-(\theta_x-\sigma)^2\over 4\sigma} &
{r \over \sigma+1} 
\cr
{(\theta_0^2-(\theta_x-\sigma)^2)(\theta_0^2-(\theta_x+\sigma)^2)\over 
16\sigma^2(\sigma-1)r}
&
{(\theta_x+\sigma)^2-\theta_0^2\over 4\sigma(\sigma+1)}
},~~~~~r\neq 0.
$$
Note that $r$ in (\ref{hathatAx}) is constant, because the monodromy of (\ref{system0}) depends on $r$.

\vskip 0.2 cm 
The matching relation $
\Psi_{OUT}(\lambda)\sim K_0(x) \Psi_0\left({\lambda/ x}\right)$, 
reads: 
$$
\left\{
\pmatrix{0 & -{\rho\over 1+\theta_1-\theta_\infty}
\cr\cr
0 & 0 
}
{1\over \lambda} + 
\pmatrix{ 1 & {\rho ~\theta_1(\theta_1-\theta_\infty+2)-2r_1~(\theta_1-\theta_\infty+1)\over 2(\theta_1-\theta_\infty)(\theta_1-\theta_\infty+1)}
\cr
\cr
0
 & 
1 }
+O(\lambda)
\right\}~\sim~~~~~~~~~~~~~~~~~~
$$
$$
~~~~~~~~~~~~~~~~~~~~\sim~
K_0(x)~ \left[
I+ K_1{x\over \lambda}+...
\right]
\pmatrix{
x^{\theta_\infty-\theta_1 \over 2} &  0 
\cr
0 &
x^{\theta_1-\theta_\infty\over 2}
}.
$$
Namely:
$$
\left\{
\matrix{
\pmatrix{ 1 & {\rho ~\theta_1(\theta_1-\theta_\infty+2)-2r_1~(\theta_1-\theta_\infty+1)\over 2(\theta_1-\theta_\infty)(\theta_1-\theta_\infty+1)}
\cr
\cr
0 & 1 }\sim K_0(x)~\pmatrix{
x^{\theta_\infty-\theta_1 \over 2} &  0 
\cr
0 &
x^{\theta_1-\theta_\infty\over 2}
}
\cr
\cr
\cr
\pmatrix{0 & -{\rho\over 1+\theta_1-\theta_\infty}
\cr\cr
0 & 0 }\sim ~x~ K_0(x)~F_1^{(0)}~\pmatrix{
x^{\theta_\infty-\theta_1 \over 2} &  0 
\cr
0 &
x^{\theta_1-\theta_\infty\over 2}
}~~~~~~~~~~~~~~~~
}
\right.
$$
The first equation above is:
$$
  K_0(x)\sim  \pmatrix{ 1 & {\rho
  ~\theta_1(\theta_1-\theta_\infty+2)-2r_1~(\theta_1-\theta_\infty+1)\over
  2(\theta_1-\theta_\infty)(\theta_1-\theta_\infty+1)} 
\cr
\cr
0 & 1 }~\pmatrix{
x^{\theta_1-\theta_\infty \over 2} &  0 
\cr
0 &
x^{\theta_\infty-\theta_1\over 2}
}.
$$
We substitute this in the second equation. For simplicity,  denote $K_{ij}$ the
matrix elements of  
$K_1$. We obtain:
$$
\pmatrix{0 & -{\rho\over 1+\theta_1-\theta_\infty}
\cr\cr
0 & 0 }
\sim
% \pmatrix{ 1 & {\rho
% ~\theta_1(\theta_1-\theta_\infty+2)-2r_1~(\theta_1-\theta_\infty+1)\over 
%2(\theta_1-\theta_\infty)(\theta_1-\theta_\infty+1)}  
%\cr
%\cr
%0 & 1 } 
%\pmatrix{
%x~K_{11} & K_{12} x^{\theta_1-\theta_\infty+1} \cr
%\cr
%K_{21} x^{\theta_\infty-\theta_1+1} & x~K_{22}}
%$$
%$$
%=
   \pmatrix{x~K_{11}+ (*)~K_{21}
~x^{\theta_\infty-\theta_1+1}
&
K_{12}~x^{\theta_1-\theta_\infty+1}+(*)~x~K_{22}
\cr
\cr
K_{21}~x^{\theta_\infty-\theta_1+1} 
&
x~K_{22}
},
$$
where:  
$$
(*):=
 {\rho ~\theta_1(\theta_1-\theta_\infty+2)-2r_1~(\theta_1-\theta_\infty+1)\over 2(\theta_1-\theta_\infty)(\theta_1-\theta_\infty+1)}.
$$
\vskip 0.2 cm 
\noindent
The element $(2,1)$ must vanish. This occurs iff: 
\be
    x^{\theta_\infty-\theta_1+1}\to 0, ~~~\hbox{ for } x\to 0;~~\Longleftrightarrow~~
\Re (\theta_\infty-\theta_1)>-1.
\label{conditionmatch}
\ee
\vskip 0.2 cm
\noindent
We substitute this result in the element $(1,1)$ and then we impose
that it vanishes: 
$$
 {\rho
 ~\theta_1(\theta_1-\theta_\infty+2)-2r_1~(\theta_1-\theta_\infty+1)\over
 2(\theta_1-\theta_\infty)(\theta_1-\theta_\infty+1)} 
~K_{21} ~x^{\theta_\infty-\theta_1+1}\to 0. 
$$
This implies that $\rho
~\theta_1(\theta_1-\theta_\infty+2)-2r_1~(\theta_1-\theta_\infty+1)=
o(x^{\theta_1-\theta_\infty-1})$.  

\noindent
From the element $(1,2)$ we have: 
\be
-{\rho\over 1+\theta_1-\theta_\infty}\sim K_{12}~x^{\theta_1-\theta_\infty+1}
+o(x^{\theta_1-\theta_\infty})~K_{22}.
\label{maldipancia}
\ee
This relation may be satisfied in  two ways: the first is that $\rho=\rho(x)=o(
x^{\theta_1-\theta_\infty})$. The second is:  
\be
\rho=\rho(x)\sim-(1+\theta_1-\theta_\infty)K_{12}x^{\theta_1-\theta_\infty+1}
~~~~=
-r x^{\theta_1-\theta_\infty+1}.
\label{rho}
\ee
In both cases $\rho$ is a function of $x$. We are going to prove that
the monodromy of $\Psi_{OUT}$ is independent of $\rho(x)$ (namely, of $x$)
if and only if:  
\be
\rho= r_1~{\theta_1-\theta_\infty+1\over \theta_1}.
\label{rhor1}
\ee
This fact  rules out the first possibility, because  (\ref{maldipancia})
 becomes:
$$
{\rho(x)\over \theta_\infty-\theta_1-1} \sim K_{12}~ x^{\theta_1-\theta_\infty+1} 
+
\hbox{\rm constant}\times\rho(x)~x,
$$
so the last term in the r.h.s. is a higher order correction and $\rho$ is given by (\ref{rho}). Before proving 
(\ref{rhor1}), we complete the matching procedure.
 Using (\ref{rhor1}), we find:
$$
K_0(x)\sim \pmatrix{1 & g ~x^{\theta_1-\theta_\infty+1}
\cr
\cr
0 & 1 }
~\pmatrix{ x^{\theta_1-\theta_\infty\over  2} & 0 
\cr
\cr
0 & 
 x^{\theta_\infty-\theta_1\over  2}
},~~~~ g:=-{r\over \theta_1-\theta_\infty+1}
$$
We are ready to compute $\hat{A_0}=K_0 \hat{\hat{A_0}} {K_0}^{-1}$, 
where $ \hat{\hat{A_0}}
$ is (\ref{hathatA0}), for $\sigma=\theta_1-\theta_\infty$:
\be
\hat{A_0}=\pmatrix{(\hat{\hat{A_0}})_{11}+ g(\hat{\hat{A_0}})_{21} ~x
&
\left({(\hat{\hat{A_0}})_{12}\over x}-2g(\hat{\hat{A_0}})_{11}-  g^2
(\hat{\hat{A_0}})_{21}~x\right)~x^{\theta_1-\theta_\infty+1}
\cr
\cr
\Bigr((\hat{\hat{A_0}})_{21} ~x\Bigl)~{1\over x^{\theta_1-\theta_\infty+1}} 
&
 -(\hat{\hat{A_0}})_{11}- g(\hat{\hat{A_0}})_{21} ~x
}.
\label{hatA0firstorder}
\ee
The first term of each matrix element certainly 
contains no error.

\vskip 0.2 cm
We can now  substitute 
$(\hat{A_0})_{12}= \left(-{r\over
  x}+O(1)\right)~x^{\theta_1-\theta_\infty+1}$ and  
 $\hat{A_1}=- r_1(x)\sim {\theta_1\over \theta_1-\theta_\infty+1}r
x^{\theta_1-\theta_\infty+1}$ into (\ref{leadingtermaprile}), and find:
$$
y(x) \sim
{\theta_1-\theta_\infty+1\over 1-\theta_\infty},~~~~~~~x\to 0.
$$

\vskip 0.5 cm
\noindent
{\it Proof of (\ref{rhor1})}.~~We compute the monodromy of the
solution $\Psi_{OUT}$. At $\lambda=0$ this is given by the matrix
$M_0^{OUT}=\exp\{i\pi(\theta_1-\theta_\infty)\sigma_3\}$, obtained by
the expansion of $\Psi_{OUT}$ at $\lambda=0$ as we did above, with the
choice $C_1= e^{i{\pi\over 2}\theta_1}=1/C_2$.  

Let us study the monodromy at $\lambda=1$. 
We need to  expand  (\ref{psi2}) and (\ref{psi1}) at
$\lambda=1$. First of all, we  
use the contiguity relation:  
{\small
$$
F(1-\theta_1,\theta_\infty-\theta_1,\theta_\infty-\theta_1+1;\lambda)= 
{\theta_1-\theta_\infty\over \theta_\infty-1}~{1\over \lambda}~
\left[
(1-\lambda)^{\theta_1}-F(1-\theta_1,\theta_\infty-\theta_1-1,\theta_\infty-\theta_1;\lambda)
\right].
$$
}
This is used to  rewrite $\psi_1$:
$$
\psi_1=C_1\lambda^{\theta_1-\theta_\infty\over2}(\lambda-1)^{-{\theta_1\over 2}}
+
C_2 e^{i\pi\theta_1}~\lambda^{\theta_\infty-\theta_1\over 2}(\lambda-1)^{-{\theta_1\over 2}}~\left[
{\rho-r_1\over \theta_\infty-1}~{(1-\lambda)^{\theta_1}\over \lambda}+
\right.
$$
$$
\left.
+\left(
{\rho\over \theta_\infty-\theta_1-1}+{r_1-\rho\over \theta_\infty-1}\right)
~{1\over \lambda} ~ F(1-\theta_1,\theta_\infty-\theta_1-1,\theta_\infty-\theta_1;\lambda)
\right]
$$
Then, we substitute in $\psi_1$ the following connection formula:
$$
 F(1-\theta_1,\theta_\infty-\theta_1-1,\theta_\infty-\theta_1;\lambda)
=
{\Gamma(\theta_1)\Gamma(\theta_\infty-\theta_1)\over \Gamma(\theta_\infty-1)}
~F(1-\theta_1,\theta_\infty-\theta_1-1,1-\theta_1;1-\lambda)+
$$
$$
~~~~~~~~~~~~~~~~
+{\Gamma(-\theta_1)\Gamma(\theta_\infty-\theta_1)\over \Gamma(1-\theta_1)\Gamma(\theta_\infty-\theta_1-1)}~(1-\lambda)^{\theta_1}F(1,\theta_\infty-1,1+
\theta_1:1-\lambda)
,~~~~~\theta_1,\theta_\infty\not\in{\bf Z}.
$$
Thus, $\psi_1$ has the following structure, when $\lambda\to 1$:
$$
\psi_1= C_1 e^{-i{\pi\over 2}\theta_1}~\lambda^{\theta_1-\theta_\infty\over 2}(1-\lambda)^{-{\theta_1\over 2}}~+C_2e^{i{\pi \over 2}\theta_1}~\lambda^{\theta_\infty-\theta_1\over 2} 
\left[
(1-\lambda)^{\theta_1\over 2}\sum_{n=0}^{\infty}a_n(1-\lambda)^n +
\right.
$$
$$
\left. +
\left(
{\rho\over \theta_\infty-\theta_1-1}+{r_1-\rho\over \theta_\infty-1}\right)
(1-\lambda)^{-{\theta_1\over 2}}\sum_{n=0}^{\infty} b_n(1-\lambda)^n
\right],
$$
where $a_n, b_n\neq  0$ are the coefficients that follow from 
 the expansion $1/\lambda$ and 
 the hyper-geometric functions at $\lambda=1$. When 
$\lambda^{\pm{\theta_\infty-\theta_1\over 2}}$ is also expanded at
 $\lambda=1$, 
we find that    $\Psi_{OUT}$ has the following structure:
{\small
$$
\Psi_{OUT}=\pmatrix{
(\hbox{series}_1)~(1-\lambda)^{-{\theta_1\over 2}} & 
(\hbox{series}_2)~(1-\lambda)^{{\theta_1\over 2}}+\left(
{\rho\over \theta_\infty-\theta_1-1}+{r_1-\rho\over \theta_\infty-1}\right)
(\hbox{series}_3)~(1-\lambda)^{-{\theta_1\over 2}}
\cr
\cr
0 &(\hbox{series}_4)~(1-\lambda)^{\theta_1\over 2}
}.
$$}
Here ``series'' means a series of the form $\sum_{n=0}^\infty c_n (1-\lambda)^n$, with $c_0\neq 0$. We just give the dominant term, with the choice 
  $C_1=e^{i{\pi\over 2}\theta_1}=1/C_2$:
$$
\Psi_{OUT}=
[I+O(1-\lambda)] \times
$$
$$
\times~\pmatrix{
(1-\lambda)^{-{\theta_1\over 2}}
&
\left[
{\rho-r_1\over\theta_\infty-1}+  (**){\Gamma(-\theta_1)\Gamma(\theta_\infty-
\theta_1)\over \Gamma(1-\theta_1)\Gamma(\theta_\infty-\theta_1-1)}
\right](1-\lambda)^{\theta_1\over 2}
+
(**){\Gamma(\theta_1)\Gamma(\theta_\infty-\theta_1)\over \Gamma(\theta_\infty-1)}(1-\lambda)^{-{\theta_1\over 2}}
\cr
\cr
0
&
(1-\lambda)^{\theta_1\over 2}
},
$$
where:
$$
(**)= {\rho\over \theta_\infty-\theta_1-1}+{r_1-\rho\over \theta_\infty-1}.
$$
Let $v,w$ be non zero arbitrary numbers. We conclude that:
$$
\Psi_{OUT}= \pmatrix{ v & w 
\cr
0
&
 w\left[
{\rho-r_1\over\theta_\infty-1}+  (**){\Gamma(-\theta_1)\Gamma(\theta_\infty-
\theta_1)\over \Gamma(1-\theta_1)\Gamma(\theta_\infty-\theta_1-1)}
\right]^{-1}
} ~[I+O(1-\lambda)]~\times
$$
$$
\times
\pmatrix{ (1-\lambda)^{-{\theta_1\over 2}} & 0 
\cr
0 &  (1-\lambda)^{\theta_1\over 2}
}
\pmatrix{   {1\over v} &
{1\over v} 
(**){\Gamma(\theta_1)\Gamma(\theta_\infty-\theta_1)\over \Gamma(\theta_\infty-1)}
\cr
0
&
{1\over w} \left[
{\rho-r_1\over\theta_\infty-1}+  (**){\Gamma(-\theta_1)\Gamma(\theta_\infty-
\theta_1)\over \Gamma(1-\theta_1)\Gamma(\theta_\infty-\theta_1-1)}
\right]
}
$$
The monodromy matrix is then: 
$$
M_1^{OUT}=
Q^{-1} e^{-i\pi\theta_1\sigma_3} Q,
~~~
Q:=\pmatrix{   {1\over v} &
{1\over v} 
(**){\Gamma(\theta_1)\Gamma(\theta_\infty-\theta_1)\over \Gamma(\theta_\infty-1)}
\cr
\cr
0
&
{1\over w} \left[
{\rho-r_1\over\theta_\infty-1}+  (**){\Gamma(-\theta_1)\Gamma(\theta_\infty-
\theta_1)\over \Gamma(1-\theta_1)\Gamma(\theta_\infty-\theta_1-1)}
\right]
}.
$$
This is independent of $\rho(x)$  if and only if $(**)=0$. This proves
(\ref{rhor1}).  \qed

\vskip 0.3 cm
(\ref{rhor1}) simplifies considerably the structure of 
(\ref{psi1}):
$$
\psi_1= C_1~\lambda^{\theta_1-\theta_\infty\over 2}~(\lambda-1)^{-{\theta_1\over 2}}
-
C_2~{r_1\over \theta_1}\lambda^{{\theta_\infty-\theta_1\over 2}-1}
~(\lambda-1)^{\theta_1\over 2}.
$$
With the choice $C_1=\exp\{i{\pi\over 2}\theta_1\}=1/C_2$, we finally have:
$$
\Psi_{OUT}= \pmatrix{ \lambda^{\theta_1-\theta_\infty\over 2}~
  (1-\lambda)^{-{\theta_1\over 2}} 
&
-{r_1(x)\over \theta_1}~ \lambda^{{\theta_\infty-\theta_1\over 2}-1}~
(1-\lambda)^{\theta_1\over 2} 
\cr
\cr
0
&
 \lambda^{\theta_\infty-\theta_1\over 2}(1-\lambda)^{\theta_1\over 2}
}.
$$

\vskip 0.2 cm
\noindent
{\it Note:} 
The monodromy at $\lambda=0,1$ is diagonal: $
M_0^{OUT}= \exp\left\{i\pi(\theta_1-\theta_\infty)\sigma_3\right\}$, 
$M_1^{OUT}= \exp\left\{-i\pi\theta_1\sigma_3\right\}$. It is
{\it independent of } $r_1$. 
In case {\bf (I)} -- namely, $r_1$ is a non-zero constant and 
the system for $\Psi_1$ is fuchsian --  if $M_0^{OUT}$ is in
diagonal form, then 
$M_1^{OUT}$ is upper triangular and depends on $r_1$.

\subsection{Higher Orders Terms}
We may repeat the same procedure of section  \ref{higherorder}. We
write a Taylor expansion $y(x)=
(\theta_1-\theta_\infty+1)/(1-\theta_\infty)+\sum_{n\geq 1} b_n x^n$,
substitute it into (PVI) and determine recursively the $b_n$'s. Then,
we may 
substitute the result in the matrix elements of $A_i$, $i=0,x,1$,
according to the formulae of \cite{JMU}, and find the higher orders of
the matrix elements. We  just give the result. 
{\small
$$
(A_x)_{11}=~-(A_x)_{22}= ~
{(\theta_\infty-\theta_1)^2+\theta_x^2-\theta_0^2
\over 
4(\theta_1-\theta_\infty)  } ~+~~~~~~~~~~~~~~~~~~~~~~~~~~~~~~~~~~~
$$
$$
+~{\theta_1\over 8}
{[(\theta_1-\theta_\infty)^2-(\theta_0-\theta_x)^2]
[(\theta_1-\theta_\infty)^2-(\theta_0+\theta_x)^2]
\over 
(\theta_1-\theta_\infty)^2[(\theta_1-\theta_\infty)^2-1]} ~x+O(x^2).
$$
\vskip 0.2 cm
$$
(A_x)_{12}=-r
\left\{
{1\over x} -
{\theta_1
[(\theta_1-\theta_\infty+2)(\theta_1-\theta_\infty)+\theta_0^2-\theta_x^2]
\over 
2(\theta_1-\theta_\infty)(\theta_\infty-\theta_1-1)}+O(x)
\right\}
~x^{\theta_1-\theta_\infty+1}.~~~~~~~~~~~~~~~~~
$$
\vskip 0.2 cm
$$
(A_x)_{21}={1\over r}~
{[(\theta_1-\theta_\infty)^2-(\theta_0-\theta_x)^2][(\theta_1-\theta_\infty)^2-(\theta_0+\theta_x)^2]\over 16(\theta_\infty-\theta_1)^2}
\Bigr\{
x~+~~~~~~~~~~
$$
$$
~~~~~~~~\left. -~{\theta_1[(\theta_\infty-\theta_1)(\theta_\infty-\theta_1+2)+\theta_0^2-\theta_x^2]\over 2(\theta_\infty-\theta_1)^3(\theta_\infty-\theta_1+1)}~x^2+O(x^3)
\right\}~
{1 \over x^{\theta_1-\theta_\infty+1}}.
$$
\vskip 0.3 cm
$$
(A_0)_{11}=~-(A_0)_{22}=~
{(\theta_\infty-\theta_1)^2+\theta_0^2-\theta_x^2
\over 
4(\theta_1-\theta_\infty)}~+~~~~~~~~~~~~~~~~~~~~~~~~~~~~~~~~~~~~~~~~~~~~~~~~~~~
$$
$$
-~
{\theta_1\over 8}
{
[(\theta_1-\theta_\infty)^2-(\theta_0-\theta_x)^2][(\theta_1-\theta_\infty)^2-(\theta_0+\theta_x)^2]
\over 
(\theta_1-\theta_\infty)^2[(\theta_1-\theta_\infty)^2-1]} ~x +O(x^2).
$$
\vskip 0.2 cm
$$
(A_0)_{12}=
r
\left\{
{1\over x}  -
{ \theta_1[(\theta_1-\theta_\infty)^2+\theta_0^2-\theta_x^2]
\over 2(\theta_1-\theta_\infty)(\theta_\infty-\theta_1-1)}
+O(x)
\right\}~x^{\theta_1-\theta_\infty+1}.~~~~~~~~~~~~~~~~~~~~~~~
$$
\vskip 0.2 cm
$$
(A_0)_{21}=
-{1\over r}
{[(\theta_1-\theta_\infty)^2-(\theta_0-\theta_x)^2][(\theta_1-\theta_\infty)^2-(\theta_0+\theta_x)^2]
\over 
16(\theta_\infty-\theta_1)^2}~
\Bigr\{~x~+~~~~~~~~~~~~
$$
$$
\left.
-~ {\theta_1[(\theta_1-\theta_\infty)^2+\theta_0^2-\theta_x^2]
\over
2(\theta_\infty-\theta_1)^3(\theta_\infty-\theta_1+1)}~x^2~+O(x^3)
\right\}~{1 \over x^{\theta_1-\theta_\infty+1}}.
$$
\vskip 0.3 cm
$$
(A_1)_{11}=-(A_1)_{22}=
-{\theta_1\over 2} -{\theta_1[(\theta_1-\theta_\infty)^2-(\theta_0-\theta_x)^2]
[(\theta_1-\theta_\infty)^2-(\theta_0+\theta_x)^2]\over 16
[(\theta_\infty-\theta_1)^2-1](\theta_\infty-\theta_1)^2}~x^2+O(x^3).
$$
 \vskip 0.2 cm
$$
(A_1)_{12}=-r{\theta_1\over \theta_\infty-\theta_1-1}\left\{
1+
{(\theta_1+1)~[
(\theta_\infty-\theta_1)(\theta_\infty-\theta_1-2)+\theta_0^2-\theta_x^2]
\over 
2 (\theta_\infty-\theta_1)
(\theta_\infty-\theta_1-2)}~x~+O(x^2)
\right\}~x^{\theta_1-\theta_\infty+1}.
$$
\vskip 0.2 cm
$$
(A_1)_{21}=\left\{{\theta_1[(\theta_1-\theta_\infty)^2-(\theta_0-\theta_x)^2]
[(\theta_1-\theta_\infty)^2-(\theta_0+\theta_x)^2]\over 16 r
(\theta_\infty-\theta_1+1)(\theta_\infty-\theta_1)^2}~x^2 
\right\}
~{1\over x^{\theta_1-\theta_\infty+1}} .
$$
}

\vskip 0.3 cm
\noindent
The above  leading terms of $y(x)$ are related to the above formulae 
through (\ref{leadingtermaprile}). The truncation of  $(A_0)_{12}$ and
$(A_1)_{12}$ above  is enough to reproduce 
 the first two terms of solution (\ref{form1}): 
$$
y(x)= {\theta_1-\theta_\infty+1\over 1-\theta_\infty}
+
{\theta_1 [(\theta_1-\theta_\infty)^2+\theta_x^2-\theta_0^2+2\theta_1
-2\theta_\infty]\over 2 (1-\theta_\infty)(\theta_\infty-\theta_1)(\theta_1-\theta_\infty+2)}~x~+O(x^2).
$$

\subsection{Transpose Case and General Result}
\label{transpgeneral}

 In the above computations, the condition $
x^{\theta_\infty-\theta_1+1}\to 0$  
was necessary to do the matching. We can repeat the matching procedure
starting from  the transpose matrices: 
$$
\hat{A}_1=\pmatrix{-{\theta_1\over 2} & 0 
\cr 
              -r_1 & {\theta_1\over 2}},
~~
\hat{A}_0+\hat{A}_x= \pmatrix{{\theta_1-\theta_\infty\over 2} & 0 
\cr
                                 r_1 &- {\theta_1-\theta_\infty\over 2}
},~~~x~A_x=\pmatrix{0 & 0 \cr \rho(x) & 0}.
$$
 The procedure is exactly the same, with the necessary condition:
$$
x^{\theta_1-\theta_\infty+1}\to 0,~~~ \hbox{ when }x\to 0.
$$
As a result we obtain again  exactly the matrices $A_0$, $A_x$, $A_1$ above. The reader can convince himself without doing any computation, simply looking at 
 the structure of the first terms of 
above matrices. For example, let us have a look at $A_x$. Denote the
constant terms with letters $c_1$,  
$c_2$, ..., etc. We have:
$$
A_x= \pmatrix{  c_1+... & \left\{{c_2\over x}+...\right\}~x^{\theta_1-\theta_\infty+1}
\cr\cr
\Bigl\{c_3~x~+...\Bigr\}~{1\over x^{\theta_1-\theta_\infty+1}} & -c_1+... 
}
=~~~~~~~~~~~~~~~~~~~~~~~~
$$
\vskip 0.2 cm
$$
~~~~~~~~~~~~~~~~~~~~~~~~~~~~=
 \pmatrix{  c_1+... & \Bigl\{c_2~x+...\Bigr\}~{1\over x^{\theta_\infty-
\theta_1+1}}
\cr\cr
\left\{{c_3\over x}~+...\right\}~x^{\theta_\infty-\theta_1+1} & -c_1+... 
}.$$
The role of $\theta_\infty-\theta_1+1$ and $\theta_1-\theta_\infty+1$
is just exchanged.  
By continuity, the matrices $A_0$, $A_x$, $A_1$ computed above hold
 {\it for any value of $\theta_1-\theta_\infty \not \in {\bf Z}$,
   $\theta_\infty\neq 1$}. 

\vskip 0.3 cm
The  matching procedure can be repeated in the same way 
in case $\sigma=-(
\theta_\infty-\theta_1)$ and in the case
$\sigma=\pm(\theta_\infty+\theta_1)$, which yields (\ref{riuffa}). 
For this last cases the results are 
 obtained just by the substitution $\theta_1\mapsto -\theta_1$.

\subsection{Monodromy Data}

We compute the monodromy data for the case
$\sigma=\theta_1-\theta_\infty$, all the other cases being
analogous. In this case, the matching has been realized by:
$$
\Psi_{OUT}(\lambda,x):=\pmatrix{ 
\lambda^{\theta_1-\theta_\infty\over 2}(1-\lambda)^{-{\theta1\over 2}} & 
-{r_1(x)\over \theta_1}  \lambda^{{\theta_\infty-\theta_1\over 2}-1}(1-\lambda)^{{\theta1\over 2}}    
\cr\cr
0 
& 
\lambda^{\theta_\infty-\theta_1\over 2}(1-\lambda)^{{\theta1\over 2}}
}.
$$
$$
\Psi_{IN}(\lambda,x)= K_0(x)\Psi_0\left({\lambda\over x}\right),~~~~~
\Psi_0(\mu)=\left[I+O\left({1\over
    \mu}\right)\right]~\mu^{{\theta_1-\theta_\infty\over
    2}\sigma_3},~~~\mu\to
\infty.
$$ 
Let $\Psi(\lambda)$ denote the solution of the system (\ref{SYSTEM}) of (PVI),
such that: 
$\Psi(\lambda)=[I+O(\lambda^{-1})]~\lambda^{-{\theta_\infty\over
    2}\sigma_3}$, $\lambda\to\infty$, $\theta_\infty\not\in{\bf Z}$. 

\vskip 0.5 cm
\noindent
{\bf MATCHING $\Psi\leftrightarrow \Psi_{OUT}$.}

\vskip 0.2 cm
With the choice $1-\lambda=e^{-i\pi}(\lambda-1)$ ($1-\lambda>0$ for
$0<\lambda<1$) we have:
$$
\Psi_{OUT}= \left[I + O\left({1\over
    \lambda}\right)\right]~\lambda^{-{\theta_\infty\over 2}\sigma_3}~
e^{i{\pi\over 2} \theta_1\sigma_3}
,~~~\lambda\to\infty
$$
Therefore, the correct choice for $\Psi_{OUT}$, which matches with
$\Psi$ is:
$$
\Psi_{OUT}^{Match}:=\Psi_{OUT}
C_{OUT},~~~~~C_{OUT}:=e^{-i{\pi\over 2}\theta_1\sigma_3}.
$$
As a consequence we obtain:
\vskip 0.2 cm
\noindent
{\bf i)}  the monodromy of $\Psi$ at $\lambda=1,\infty$, which is equal to the
monodromy of $\Psi_{OUT}^{Match}$:
$$
  M_1= e^{-i\pi \theta_1\sigma_3},~~~M_\infty=e^{-i\pi \theta_\infty \sigma_3}; 
$$
{\bf ii)} the correct choice for $
\Psi_{IN}^{Match}:= K_0(x)~\Psi_0~ C_{OUT}
$.

\vskip 0.5 cm
\noindent
{\bf MATCHING $\Psi\leftrightarrow \Psi_{IN}^{Match}$} 
\vskip 0.2 cm
This is realized by construction. A consequence, we can compute the
monodromy of $\Psi$ at $\lambda=0,x$. In order to do this, we need the 
local behavior of  $\Psi_{IN}^{Match}(\lambda,x)$  at
$\lambda=0,x$. 
We start with $\Psi_0$, recalling that:
$$
\Psi_0(\mu) = \mu^{\theta_0\over 2}(\mu-1)^{\theta_x\over 2}
\Phi_0(\mu),~~~~~
\Phi_0(\mu)= \pmatrix{ \varphi_1(\mu) & \varphi_2(\mu) \cr \xi_1(\mu) &
  \xi_2(\mu)}.
$$
Here, $\varphi_1$, $\varphi_2$ are two independent solutions of a
Gauss hyper-geometric equation (see (\ref{hypergeom1}) in Appendix 1):
$$
\mu(1-\mu)~ {d^2 \varphi \over d\mu^2} +\bigl(1+c-(a+[b+1]+1)~\mu \bigr)~ {d\varphi\over d\mu}
-a(b+1)~\varphi=0,
$$
where $
a:={\theta_0\over 2}+{\theta_x\over 2} +{\theta_\infty\over
  2}-{\theta_1\over 2}$, 
$b+1:=  {\theta_0\over 2}+{\theta_x\over 2} +{\theta_1\over
  2}-{\theta_\infty\over 2}+1$, 
$c+1:=\theta_0+1$. The functions 
$\xi_1$ and $\xi_2$ are obtained from  $\varphi_1$ and $\varphi_2$ by: 
\be
\xi= {1\over r}\left[
\mu(1-\mu)~{d\varphi\over d\mu} ~-a\left(
\mu+{b-c \over a-b}
\right)~\varphi
\right],
\label{aprilexi}
\ee
In order to have
 generic solutions (i.e. non-logarithmic 
solutions) of the hyper-geometric equation, we must require:
$$
\theta_1-\theta_\infty,~\theta_0,~\theta_x \not \in {\bf Z}.
$$
Then, we have the following sets of independent solutions at
$\mu=0,1,\infty$ respectively (we denote by $F$ the Gauss hyper-geometric
function):
$$
\left\{
\matrix{
\varphi_1^{(0)}=
F\left({\theta_0\over 2}+{\theta_x\over 2} +{\theta_\infty\over
  2}-{\theta_1\over 2},~
{\theta_0\over 2}+{\theta_x\over 2} +{\theta_1\over
  2}-{\theta_\infty\over 2}+1, ~1+\theta_0;~\mu 
\right),
\cr
\cr
\varphi_2^{(0)}=\mu^{-\theta_0}~
F\left(
{\theta_x\over 2}-{\theta_0\over 2} +{\theta_\infty\over
  2}-{\theta_1\over 2},~
{\theta_x\over 2}-{\theta_0\over 2} +{\theta_1\over
  2}-{\theta_\infty\over 2}+1,~ 1-\theta_0;~\mu 
\right).
}
\right. 
$$

$$
\left\{\matrix{
\varphi_1^{(1)}=
F\left({\theta_0\over 2}+{\theta_x\over 2} +{\theta_\infty\over
  2}-{\theta_1\over 2},
~{\theta_0\over 2}+{\theta_x\over 2} +{\theta_1\over
  2}-{\theta_\infty\over 2}+1,~ 1+\theta_x;~1-\mu 
\right),
\cr
\cr
\varphi_2^{(1)}=
(1-\mu)^{-\theta_x}~
F\left({\theta_0\over 2}-{\theta_x\over 2} +{\theta_\infty\over
  2}-{\theta_1\over 2},~
{\theta_0\over 2}-{\theta_x\over 2} +{\theta_1\over
  2}-{\theta_\infty\over 2}+1,~ 1-\theta_x;~1-\mu 
\right).
}
\right. 
$$

$$
\left\{
\matrix{
\varphi_1^{(\infty)}=
\mu^{-{\theta_0\over 2}-{\theta_x\over 2} -{\theta_\infty\over
  2}+{\theta_1\over 2}}  ~
F\left({\theta_0\over 2}+{\theta_x\over 2} +{\theta_\infty\over
  2}-{\theta_1\over 2},~
{\theta_x\over 2}-{\theta_0\over 2} +{\theta_\infty\over
  2}-{\theta_1\over 2}, ~\theta_\infty-\theta_1;~{1\over\mu} 
\right),
\cr
\cr
\varphi_2^{(\infty)}=
\mu^{-{\theta_0\over 2}-{\theta_x\over 2} -{\theta_1\over
  2}+{\theta_\infty\over 2}-1}~
F\left({\theta_0\over 2}+{\theta_x\over 2} +{\theta_1\over
  2}-{\theta_\infty\over 2}+1,~
{\theta_x\over 2}-{\theta_0\over 2} +{\theta_1\over
  2}-{\theta_\infty\over 2}+1, ~2+\theta_1-\theta_\infty;~{1\over \mu} 
\right).
}
\right. 
$$

The connection formulae can be found in any book on special functions:
$$
\left[\varphi_1^{(0)},\varphi_2^{(0)}\right]
=
\left[\varphi_1^{(1)},\varphi_2^{(1)}\right]C_{01},~~~
-\pi<\hbox{arg}~(1-\mu)<\pi .
$$
$$
\left[\varphi_1^{(0)},\varphi_2^{(0)}\right]
=
\left[\varphi_1^{(\infty)},\varphi_2^{(\infty)}\right]C_{0\infty},~~~
0<\hbox{arg}~\mu<2\pi.
$$
Here the connection matrices $C_{01}$, $C_{0\infty}$ are 
(\ref{C01pasqua}) and 
(\ref{C0inftypasqua}) respectively. 

\vskip 0.3 cm
 From the Taylor expansion of the hyper-geometric functions
in $\varphi_i^{(\infty)}$ and  (\ref{aprilexi}) we compute:
$$
\mu^{\theta_0\over 2}
(\mu-1)^{\theta_x\over 2}
\varphi_1^{(\infty)}=\left[1+O\left({1\over
  \mu}\right)\right]\mu^{\theta_1-\theta_\infty \over 2},~~~
\mu^{\theta_0\over 2}
(\mu-1)^{\theta_x\over 2}
\varphi_2^{(\infty)}=~{1\over \mu}~\left[1+O\left({1\over
  \mu}\right)\right]\mu^{\theta_\infty-\theta_1 \over 2},
$$
$$
\mu^{\theta_0\over 2}
(\mu-1)^{\theta_x\over 2}
\xi_1^{(\infty)}=\left[1+O\left({1\over
  \mu}\right)\right]\mu^{\theta_1-\theta_\infty \over 2},~~~
\mu^{\theta_0\over 2}
(\mu-1)^{\theta_x\over 2}
\xi_2^{(\infty)}=~{\theta_1-\theta_\infty+1\over r_0}~\left[1+O\left({1\over
  \mu}\right)\right]\mu^{\theta_\infty-\theta_1 \over 2},
$$
It follows that the matrix $\Psi_0(\mu)= \mu^{\theta_0\over 2}(\mu-1)^{\theta_x\over
  2}\Phi_0(\mu)$ with the prescribed behavior  $\left[I+O\left({1\over
  \mu}\right)\right]\mu^{{\theta_1-\theta_\infty \over 2}\sigma_3}$ 
at $\mu=\infty$ is: 
$$
\Psi_0(\mu)=~ \mu^{\theta_0\over 2}(\mu-1)^{\theta_x\over
  2} ~
\pmatrix{\varphi_1^{(\infty)} & \varphi_2^{(\infty)}
\cr
\xi_1^{(\infty)} & \xi_2^{(\infty)}
}
~\pmatrix{1 & 0 \cr 0 & {r\over \theta_1-\theta_\infty+1}
}
$$
Let:
$$
C:=\pmatrix{1 & 0 \cr 0 & {r\over \theta_1-\theta_\infty+1}
}C_{OUT}
\equiv
\pmatrix{1 & 0 \cr 0 & {r\over \theta_1-\theta_\infty+1}
}e^{-i{\pi\over 2}\theta_1\sigma_3}.
$$
We conclude that:
$$
\Psi_{IN}^{Match}(\lambda,x)=K_0(x)~ 
\left({\lambda\over x}\right)^{\theta_0\over 2}\left({\lambda\over x}-1\right)^{\theta_x\over
  2} ~
\pmatrix{\varphi_1^{(\infty)}\left({\lambda\over x}\right) & \varphi_2^{(\infty)}\left({\lambda\over x}\right)
\cr
\cr
\xi_1^{(\infty)}\left({\lambda\over x}\right) & \xi_2^{(\infty)}\left({\lambda\over x}\right)
}C,~~~~~~~~~~~~~~~~~~~~~~~~~
$$
$$
=K_0(x)~ 
\left({\lambda\over x}\right)^{\theta_0\over 2}\left({\lambda\over x}-1\right)^{\theta_x\over
  2} ~
\pmatrix{\varphi_1^{(0)}\left({\lambda\over x}\right) & \varphi_2^{(0)}\left({\lambda\over x}\right)
\cr
\cr
\xi_1^{(0)}\left({\lambda\over x}\right) & \xi_2^{(0)}\left({\lambda\over x}\right)
}C_{0\infty}^{-1}C,
$$
$$
=K_0(x)~ 
\left({\lambda\over x}\right)^{\theta_0\over 2}\left({\lambda\over x}-1\right)^{\theta_x\over
  2} ~
\pmatrix{\varphi_1^{(1)}\left({\lambda\over x}\right) & \varphi_2^{(1)}\left({\lambda\over x}\right)
\cr
\cr
\xi_1^{(1)}\left({\lambda\over x}\right) & \xi_2^{(1)}\left({\lambda\over x}\right)
}C_{01}C_{0\infty}^{-1}C.
$$
The behaviors of the above matrix at $\lambda=x,~0$ is easily computed
from:
$$
\mu^{\theta_0\over 2}(\mu-1)^{\theta_x\over 2}
\varphi_1^{(1)}~=~~~~~~~~~~~~~~~~~~~~~~~~~~~~~~~~~~~~~~~~~~~~~~~~~~~~~~~~~~~~~~~~~~~~~~~~~~~~~~~~
$$
$$
~~~~~~~~~~=\mu^{\theta_0\over 2}(\mu-1)^{\theta_x\over 2} F(...;1-\mu)=~
x^{-{\theta_0\over 2}-{\theta_x\over 2}}F\left(...;{x-\lambda\over
  x}\right) ~\lambda^{\theta_0\over 2}(\lambda-x)^{\theta_x\over 2}
,
$$
$$
\mu^{\theta_0\over 2}(\mu-1)^{\theta_x\over 2}
\varphi_2^{(1)}~=~~~~~~~~~~~~~~~~~~~~~~~~~~~~~~~~~~~~~~~~~~~~~~~~~~~~~~~~~~~~~~~~~~~~~~~~~~~~~~~~
$$
$$
~~~~~~~=\mu^{\theta_0\over 2}(\mu-1)^{\theta_x\over 2}(1-\mu)^{-{\theta_x\over 2}}
=~(-1)^{\theta_x}
x^{-{\theta_0\over 2}+{\theta_x\over 2}}F\left(...;{x-\lambda\over
  x}\right) ~\lambda^{\theta_0\over 2}(\lambda-x)^{-{\theta_x\over 2}}
;
$$
and:
$$
\mu^{\theta_0\over 2}(\mu-1)^{\theta_x\over 2}
\varphi_1^{(0)}~=~~~~~~~~~~~~~~~~~~~~~~~~~~~~~~~~~~~~~~~~~~~~~~~~~~~~~~~~~~~~~~~~~~~~~~~~~~~~~~~~
$$
$$=
\mu^{\theta_0\over 2}(\mu-1)^{\theta_x\over 2} F(...;\mu)= 
x^{-{\theta_0\over 2}-{\theta_x\over 2}}~
F\left(...;{\lambda\over x}\right)(\lambda-x)^{\theta_x\over 2} \lambda^{\theta_0\over 2},
$$
$$
\mu^{\theta_0\over 2}(\mu-1)^{\theta_x\over 2}
\varphi_2^{(0)}~=~~~~~~~~~~~~~~~~~~~~~~~~~~~~~~~~~~~~~~~~~~~~~~~~~~~~~~~~~~~~~~~~~~~~~~~~~~~~~~~~
$$
$$
=\mu^{-{\theta_0\over 2}}(\mu-1)^{\theta_x\over 2} F(...;\mu)= 
x^{{\theta_0\over 2}-{\theta_x\over 2}}~
F\left(...;{\lambda\over x}\right) (\lambda-x)^{\theta_x\over 2}\lambda^{-{\theta_0\over 2}}.
$$
From these we compute:
$$
K_0(x)~ 
\mu^{\theta_0\over 2}
(\mu-1)^{\theta_x\over 2} ~
\pmatrix{\varphi_1^{(0)} 
& \varphi_2^{(0)}
\cr
\xi_1^{(0)} &
\xi_2^{(0)}
}=
\psi_0^{IN}(x)(I+O(\lambda))~\lambda^{{\theta_0\over
    2}\sigma_3},~~~\lambda\to 0
$$
$$
K_0(x)~ 
\mu^{\theta_0\over 2}(\mu-1)^{\theta_x\over
  2} ~
\pmatrix{\varphi_1^{(1)} & \varphi_2^{(1)}
\cr
\xi_1^{(1)}& \xi_2^{(1)}
}=\psi_x^{IN}(x)(I+O(\lambda-x))~(\lambda-x)^{{\theta_x\over
    2}\sigma_3},~~~\lambda\to x.
$$
Here, we don't need to explicitly give the invertible matrices 
$\psi_0^{IN}(x)$ and $\psi_x^{IN}(x)$. 
From the above procedure, we find:
$$
M_0=C^{-1}~\Bigl(C_{0\infty} ~\exp\{i\pi\theta_0\sigma_3\}~
C_{0\infty}^{-1}\bigr) ~C,
~~~
M_x=C^{-1}\Bigl[C_{0\infty}~\Bigl(~ C_{01}^{-1}~ 
\exp\{i\pi\theta_x\sigma_3\}
 ~ C_{01}\Bigr)~
C_{0\infty}^{-1}\Bigr]~C.
$$
We finally observe that $M_1$ and $M_\infty$ are diagonal, so they are
invariant for the conjugation $M\mapsto CMC^{-1}$. With this in mind,
we get the result of theorem \ref{thMONODROMY}, point a).

We stress that 
 $\hat{A}_0+\hat{A}_x$ and $\hat{A}_1$ 
are upper (or lower) triangular matrices, and {\it the group generated
 by $M_xM_0$ and $M_1$ is reducible}.  

%%%%%%%%%%%%%%%%%%%%%%%%%%%%%%%%%%

%%%%%%%%%%%%%%
%%%%%%%%%
%%%%%%%       CANCELED IN THE PREPRINT
%%%%%%%%%
%%%%%%%%%%%%%%
%
%In section \ref{semireducible} we saw that if 
% the conditions  $\sigma=\pm(\theta_\infty-\theta_1)$ or $\sigma=\pm
% (\theta_\infty+\theta_1$ occur, then $\hat{A}_0+\hat{A}_x$ and
% $\hat{A}_1$  
%are upper (or lower) triangular matrices, and {\it the group
% generated  by $M_xM_0$ and $M_1$ is reducible}.  

%The matrices (\ref{hathatA0}) and (\ref{hathatAx}) are 
 %upper triangular (or lower triangular, by a redefinition of $r_0$) iff
% $\sigma= \pm(\theta_0-\theta_x)$, $\sigma=\pm (\theta_0+\theta_x)$. 
%Theferore, a necessary condition for the  system (\ref{SYSTEM}) to have
 %a reducible monodromy group is that one of the following cases occur:  
% $\theta_1-\theta_\infty=\pm(\theta_0+\theta_x)$,
% $\theta_1-\theta_\infty=\pm(\theta_0-\theta_x)$,
% $\theta_1+\theta_\infty=\pm(\theta_0+\theta_x)$,
% $\theta_1+\theta_\infty=\pm(\theta_0-\theta_x)$.  
%The matching procedure in the complete reducible case 
%works exactly as in section \ref{semireducible}, and the expansions
%(\ref{solpiu}) (\ref{solmeno}) hold with the specific values of the
%$\theta$'s. See section \ref{reduciblePASQUA} for the reducible
%monodromy case. 
%%%%%%%%%%%%%%%%%%%%

\section{Case $\theta_\infty=1$, $\theta_1=0$, $\sigma=\pm
  1$.~Solution (\ref{form3}).}
\label{supersemireducible}

 In section \ref{semireducible} we imposed that
 $\theta_\infty\neq 1$ and $\theta_\infty-\theta_1\not\in{\bf Z}$. Here we consider 
the case  $\theta_\infty= 1$, $\theta_\infty-\theta_1=1$. We have:
$$
\hat{A_0}+\hat{A_x}=\pmatrix{-{1\over 2} & r_1 \cr 
                              0 & {1\over 2} }
, ~~~~~
\hat{A_1}= \pmatrix{ 0 & -r_1 \cr 0 & 0 }.
$$
Also the transpose matrices are possible. 
 For $|\lambda|\leq |x|^{\delta_{IN}}$, we use the
 reduction  (\ref{system0}) (namely, (\ref{fuchsianSYSTEMIN})).  
 The matrices $\hat{\hat{A_0}}$, $\hat{\hat{A_x}}$ are given by
 (\ref{hathatA0}) and (\ref{hathatAx}) with the substitution
 $\sigma=1$ or $-1$. For definiteness, let us take $\sigma=1$ in the
 following.    

\vskip 0.2 cm 
 For $|\lambda|\geq |x|^{\delta_{OUT}}$, we approximate (\ref{SYSTEM}) with: $
{d\Psi_{OUT}\over d\lambda} = \left[
{xA_x\over \lambda^2}+ {A_0+A_x\over \lambda} 
+{A_1\over \lambda-1}
\right]\Psi_{OUT}
$. 
For definiteness, let us consider the case when $\hat{A}_1$ and $\hat{A_0}+\hat{A_x}$ are upper triangular. Again, we make the  hypothesis that the leading terms in $xA_x$ define an  upper triangular matrix:
$$
 xA_x=: \pmatrix{ 0 & \rho \cr 0 & 0 }~+\hbox{ higher orders}.
$$
Therefore, we will study: 
$$
{d\Psi_{OUT}\over d\lambda} = \left[
{1\over \lambda^2} \pmatrix{ 0 & \rho \cr 0 & 0 }+ {1\over \lambda}
\pmatrix{-{1\over 2} & r_1 \cr  
                              0 & {1\over 2} }
+{1\over \lambda-1} \pmatrix{ 0 & -r_1 \cr 0 & 0 }
\right]\Psi_{OUT}
$$
This is a reducible system. To solve it, we write $\Psi_{OUT}$ in
vector notation  
$\Psi_{OUT}= \pmatrix{\psi_1\cr\psi_2}$. The system becomes:
$$
{d\psi_1\over d\lambda}= -{1\over 2\lambda} \psi_1 +\left(
{\rho\over\lambda^2}+{r_1\over \lambda}-{r_1\over \lambda-1}
\right)\psi_2,~~~~~~{d\psi_2\over d\lambda} = {1\over 2\lambda} \psi_2.
$$
The solution obtained by variation of parameters is:
$$
\psi_1= C_1\lambda^{-{1\over 2}}+C_2 \lambda^{-{1\over 2}}(\rho\ln\lambda -
r_1 \ln(\lambda-1)),~~~~~\psi_2=C_2\lambda^{1\over 2},~~~~~C_1,C_2\in {\bf C}.
$$
We can choose the following fundamental matrix:
$$
{\small
\Psi_{OUT}=\pmatrix{\lambda^{-{1\over 2}} & \lambda^{-{1\over 2}} \bigl(
\rho \ln\lambda
-r_1 \ln(\lambda-1)\bigr)
\cr\cr
0 & \lambda^{1\over 2}
}~=
\pmatrix{\lambda^{-{1\over 2}} & 0 \cr 0 &  \lambda^{1\over 2}}
\pmatrix{1 & \rho~\ln \lambda \cr 0 & 1 }
\pmatrix{ 1 & -r_1 \ln(\lambda-1) \cr 0 & 1 }. 
}
$$
Its monodromy  relative to the loops $\lambda\mapsto \lambda
e^{2\pi i}$ and $(\lambda-1)\mapsto (\lambda-1) e^{2\pi i}$ is respectively:
$$\Psi_{OUT}\mapsto\Psi_{OUT} \pmatrix{-1
  & -2\pi i \rho  
\cr
0 & -1}.~~~~~
\Psi_{OUT}\mapsto\Psi_{OUT}
\pmatrix{1 & -2\pi i r_1 
\cr
0 & 1}.
$$
Therefore, $\rho$ and $r_1$ {\it must be independent of} (the monodromy preserving deformation) $x$. 
 We observe that any fundamental matrix solution of the form:
$$
 \Psi_{OUT}\pmatrix{ 1 & f(x) \cr 0 & 1},
$$ 
has the same monodromy of $\Psi_{OUT}$, 
for any arbitrary function of $x$. This fact will be used soon, with
the choice $f(x)= -\rho \ln x$.  

\vskip 0.2 cm
 With the choice of the branch
of $\ln(\lambda-1)= \ln(1-\lambda)+i\pi $
(i.e. $(\lambda-1)=e^{i\pi}(1-\lambda)$, $1-\lambda>0$ for
$0<\lambda<1$), it is convenient to redefine $\Psi_{OUT}$ by: 
$$
\Psi_{OUT}:= \pmatrix{\lambda^{-{1\over 2}} & 0 \cr 0 &  \lambda^{1\over 2}}
\pmatrix{1 & \rho~\ln \lambda \cr 0 & 1 }
\pmatrix{ 1 & -r_1 \ln(\lambda-1) \cr 0 & 1 }
\pmatrix{1 & i\pi r_1 \cr 0 & 1}~= 
$$
$$
=~\pmatrix{1 & -{r_1\over \lambda}\ln(1-\lambda) \cr
           0 & 1}
\pmatrix{\lambda^{-{1\over 2}} & 0 \cr 0 & \lambda^{1\over 2}}
\pmatrix{ 1 & \rho ~\ln\lambda \cr 0 & 1}.
$$ 
Therefore:
$$
\Psi_{OUT}= \left[
\pmatrix{1 & r_1 \cr 0 & 1} + \pmatrix{  0 & i\pi r_1 \cr 0 & 0}
\sum_{n=1}^\infty{\lambda^n\over n+1}
\right]
~\pmatrix{\lambda^{-{1\over 2}} & 0 \cr 0 & \lambda^{1\over 2}}
\pmatrix{ 1 & \rho~ \ln\lambda \cr 0 & 1},~~~\lambda\to 0
$$

\subsection{Matching}

The solution $\Psi_0\left({\lambda\over x}\right)$ has been 
introduced in section \ref{naive}. 
$$
\Psi_0\left({\lambda\over x}\right)= \left[
I+O\left({x\over \lambda}\right)
\right]
\pmatrix{
\left({\lambda\over x}\right)^{1\over 2} & 0 \cr 0 & 
 \left({\lambda\over x}\right)^{-{1\over 2}}
}
\pmatrix{1 & 0 \cr 
         R~\ln\left({\lambda\over x}\right) & 1},
~~~R:=(\hat{\hat{A_x}})_{21}.
$$
 Some adjustments are necessary. 
Let us consider the permutation matrix $P:=\pmatrix{0 & 1 \cr 1 & 0}$ and 
redefine: 
$$\Psi_0 \mapsto \Psi_0 ~P \pmatrix{1 & 0 \cr 0 & R^{-1}},
~~~~~~~\Psi_{OUT}\mapsto \Psi_{OUT} \pmatrix{1 & -\rho \ln x \cr 
0 & 1} \pmatrix{1 & 0 \cr 0 & \rho^{-1}}.
$$
%%%%%%%%%%%%%
%Namely, we choose the fundamental solution: 
%$$
%\Psi_0= \left[
%I+O\left({x\over \lambda}\right)
%\right]
%\pmatrix{
%\left({\lambda\over x}\right)^{1\over 2} & 0 \cr 0 & 
% \left({\lambda\over x}\right)^{-{1\over 2}}
%}
%\pmatrix{1 & 0 \cr 
%         R~\ln{\lambda\over x} & 1}~P
%~=
%$$
%$$
%~~~~~~~~~~~~~=~ P~\left[
%I+O\left({x\over \lambda}\right)
%\right]
%\pmatrix{
%\left({\lambda\over x}\right)^{-{1\over 2}} & 0 \cr 0 & 
% \left({\lambda\over x}\right)^{1\over 2}}
%\pmatrix{
%1 & R~ \ln {\lambda \over x} 
%\cr
%0 & 1}.
%$$
%We also redefine $\Psi_{OUT}\mapsto \Psi_{OUT} \pmatrix{1 & -\rho \ln x \cr 
%0 & 1}$. Namely, we choose the fundamental matrix:
%$$
%\Psi_{OUT}= ~\pmatrix{1 & -{r_1\over \lambda}\ln(1-\lambda) \cr
%           0 & 1}
%\pmatrix{\lambda^{-{1\over 2}} & 0 \cr 0 & \lambda^{1\over 2}}
%\pmatrix{ 1 & \rho ~\ln{\lambda\over x} \cr 0 & 1}~=
%$$
%$$
%=~ \left[
%\pmatrix{1 & r_1 \cr 0 & 1} + \pmatrix{  0 & i\pi r_1 \cr 0 & 0}
%\sum_{n=1}^\infty{\lambda^n\over n+1}
%\right]
%~\pmatrix{\lambda^{-{1\over 2}} & 0 \cr 0 & \lambda^{1\over 2}}
%\pmatrix{ 1 & \rho~ \ln{\lambda\over x} \cr 0 & 1},~~~\lambda\to 0.
%$$
%%
%
%
%
% 
 As we have already observed, this re-definition does not affect the monodromy 
of $\Psi_{OUT}$, which is independent of $x$. 
%
%
%
%
%
%Now, write: 
%{\small
%$$
%\pmatrix{1 & \rho \ln {\lambda\over x} \cr  0 & 1}=\pmatrix{ 1 & 0 \cr 0 &
%{1\over \rho}}\pmatrix{1 &  \ln {\lambda\over x} \cr  0 & 1}\pmatrix{
% 1 & 0 \cr 0 & 
%\rho},
%~~~~~
%\pmatrix{ 1 & R \ln{\lambda \over x} \cr 0 & 1}=\pmatrix{ 1 & 0 \cr 0 &
%{1\over R}}\pmatrix{1 &  \ln {\lambda\over x} \cr  0 & 1}\pmatrix{ 1
% & 0 \cr 0 & 
%R}.
%$$
%}
%Then
%%%%%%%%%%%%%%%
%%
%
%
%
% Re-define again:
%$$
%\Psi_{OUT}\mapsto \Psi_{OUT} \pmatrix{ 1 & 0 \cr 0 &
%{1\over \rho}},
%~~~~~
%\Psi_0\mapsto \Psi_0  \pmatrix{ 1 & 0 \cr 0 &
%{1\over R}}.
%$$
The matching relation $\Psi_{OUT}(\lambda)\sim K_0(x)\Psi_0
\left({\lambda/ x}\right)$ becomes: 
$$
\pmatrix{ 1 & r_1 \cr 0 & 1} \pmatrix{ 1 & 0 \cr 0 &
{1\over \rho}}
\pmatrix{ \lambda^{-{1\over 2}} & 0 \cr 0 & \lambda^{1\over 2}}
\pmatrix{ 1 &  ~\ln{\lambda\over x} \cr 0 & 1}
\sim~~~~~~~~~~~~~~~~~~~~~~~~~~~~~~~~~~~~~~~~
$$
$$~~~~~~~~~~~~~~~~~~~~~~~\sim
K_0(x)~P\pmatrix{ 1 & 0 \cr 0 &
{1\over R}}\pmatrix{x^{1\over 2} & 0 \cr 0 & x^{-{1\over 2}}} 
\pmatrix{ \lambda^{-{1\over 2}} & 0 \cr 0 & \lambda^{1\over 2}}
\pmatrix{ 1 &  \ln{\lambda \over x} \cr 0 & 1}.
$$
The matching is thus realized, with the choice: 
$$
K_0(x)\sim \pmatrix{1 & {r_1\over \rho} \cr 0 & {1\over \rho}}
\pmatrix{x^{-{1\over 2}} & 0 \cr 0 & R~ x^{1\over 2}}P
~
\equiv \pmatrix{{r_1\over \rho} & 1 \cr {1\over \rho} & 0} \pmatrix{ 
R~x^{1\over 2} & 0 \cr 0 & x^{-{1\over 2}}
}.
$$
It follows that:
$$\hat{A_0}= \pmatrix{{\theta_x^2-\theta_0^2-1\over 4} +{r_1\over \rho}
{[1-(\theta_0-\theta_x)^2][1-(\theta_0+\theta_x)^2]\over 16}
x
&
-{\rho\over x} +r_1{\theta_0^2-\theta_x^2+1\over 2} -{r_1^2\over \rho}
 {[1-(\theta_0-\theta_x)^2)][1-(\theta_0+\theta_x)^2]\over 16}x
\cr\cr
{1\over \rho}{[1-(\theta_0-\theta_x)^2)][1-(\theta_0+\theta_x)^2]\over 16}x
&
-{\theta_x^2-\theta_0^2-1\over 4} -{r_1\over \rho}
{[1-(\theta_0-\theta_x)^2][1-(\theta_0+\theta_x)^2]\over 16}
x
},
$$
$$
\hat{A_x}=\pmatrix{
{\theta_0^2-\theta_x^2-1\over 4} -{r_1\over \rho}{[1-(\theta_0-\theta_x)^2][1-(\theta_0+\theta_x)^2]\over 16} x
&
{\rho\over x} +r_1{\theta_x^2-\theta_0^2+1\over 2} +{r_1^2\over \rho}
{[1-(\theta_0-\theta_x)^2][1-(\theta_0+\theta_x)^2]\over 16}x
\cr\cr
-{1\over \rho}{[1-(\theta_0-\theta_x)^2][1-(\theta_0+\theta_x)^2]\over 16}x
&
 -{\theta_0^2-\theta_x^2-1\over 4} +{r_1\over \rho}{[1-(\theta_0-\theta_x)^2][1-(\theta_0+\theta_x)^2]\over 16} x
}.
$$
The leading term of each matrix element certainly contains no error. 
Note that $r$ has simplified. 
Actually, the constant $\rho$ plays the role of $r$. 

\vskip 0.2 cm
We will not repeat again the discussion for the higher order
terms. The matrix elements are Taylor expansions, corresponding to a
Taylor expanded $y(x)$, the convergence of which is proved 
by a Briot-Bouquet argument. 
The first two leading terms of $A_0$ and $A_x$ above are actually correct. 
In particular, we need:  
$$
(A_0)_{12}=   -{\rho\over x} +r_1{\theta_0^2-\theta_x^2+1\over 2} +O(x).
$$
The leading terms of $A_1$  are:
$$
A_1= 
\pmatrix{{r_1\over \rho}{[1-(\theta_0-\theta_x)^2][1-(\theta_0+\theta_x)^2]\over 32} x^2+ O(x^3) 
&
 -r_1 +r_1 {\theta_0^2-\theta_x^2-1\over 2}x+O(x^2)
\cr\cr
{r_1\over \rho^2} {[1-(\theta_0-\theta_x)^2]^2[1-(\theta_0+\theta_x)^2]^2
\over 1024}x^4+O(x^5)
&
-{r_1\over \rho}{[1-(\theta_0-\theta_x)^2][1-(\theta_0+\theta_x)^2]\over 32} x^2+ O(x^3) 
}.
$$
 The above truncations of  $
(A_1)_{12}$ and $(A_0)_{12}$  corresponding to the first two terms of
 the Taylor expansion of the solution (\ref{form3}), 
through (\ref{leadingtermaprile}):
 \be
y(x)~=~a ~+~{1-a\over 2} (1+\theta_0^2-\theta_x^2)~x~+~O(x^2),
~~~~~~~\hbox{ where } a:=\left(1-{r_1\over \rho}\right)^{-1}.
\label{aprileA}
\ee
Observe that $a$ depends on the monodromy datum ${r_1\over
  \rho}$.

\subsection{Monodromy}

The matching $\Psi_{IN}\leftrightarrow\Psi_{OUT}$ has been realized by:
$$
\Psi_{OUT}= \pmatrix{\lambda^{-{1\over 2}}& 0 \cr 0 & \lambda^{1\over
    2} } 
\pmatrix{1 & \rho ~\ln\lambda \cr 0 & 1} 
\pmatrix{ 1 & -r_1~\ln(\lambda-1)
\cr 0 & 1}
\pmatrix{1 & -\ln x ~+i\pi{r_1\over \rho} \cr 0 & {1\over \rho}},
$$
and:
$$
\Psi_{IN}(\mu)=K_0(x) \Psi_0(\mu) \pmatrix{0& R^{-1} \cr 1 &
  0},~~~~~
R= {[(\theta_0+\theta_x)^2-1][(\theta_0-\theta_x)^2-1]\over 16 r}.
$$
\vskip 0.2 cm
$$
\Psi_0(\mu)= \left[I+O\left({1\over \mu} \right)\right]
\pmatrix{ \mu^{1\over 2} & 0 \cr 
0 & \mu^{-{1\over 2}}}
\pmatrix{ 1 & 0 \cr R ~\ln\mu & 1},
~~~~~~~\mu\to\infty.
$$
Namely:
$$
\Psi_{IN}\sim K_0(x) P
\pmatrix{
x^{1\over 2} & 0 \cr 0 & R^{-1}~ x^{-{1\over 2}}
}
~\pmatrix{\lambda^{-{1\over 2}} & 0 \cr 0 & \lambda^{1\over 2}}
\pmatrix{ 1 & \ln{\lambda \over x} \cr 
0 & 
1 
},~~~~~~~~~~~~~~~~~~~~~~~~~~~~~~~~
$$
$$
~~~~~~~~~~~~~~~~~~~~~~~~~~~~~~~~~~~~~
= \pmatrix{1 & {r_1\over \rho} \cr 
                0 & {1\over \rho}}
~
~\pmatrix{\lambda^{-{1\over 2}} & 0 \cr 0 & \lambda^{1\over 2}}
\pmatrix{ 1 & \ln{\lambda \over x} \cr 
0 & 
1 
},
$$
for $\mu=\lambda/x\to\infty$, $|\lambda|<|x|^{\delta_{IN}}$,
$x\to 0$.

\vskip 0.5 cm
\noindent
{\bf MATCHING $\Psi\leftrightarrow \Psi_{OUT}$.}
\vskip 0.2 cm
\noindent
The correct choice of $\Psi_{OUT}$ must match with: 
$$
\Psi= \left[I+O\left({1\over \lambda}
\right)\right]\lambda^{-{\theta_\infty\over
    2}\sigma_3}~\lambda^L,~~~\lambda\to\infty;~~~~~\theta_\infty=1; 
$$
where: $L=\pmatrix{0 &\rho-r_1 \cr 0 & 0}$.  This form of $L$ follows from 
the standard theory of Fuchsian systems, and from the expansion of 
the system (\ref{SYSTEM}) at $\nu:={1\over\lambda}\to 0$:
$$
{d\Psi\over d\nu} = \left[{\sigma_3\over 2}
{1\over \nu}~-~(A_1+xA_x)+O(\nu)\right]\Psi,~~~~~\nu\to \infty;
$$
Thus $L_{12}=-(A_1+xA_x)_{1,2}|_{x=0}\equiv r_1-\rho$~
(this is computed from
the expansions of $A_1$ and $A_x$ at $x=0$ obtained before).   
We expand $\Psi_{OUT}$ at $\lambda=\infty$. We easily get: 
$$
\Psi_{OUT}= [I+O(\lambda^{-1})]~\lambda^{-{1\over 2}\sigma_3}
\pmatrix{ 1 & (\rho-r_1)~\ln\lambda \cr 0 & 1} 
~
\pmatrix{ 1 & i\pi{r_1\over \rho}-\ln x \cr
0
&
{1\over \rho}
}.
$$
Therefore, the correct choice is: 
$$
\Psi_{OUT}^{Match}:= \Psi_{OUT}
\pmatrix{ 1 & i\pi{r_1\over \rho}-\ln x \cr
0
&
{1\over \rho}
}^{-1}
=
 \lambda^{-{1\over 2}\sigma_3} \pmatrix{ 1 & \rho \ln\lambda \cr
0 & 1 } 
\pmatrix{ 1 & - r_1\ln(\lambda -1) \cr 0 & 1
}.
$$
As a consequence, we obtain:

\vskip 0.2 cm
\noindent
{\bf i)} the monodromy of $\Psi$ at  $\lambda=1,\infty$. This coincides with
that of $ \Psi_{OUT}^{Match}$, which is easily computed from the local
behavior: 
$$
\Psi_{OUT}^{Match}(\lambda)= [I+O(\lambda-1)] 
\pmatrix{1 & -r_1 \ln(\lambda-1) \cr 0 & 1 },~~~~~\lambda\to 1.
$$
Thus, for $(\lambda-1)\mapsto (\lambda-1) \exp\{2\pi i\}$ and
$\lambda \mapsto \lambda  \exp\{2\pi i\}$ the
monodromy is:
$$
M_1= \pmatrix{1 & -2\pi i r_1 \cr 0 & 1}, ~~~~~M_{\infty}= \pmatrix{
-1 & 2\pi i (r_1-\rho) \cr 0 & -1}.
$$ 

\vskip 0.2 cm
\noindent
{\bf ii)} the correct choice of $\Psi_{IN}$, which matches with
$\Psi$. This is:
$$
\Psi_{IN}^{Match}(\lambda,x)= \Psi_{IN}(\lambda,x)\pmatrix{ 1 &
  i\pi{r_1\over \rho}-\ln x \cr 
0
&
{1\over \rho}
}^{-1}= \Psi_{IN}(\lambda,x)~\pmatrix{ 1 & \rho~ \ln x - i \pi r_1 
\cr  
0 & \rho
}=
$$
$$
=~ 
K_0(x) \Psi_0(\mu)~  C_{IN},~~~~~~ 
C_{IN}=\pmatrix{0 & \rho/R \cr
                            1 &  \rho \ln x - i \pi r_1}
.
$$

\vskip 0.2 cm
\noindent
{\it Important Remark:} $C_{IN}$ depends on $x$. But the
monodromy must be independent of $x$. For this reason, we will have to
consider the substitution 
$$
\Psi_{IN}^{Match}\mapsto 
\Psi_{IN}^{Match} \bigl(C_{IN}\bigr)^{-1}= K_0(x)
\Psi_0(\mu).$$ 
 This makes the
monodromy at $\lambda=0,x$ independent of $x$. 
The corresponding
transformation:  $$
\Psi_{OUT}^{Match}\mapsto 
\Psi_{OUT}^{Match} {C_{IN}}^{-1},
$$ 
changes $M_1$, $M_\infty$, but does not introduce a dependence on
$x$. Namely:
{\small
$$
 M_1\mapsto C_{IN}M_1 {C_{IN}}^{-1}= 
\pmatrix{ 1 & 0
 \cr 
-2\pi i {r_1\over \rho}~R & 1 }
,~~~
 M_\infty\mapsto C_{IN}M_\infty {C_{IN}}^{-1}= 
\pmatrix{- 1 & 0
 \cr 
2\pi i \left({r_1\over \rho}-1\right)~R & -1 }
.
$$
}

\vskip 0.5 cm 
\noindent
{\bf MATCHING $\Psi\leftrightarrow \Psi_{IN}$.}
\vskip 0.2 cm
\noindent
$\Psi$ and $\Psi_{IN}^{Match}(\lambda,x)$ are matching by
construction. 
The monodromy of $\Psi$ at $\lambda=0,x$
coincides with that of $\Psi_{IN}^{Match}(\lambda,x)$.
 In order to compute it, we write: 
$$
\Psi_0(\mu)=\mu^{\theta_0\over 2}(\mu-1)^{\theta_x\over 2}\Phi_0(\mu),~~~~~
\Phi_0(\mu)=\pmatrix{\varphi_1 & \varphi_2 \cr \xi_1 &\xi_2}.
$$
 $\Phi_0$ satisfies (\ref{systemPhi0}) with $\sigma=1$. It 
is expressed in terms of two independent solutions $\varphi_1$,
 $\varphi_2$ of
the Gauss hyper-geometric equation (see Appendix 1):
$$
\mu(1-\mu)~ {d^2 \varphi \over d\mu^2} +\bigl(1+c-(a+[b+1]+1)~\mu
\bigr)~ {d\varphi\over d\mu} 
-a(b+1)~\varphi=0,
$$
where, 
$$
a={\theta_0\over 2}+{\theta_x\over 2}-{1\over 2},~~~ 
b+1= {\theta_0\over 2}+{\theta_x\over 2}+{3\over 2},~~~ 
c+1= \theta_0+1.
$$    
From $\varphi_1$ and $\varphi_2$ we compute: 
$$
\xi_i= {1\over r}\left[
\mu(1-\mu)~{d\varphi_i\over d\mu} ~-a\left(
\mu+{b-c \over a-b}
\right)~\varphi_i
\right],~~~i=1,2.
$$

\vskip 0.2 cm
We need a compete set of solutions at $\mu=0,1,\infty$. Since $a-b$ is
an integer, we are in a logarithmic case. We  briefly
explain some preliminary facts. 
Let us consider a Gauss Hyper-geometric equation in
standard form:
$$
\mu~(1-\mu)~ {d^2 \varphi \over d\mu^2} +\bigl[\gamma-(\alpha+\beta+1)~\mu \bigr]~ {d\varphi\over d\mu}
-\alpha\beta~\varphi=0
$$
($\alpha,\beta,\gamma$ here are not   the coefficients of (PVI)! We
are just using the same symbols only here). 
Logarithmic solutions at $\mu=0$ may occur only if $\gamma\in{\bf Z}$.
Logarithmic solutions at $\mu=1$ may occur only if 
$\alpha+\beta-\gamma\in{\bf Z}$. 
Logarithmic solutions at $\mu=\infty$ may occur only if
$\alpha-\beta\in{\bf Z}$. Several sub-case must be distinguished, and
this is not the place to discuss them. 

In our case $\alpha=a$, $\beta=b+1$ and $\gamma=c+1$. Therefore
 $\alpha-\beta = -2$.  This case is logarithmic. Two independent
 solutions are:
{\small
$$
\varphi_1^{(\infty)}=
\mu^{-\beta}g_1\left(\beta,\beta-\gamma+1,1+\beta-\alpha;~{1\over
  \mu}\right)
=[\hbox{ in our case }]
~~\mu^{-\alpha-2}g_1\left(\beta,\beta-\gamma+1,3;~{1\over
  \mu}\right),
$$
$$
\varphi_2^{(\infty)}=
\mu^{-\beta}F\left(\beta,\beta-\gamma+1,1+\beta-\alpha;~{1\over
  \mu}\right)
=[\hbox{ in our  case }]~~\mu^{-\beta}F\left(\beta,\beta-\gamma+1,3;~{1\over
  \mu}\right).
$$
}
Here $F$ is the Gauss hyper-geometric function and $g_1$ is a logarithmic
solution introduced in Norlund's paper \cite{Norlund}, page 7:
$$
g_1(u,v,w;z)= \sum_{n=1}^{w-1} (-1)^{n-1}(n-1)!~{(u)_{-n}
(v)_{-n}\over (w)_{-n}} {1\over z^n}~+ F(u,v,w; ~z)~\ln(-z)~+
$$
$$+\sum_{n=0}^\infty {(u)_n(v)_n\over n!(w)_n}
\bigl[\psi(1-u-n)+\psi(v+n)-\psi(w+n)-\psi(1+n)
\bigr]z^n,
$$
where
$$
\psi(x):={d \over dx}\ln \Gamma(x),~~~~~
% \hbox{ with poles } x=0,-1,-2,...
x\neq 0,-1,-2,-3,...~;~~~~|z|<1, ~~\ln(-z)<0 \hbox{ for } -1<z<0.
$$
%$$
%\psi(x):={d \over dx} \Gamma(x), \hbox{ with poles } x=0,-1,-2,...
%$$
Note that the first sum $\sum_{n=1}^{w-1} (-1)^{n-1}(n-1)!~{(u)_{-n}
(v)_{-n}\over (w)_{-n}} {1\over z^n}$ is a polynomial in
$1/z$.

%%%%%%%%%%%%%%%%%%%%%%%%%%%%%%%%
%%%%%%%%%%%%%%%
%%%%%%%%%%         CANCELED IN THE PREPRINT
%%%%%%%%%%%%%%%
%%%%%%%%%%%%%%%%%%%%%%%%%%%%%%%%
%\vskip 0.2 cm 
%If $\gamma\not \in {\bf Z}$, two independent solutions at $\mu=0$ are:
%$$
%\varphi_1^{(0)}= F(\alpha,\beta,\gamma;\mu),
%$$
%$$
%\varphi_2^{(0)}=\mu^{1-\gamma}
%F(\alpha-\gamma+1,\beta-\gamma+1,2-\gamma;~\mu).
%$$
%
%If $\alpha+\beta-\gamma\not\in{\bf Z}$, two independent solutions at 
%$\mu=1$ are:
%$$
%\varphi_1^{(1)}= F(\alpha,\beta,\alpha+\beta-\gamma+1;~1-\mu),
%$$
%$$
%\varphi_2^{(1)}= (1-\mu)^{\gamma-\alpha-\beta} F(\gamma-\beta,
%\gamma-\alpha,1+\gamma-\alpha-\beta;~1-\mu).
%$$
%
If $\theta_0,\theta_x\not\in{\bf Z}$, we fall in the non-logarithmic
cases at $\mu=0,1$. Thus:
$$
\varphi_1^{(0)}= 
F\left(
{\theta_0\over 2}+{\theta_x\over 2}-{1\over 2} , {\theta_0\over 2}+{\theta_x\over 2}+{3\over 2},1+\theta_0;\mu
\right),
$$
$$
\varphi_2^{(0)}= \mu^{-\theta_0}
F\left(-
{\theta_0\over 2}+{\theta_x\over 2}-{1\over 2} , -{\theta_0\over 2}+{\theta_x\over 2}+{3\over 2},1-\theta_0;\mu
\right).
$$
\vskip 0.2 cm
$$
\varphi_1^{(1)}= 
F\left(
{\theta_0\over 2}+{\theta_x\over 2}-{1\over 2} , {\theta_0\over 2}+{\theta_x\over 2}+{3\over 2},1+\theta_x;1-\mu\right),
$$
$$
\varphi_2^{(1)}= (1-\mu)^{-\theta_x}
F\left(
{\theta_0\over 2}-{\theta_x\over 2}-{1\over 2} , {\theta_0\over 2}-{\theta_x\over 2}+{3\over 2},1-\theta_x;1-\mu
\right).
$$
\vskip 0.2 cm
$$
\varphi_1^{(\infty)}= \mu^{-{\theta_0\over 2}-{\theta_x\over 2}-{3\over
    2}}~
g_1\left(
{\theta_0\over 2}+{\theta_x\over 2} +{3\over 2} , -{\theta_0\over
  2}+{\theta_x \over 2} +{3\over 2}, 3;{1\over \mu}
\right),
$$
$$
\varphi_2^{(\infty)}= \mu^{-{\theta_0\over 2}-{\theta_x\over 2}-{3\over
    2}}~
F\left(
{\theta_0\over 2}+{\theta_x\over 2} +{3\over 2} , -{\theta_0\over
  2}+{\theta_x \over 2} +{3\over 2}, 3;{1\over \mu}
\right).
$$

\vskip 0.2 cm
The connection matrix between $[ \varphi_1^{(0)},\varphi_2^{(0)}]$
and $[\varphi_1^{(1)},\varphi_2^{(1)}]$ is a standard one:
$$
[\varphi_1^{(0)},\varphi_2^{(0)}]
=
[\varphi_1^{(1)},\varphi_2^{(1)}]C_{01},~~~~~0<\hbox{~arg~}(1-\mu)<2\pi;
$$
where $
C_{01}$ is 
(\ref{C01PASQUA}). 

\vskip 0.2 cm
 On the
other hand,  the
connection matrix between $[ \varphi_1^{(0)},\varphi_2^{(0)}]$ and 
 $[\varphi_1^{(\infty)},\varphi_2^{(\infty)}]$ is computed in 
\cite{Norlund}. Our specific case falls in the case $
\alpha-\beta=-p$, $p\geq 0$ integer. 
From \cite{Norlund}, page 27, case 11, formulae  (1) and (2), we deduce the
connection formulae:
$$
(-\mu)^{\alpha-\beta} F\left(
\beta,\beta-\gamma+1,1+\beta-\alpha;{1\over \mu}
\right)
=
{
\Gamma(\beta-\alpha+1)\Gamma(1-\gamma)
\over 
\Gamma(1-\alpha)\Gamma(\beta-\gamma+1)
}
(-\mu)^{\alpha}F(\alpha,\beta,\gamma;\mu)~+
$$
$$~~~~~~~~~~~~~~~~
+ {
\Gamma(\beta-\alpha+1)\Gamma(\gamma-1)
\over 
\Gamma(\gamma-\alpha)\Gamma(\beta)
}
(-\mu)^{\alpha-\gamma+1}F(\alpha-\gamma+1,\beta-\gamma+1,2-\gamma;\mu),
$$
\vskip 0.2 cm
$$
(-\mu)^{\alpha-\beta}g_1\left(\beta,\beta-\gamma+1,1+\beta-\alpha;{1\over
  \mu}\right) 
=~~~~~~~~~~~~~~~~~~~~~~~~~~~~~~~~~~~~~~~~~~~~~~~~~~~~~~~~~~
$$
$$
~~~~~~~~~~~
(-1)^{\alpha-\beta+1}
{
\Gamma(1+\beta-\alpha)\Gamma(1-\beta)\Gamma(\alpha-\gamma+1)
\over 
\Gamma(2-\gamma)
}
~(-\mu)^{\alpha-\gamma+1}F(\alpha-\gamma+1,\beta-\gamma+1,2-\gamma;\mu)
.
$$
Here $|\hbox{arg}~(-\mu)|<\pi$, $\ln(-\mu)>0$ for
$-\infty<\mu<-1$. The branch cut in the $\mu$-plane is
$[0,+\infty)$. As for the minus signs, 
we choose $-\mu=e^{-i\pi}\mu$, then
  $0<\arg\mu<2\pi$. Moreover, in our case
$(-1)^{\alpha-\beta}=1$. With this preparation, 
we can write the connection matrix:
$$
[\varphi_1^{(\infty)},\varphi_2^{(\infty)}]
=
[\varphi_1^{(0)},\varphi_2^{(0)}]~C_{\infty 0},~~~~~0<\arg\mu<2\pi;
$$
where $
C_{\infty0}$ is (\ref{Cinfty0PASQUA}).

\vskip 0.3 cm 
In order to write $\Psi_{IN}^{Match}$ in terms of the
$\varphi_i^{(\infty)}$, let us compute the behavior for $\mu \to
\infty$ of $\mu^{\theta_0/2}(\mu-1)^{\theta_x/2}
\pmatrix{\varphi_1^{(\infty)} &\varphi_2^{(\infty)} 
\cr \xi_1^{(\infty)} &\xi_2^{(\infty)}} $. We have:
$$
\mu^{\theta_0\over 2}(\mu-1)^{\theta_x \over 2}~\varphi_1^{(\infty)}
=
-{32\over [(\theta_0+\theta_x)^2-1][(\theta_0-\theta_x)^2-1]}
~\mu^{1\over 2} \left[1+O\left({1\over \mu}\right)\right]
~+\mu^{-{3\over 2}}\ln \mu \left[1+O\left({1\over \mu}\right)\right].
$$
$$
\mu^{\theta_0\over 2}(\mu-1)^{\theta_x \over 2}~\varphi_2^{(\infty)}
=
\mu^{-{3\over 2}} \left[1+O\left({1\over
    \mu}\right)\right].
$$
\vskip 0.2 cm 
$$
\mu^{\theta_0\over 2}(\mu-1)^{\theta_x \over 2}~\xi_1^{(\infty)}
=
-{2\over r} ~\mu^{-{1\over 2}} \ln\mu  \left[1+O\left({1\over
  \mu}\right)\right]
+
O(\mu^{-{1\over 2}}).
$$
$$
\mu^{\theta_0\over 2}(\mu-1)^{\theta_x \over 2}~\xi_2^{(\infty)}
=
{2\over r} ~\mu^{-{1\over 2}}  \left[1+O\left({1\over
  \mu}\right)\right].
$$
Therefore:
$$
\mu^{\theta_0\over 2}(\mu-1)^{\theta_x \over 2}
\pmatrix{\varphi_1^{(\infty)} &\varphi_2^{(\infty)} 
\cr \xi_1^{(\infty)} &\xi_2^{(\infty)}}
=
\left[I+O\left({1\over \mu}\right)\right]~\mu^{{1\over 2}\sigma_3}
~
\pmatrix{1 & 0 \cr
R~ \ln\mu & 1 
}
~\pmatrix{ -{2\over r R} & 0 \cr 0 & {2\over r} 
}.
$$
We conclude that 
the matrix $\Psi_0$ used in the $\Psi_{OUT}\leftrightarrow \Psi_{IN}$
matching is:
{\small
$$
\Psi_0= \mu^{\theta_0\over 2}(\mu-1)^{\theta_x \over 2}
\pmatrix{\varphi_1^{(\infty)} &\varphi_2^{(\infty)} 
\cr \xi_1^{(\infty)} &\xi_2^{(\infty)}
}~C,~~~~~~~C:=
\pmatrix{-R~{r\over 2} & 0 \cr 0 & {r\over 2}}=
\pmatrix{{[1-(\theta_0+\theta_x)^2][(\theta_0-\theta_x)^2-1]
\over 32} & 0 \cr 0 & {r\over 2}
}
.
$$
}

\vskip 0.3 cm
As we have already remarked, we do the transformation $
\Psi_{IN}^{Match}\mapsto 
\Psi_{IN}^{Match}{C_{IN}}^{-1}$, to  compute the
$x$-independent  monodromy at $\lambda=0,x$. This means that we have
to compute the monodromy of:   
$$
\Psi_{IN}^{Match}{C_{IN}}^{-1}= K_0(x)
\Psi_0\left(\lambda/x\right)=
$$
$$
= K_0(x)~ \left({\lambda\over x}\right)^{\theta_0\over 2}\left({\lambda\over x}-1\right)^{\theta_x \over 2}
\pmatrix{\varphi_1^{(\infty)} &\varphi_2^{(\infty)} 
\cr \xi_1^{(\infty)} &\xi_2^{(\infty)}
}~C,
$$
$$
= 
  K_0(x)~\left({\lambda\over x}\right)^{\theta_0\over 2}\left({\lambda\over x}-1\right)^{\theta_x \over 2}
\pmatrix{\varphi_1^{(0)} &\varphi_2^{(0)} 
\cr \xi_1^{(0)} &\xi_2^{(0)}
}~C_{\infty0}~C,
$$
$$
=  K_0(x)~ \left({\lambda\over x}\right)^{\theta_0\over
  2}\left({\lambda\over x}-1\right)^{\theta_x \over 2} 
\pmatrix{\varphi_1^{(1)} &\varphi_2^{(1)} 
\cr \xi_1^{(1)} &\xi_2^{(1)}
}~C_{01}~C_{\infty0}~C.
$$
In order to do this, we observe that from the definition of
$\varphi_1^{(0)},\varphi_2^{(0)}$ if follows that:
$$
 K_0(x)~\left({\lambda\over x}\right)^{\theta_0\over
 2}\left({\lambda\over x}-1\right)^{\theta_x \over 2} 
\pmatrix{\varphi_1^{(0)} &\varphi_2^{(0)} 
\cr \xi_1^{(0)} &\xi_2^{(0)}
}= \psi_0^{IN}(x) (I+O(\lambda))\lambda^{{\theta_0\over 2}\sigma_3}
,
$$
for $\lambda\to 0$. 
From the definition of
$\varphi_1^{(1)},\varphi_2^{(1)}$, it follows that:
$$
K_0(x)~\left({\lambda\over x}\right)^{\theta_0\over
  2}\left({\lambda\over x}-1\right)^{\theta_x \over 2} 
\pmatrix{\varphi_1^{(1)} &\varphi_2^{(1)} 
\cr \xi_1^{(1)} &\xi_2^{(1)}
}= \psi_x^{IN}(x)(I+O(\lambda-x))(\lambda-x)^{{\theta_x\over 2}\sigma_3}
,
$$
for $\lambda\to x$.  
It is not necessary here to write explicitly the invertible 
 matrices $\psi_0^{IN}(x)$ and $\psi_x^{IN}{x}$. 
The above construction implies 
that  $\Psi_{IN}^{Match}{C_{IN}}^{-1}$
has monodromy matrices at $\lambda=0,x$ respectively given by:
$$
M_0= C^{-1} \bigl(C_{\infty0}\bigr)^{-1}~
\exp\{i\pi\theta_0\sigma_3\} ~C_{\infty0} C,
$$
$$
M_x= C^{-1} \bigl(C_{\infty0}\bigr)^{-1}\bigl(C_{01}\bigr)^{-1}~
\exp\{i\pi\theta_x\sigma_3\} ~C_{01}C_{\infty0} C.
$$
These matrices coincide with the monodromy matrices of
$\Psi~{C_{IN}}^{-1}$. 

\vskip 0.3 cm
As a last simplification, we 
 consider the transformation $M\mapsto C M C^{-1}$. We
obtain the result of theorem \ref{thMONODROMY}, case c). Namely:
 $$
M_0=  \bigl(C_{\infty0}\bigr)^{-1}~
\exp\{i\pi\theta_0\sigma_3\} ~C_{\infty0},~~~~~
M_{\infty}= \pmatrix{ -1 & 0 \cr 2\pi i \left(1-{r_1\over \rho}\right) & -1 },
$$
$$
M_1= \pmatrix{ 1 & 0 \cr 2\pi i {r_1\over \rho} & 1},~~~~~M_x=\bigl(C_{\infty0}\bigr)^{-1}\bigl(C_{01}\bigr)^{-1}~
\exp\{i\pi\theta_x\sigma_3 \}~C_{01}C_{\infty0}.
$$
We observe that $r$ does not appear in the monodromy matrices. 

%%%%%%%%%%%%%%%%%%%%%%%%%%%%%%%%%%%%%%%%%%%%%%%%%%%%

%%%%%%%%%%%%%%%%%%%%%%%%%%%%%%%%%%%%%%%%%%%%%%%

%%%%%%%%%%%%%%%%%%%%%%%%%%%%%%%%%%%%%%%%%%%%%

%%% APPENDIX 

%%%%%%%%%%%%%
%%%%%%%%%        \small
%%%%%%%%%%%%%

{
%\small 

\section{ Appendix 1}

% PROPOSITION 

\bpr 
\label{matrices}
Let $B_0$, $B_1$ be $2\times 2$ matrices such that 
$$
\hbox{Eigenvalues }(B_0)=0,-c,~~~\hbox{Eigenvalues }(B_1)=0,c-a-b.
$$
and $B_0+B_1$ is either  diagonalizable: 
$$ 
B_0+B_1=\pmatrix{ -a & 0 \cr 0 & -b} ~~\hbox{ (it may happen that }a=b),
$$
or it is a  Jordan form: 
$$ 
B_0+B_1=\pmatrix{ -a & 1 \cr 0 & -a}.
$$
Then, $B_0$ and $B_1$  can be computed as in the following cases. 
Let $r$, $s$  be any complex numbers. 
\vskip 0.2 cm
\noindent 
{\bf 1) Diagonalizable case.} 

\noindent
Case $a\neq b$:
\be
B_0:= \pmatrix{ {a(b-c)\over a-b} & r 
 \cr 
    {ab(a-c)(c-b)\over r(a-b)^2} & {b (c-a)\over a-b}
},~~
B_1= \pmatrix{ {a(c-a)\over a-b} & -r 
\cr
               -(B_0)_{21} & {b(b-c)\over a-b}
},~~~r\neq 0
\label{1}
\ee
\be
\hbox{ if } a=0:~~~B_0=\pmatrix{0 & r 
                                          \cr
                                0 & -c },
~~
               B_1=\pmatrix{0 & -r \cr
                            0 & c-b}.
\label{2}
\ee
\be
\hbox{ If } b=0:~~~B_0=\pmatrix{ -c & r \cr
                                  0 & 0 },~~
                   B_1=\pmatrix{ c-a & -r \cr
                                   0 & 0}.
\label{3}
\ee
\be
\hbox{ If } a=c\neq b:~~~B_0=\pmatrix{-a & r \cr
                                       0 & 0 },~~~
                         B_1=\pmatrix{0 & -r \cr
                                      0 & -b}.
\label{4}
\ee
\be
\hbox{ If } b=c\neq a: ~~~B_0=\pmatrix{0 & r \cr
                                       0 & -b},~~
                          B_1=\pmatrix{-a & -r \cr 
                                        0 & 0 }.
\label{5}
\ee
Cases (\ref{2})--(\ref{5}) are already included in (\ref{1}).

\vskip 0.2cm
\noindent
Case $a=b$. We have two sub-cases:
\be
 \hbox{ If } a=b=c:~~~B_0=\pmatrix{ -c-s & r \cr
                                    -{s(c+s)\over r} & s},~~
                      B_1=\pmatrix{s & -r \cr
                                   {s(c+s)\over r} & -c-s}.
\label{6}
\ee
\be
\hbox{ If } a=b=0:~~~ B_0=\pmatrix{-c-s & r \cr
                                   -{s(c+s)\over r} & s},
   ~~~~~                   B_1=-B_0.
\label{7}
\ee
 The transpose matrices of all the above cases are also possible.

\vskip 0.2 cm
\noindent
{\bf 2) Jordan case.}  

\noindent
For $a\neq 0$ and $a\neq c$ we have:
\be
B_0=\pmatrix{ r & {r(r+c)\over a(a-c)} \cr 
             a(c-a) & -c-r },~~
B_1=\pmatrix{-a-r & 1-{r(r+c)\over a(a-c)} \cr
             a(a-c) & c-a+r }.
\label{8}
\ee
For $a=0$, or $a=c$, we have two possibilities:
\be
 B_0=\pmatrix{0 & r \cr 0 & -c},~~ 
                    B_1=\pmatrix{-a & 1-r \cr 0 & -a+c }; 
\label{9}
\ee
or
\be
B_0=\pmatrix{-c & r \cr 0   & 0} ,~~~B_1=\pmatrix{ c-a & 1-r \cr 0 & -a}
\label{10}
\ee
\epr

%%%%%%%%%%%%%
%%%%%
%%%       CANCELED IN THE PREPRINT
%%%%%
%%%%%%%%%%%%%
%\noindent
%{\it Proof:} Write $B_0=\pmatrix{ b_{11} & b_{12} \cr b_{21} & b_{22}}$, 
%$B_1=\pmatrix{ c_{11} & c_{12} \cr c_{21} & c_{22}}$, and solve the system 
%in eight equations and eight variables 
%$$
%\left.
%\matrix{
%\hbox{det}~ B_0=0, &~~~ \hbox{tr}~ B_0=-c, \cr
%\hbox{det}~ B_1=0, &~~~ \hbox{tr}~B_1=c-a-b, \cr
%b_{11}+c_{11}=-a, &~~~b_{22}+c_{22}=-b, \cr
%b_{12}+c_{12}=0 \hbox{ or } 1,& ~~~b_{21}+c_{21}=0. 
%}
%\right.
%$$
%\qed
%%%
%%%%%%

%%%%%%%%%%%%%%%

\bpr 
Let $B_0$ and $B_1$ be as in Proposition \ref{matrices}. 
The linear system:  
$$
{d\over dz} \pmatrix{ \varphi \cr \xi} = \left[
{B_0\over z} +{B_1\over z-1}
\right]~\pmatrix{ \varphi \cr \xi}
$$
may be reduced  to a Gauss hyper-geometric equation, in the
following cases. 

\vskip 0.2 cm
\noindent
Diagonalizable case (i.e. from (\ref{1}) to (\ref{7})):
\be
z(1-z)~ {d^2 \varphi \over dz^2} +\bigl(1+c-(a+[b+1]+1)~z \bigr)~ {d\varphi\over dz}
-a(b+1)~\varphi=0.
\label{hypergeom1}
\ee
The component $\xi$ is obtained by
 the following equalities,
 according to the different cases of Proposition \ref{matrices}.

\vskip 0.2 cm
\noindent 
Cases (\ref{1}) (\ref{2}) (\ref{3}) (\ref{4}) (\ref{5}):
\be
\xi= {1\over r}\left[
z(1-z)~{d\varphi\over dz} ~-a\left(
z+{b-c \over a-b}
\right)~\varphi
\right]
\label{xi-hypergeom1}
\ee

\vskip 0.2 cm
\noindent
Case (\ref{6}): 
$$
\xi= {1\over r} \left[
z(1-z)~{d\varphi\over dz} +(c+s-c~z)~\varphi
\right]
$$

\vskip 0.2 cm
\noindent
Case (\ref{7}):
$$
\xi={1\over r} \left[
z(1-z)~{d\varphi\over dz} + (c+s)~\varphi
\right]
$$

\vskip 0.2 cm
\noindent
Jordan case (\ref{8}): The equation for  $\varphi$ is in Gauss hypergeometic form  only when 
$r=-a$. In this case, the matrices (\ref{8}) are:
$$
B_0=\pmatrix{ -a & 1 \cr a(c-a) & a-c },~~~B_1=\pmatrix{0 & 0 \cr 
a(a-c) & c-2a}.
$$
The equation is: 
$$
z(1-z) ~{d^2\varphi\over dz^2} +\left[
1+c-(2a+1)~z
\right]~{d \varphi \over dz} 
-a^2 ~\varphi =0,
$$
$$
\xi=  z~{d\varphi\over dz}+a~\varphi.  
$$

\vskip 0.25 cm
\noindent
Jordan case (\ref{9}): for $r=1$ we get: 
$$\hbox{For  }a=0,~~~~~~~~~~~
z(1-z){d^2 \varphi \over dz^2} +(1+c-z){d\varphi\over dz}=0.
$$
$$\hbox{For  }c=a,~~~
z(1-z){d^2 \varphi \over dz^2} +(1+a-(2a+1)z){d\varphi\over dz}
+\left(-a^2-{a\over 1-z}\right)~\varphi=0.
$$
The Gauss form appears only when $a=0$.

\vskip 0.25 cm
\noindent
Jordan case (\ref{10}): for $r=1$ we get:
$$ \hbox{For }a=0,~~~
z(1-z){d^2 \varphi \over dz^2} +(1+c-z){d\varphi\over dz}
+{c\over 1-z}~\varphi=0.
$$
$$\hbox{For }
c=a,~~~
z(1-z){d^2 \varphi \over dz^2} +(1+a-(2a+1)z){d\varphi\over dz}
-a^2\varphi=0.
$$
For $r\neq 1$, we don't get a Gauss hyper-geometric form  for the equation of 
$\varphi$ in both cases (\ref{9}) and (\ref{10}). Nevertheless, the matrices are in upper triangular form, so the equation for $\xi$ is solvable by elementary integration.  
\vskip 0.2 cm
\noindent
Jordan case (\ref{8}): The equation for  $\xi$ is in Gauss hypergeometic form: 
$$
z(z-1){d^2\xi\over dz^2}+\bigr(1+c-2(a+1)z\bigl){d\xi\over dz} 
-a(a+1)\xi=0,
$$
$$
\varphi(z)={1\over a(a-c)}
\left[
z(z-1){d\xi\over dz}+(az-c-r)\xi
\right].
$$

\label{ipergeom}
\epr

\vskip 0.3 cm
 As a compendium to the above proposition, we recall that any 
irreducible representation of $\pi_1({\bf
 CP}^1\backslash \{0,1,\infty\}) 
\mapsto GL(2,{\bf C})$ can be realized as the monodromy group of a Riemann (or Gauss) equation. A reducible representation (namely, $M_0,M_1,M_\infty$ are in 
upper triangular form) can be realized by the monodromy of a $2\times 2$ 
Fuchsian system 
$$ {dY\over dz}= \left[
{B_0\over z}+{B_1\over z-1}\right]~Y
$$
where $B_0$, $B_1$ are $2\times 2$ upper triangular matrices. 
We also  state the following:

\ble
Consider a $2\times 2 $ linear system:
$$
{d Y(z)\over dz}=A(z) ~Y,
~~~~~
A(z)=\pmatrix{a(z) & b(z) \cr 0 & c(z)}
$$
such that $A(z)$ is meromorphic, with poles $a_1$, $a_2$, ..., $a_N$, $\infty$. The monodromy group is generated by $N$ upper triangular monodromy matrices:
$$
M_i= \pmatrix{\lambda_1^{(i)} & R^{(i)} 
\cr 
0 & \lambda_2^{(i)} 
},~~~i=1,2,...,N;
$$
where $\lambda_1^{(i)},\lambda_2^{(i)},R^{(i)}$ are constants (i.e. they are independent of $z$) given by:
$$
\lambda_1^{(i)}=\exp\{2\pi i~\hbox{\rm Res}~a(z)|_{a_i}\},~~~\lambda_2^{(i)}=\exp\{2\pi i~\hbox{\rm Res}~c(z)|_{a_i}\},
$$
$$
R^{(i)}=\int_z^{z\exp\{2\pi i\}}ds~ b(s)~ {u_2(s)\over u_1(s)}~+~\bigl(\lambda_1^{(i)}-\lambda_2^{(i)}\bigr)\int_{z_0}^zds~ b(s)~{u_2(s)\over u_1(s)}.
$$
$R^{(i)}$ depends on a chosen non-singular base point $z_0$, but not
on $z$. One of the non zero $R^{(i)}$ can be put equal to 1, by
re-defining $Y\mapsto YC$, $C=\pmatrix{1& 0\cr 0&1/R^{(i)}}$. 

\ele 

%%%%%%%%%%%%%

\vskip 0.3 cm
\noindent
{\it Proof:} 
Let us write $Y=\pmatrix{y_1 \cr y_2}$ and the equation in the form:
$$
{d y_1\over dz}=a(z)y_1~+b(z)y_2,~~~~~
{dy_2\over dz}= c(z)y_2.
$$
Let $z_0\neq a_i,\infty$, $i=1,...,N$. The second equation has solution $
y_2(z)= C_2~u_2$, where  $ u_2:=\exp\{\int_{z_0}^z ds
~c(s)\}$, $C_2\in{\bf C}$. 
The first equation becomes:
\be
{d y_1\over dz}=a(z)y_1~+~C_2~b(z)~u_2(z).
\label{nonho}
\ee
We solve the first equation by variation of parameters. Let $
u_1(z)=\exp\{\int_{z_0}^z ds ~a(s)\}
$ 
be a fundamental solution of the homogeneous equation ${d y_1\over dz}=
a(z)y_1$. 
We look for a solution of  (\ref{nonho}) of the 
form $y_2(z)= w(z)~u_1(z)$. 
Substitution gives:
$$
    {d w \over dz}= C_2~b(z)~{u_2(z)\over u_1(z)}~~\Longrightarrow~~
w(z)=C_2~v(z)~+~C_1,~~~~~C_1\in{\bf C}.
$$
where:
$$
v(z)=~\int_{z_0}^zds~ b(s)~{u_2(s)\over u_1(s)}
$$
The general solution of  (\ref{nonho}) is $
y_1(z)=C_1~u_1(z)~+C_2~v(z)~u_1(z)$. 
Then, a fundamental solution for the initial system can be chosen to be:
$$
Y(z)= \pmatrix{ u_1(z) & v(z)~u_1(z) \cr
\cr
                 0 & u_2(z)}
$$
We  compute the monodromy for $(z-a_i)\mapsto (z-a_i) e^{2\pi i}$. We have:
$$
u_1(z)\mapsto \lambda_1^{(i)} ~u_1(z),~~~\lambda_1^{(i)}=\exp\{2\pi i
~\hbox{\rm Res}~ a(z)|_{a_i}\}, 
 $$
\be
u_2(z)\mapsto \lambda_2^{(i)} ~u_1(z),~~~\lambda_2^{(i)}=\exp\{2\pi i
~\hbox{\rm Res}~ c(z)|_{a_i}\}. 
\label{myvirtue}
\ee
By linearity, the vector solution $\pmatrix{ v(z)~u_1(z)\cr u_2(z)}$
is transformed into a linear combination of two independent vector
solutions: $R^{(i)}\pmatrix{u_1(z) \cr 0}~+ 
S^{(i)}  \pmatrix{ v(z)~u_1(z)\cr u_2(z)}$. Moreover, $S^{(i)}$ must
coincide with  
$\lambda_2^{(i)}$, because -- by virtue of (\ref{myvirtue}) --
$u_2(z)\mapsto \lambda_2^{(i)}u_2(z)$. Namely:  
$$
\pmatrix{ v(z)~u_1(z)\cr u_2(z)}\mapsto 
~R^{(i)}~\pmatrix{u_1(z) \cr 0}~+~
\lambda_2^{(i)}  ~\pmatrix{ v(z)~u_1(z)\cr u_2(z)},~~~~~R^{(i)} \in{\bf C}. 
$$
Thus:
$$
Y(z)=\pmatrix{ u_1(z) & v(z)~u_1(z) \cr
\cr
                 0 & u_2(z)}~\mapsto ~\pmatrix{\lambda_1^{(i)} ~ u_1(z) & 
 \lambda_2^{(i)} ~v(z)~u_1(z)~+~R^{(i)}~u_1(z) 
\cr
\cr
0 &
     \lambda_2^{(i)}  ~  u_2(z)}~=
$$
$$
=~Y(z)~\pmatrix{\lambda_1^{(i)} & R^{(i)} \cr 0 & \lambda_2^{(i)}}.
$$

Let $C_i$ be a small loop around $a_i$. 
To find $R^{(i)}$, let us observe that   $
u_1(z)~v(z)\mapsto \lambda_1^{(i)} u_1(z)~\tilde{v}(z)$, 
where: 
$
 \tilde{v}(z)=v(z)~+K_i(z)$, $K_i(z):=\int_z^{z\exp\{2\pi i\}}ds~
 b(s)~ {u_2(s)\over u_1(s)}$. 
Thus:
$$
Y(z)=\pmatrix{ u_1(z) & v(z)~u_1(z) \cr
\cr
                 0 & u_2(z)}~\mapsto ~\pmatrix{\lambda_1^{(i)} ~ u_1(z) & 
 \lambda_1^{(i)} ~v(z)~u_1(z)~+~K_i(z)u_1(z) 
\cr
\cr
0 &
     \lambda_2^{(i)}  ~  u_2(z)}.
$$
We must have: $\lambda_1^{(i)}v(z)u_1(z)~+K_i(z)u_1(z)\equiv
 \lambda_2^{(i)}v(z)u_1(z)~+R^{(i)}u_1(z)$; namely:
$$
R^{(i)}=K_i(z)~+\bigl(\lambda_1^{(i)}-\lambda_2^{(i)}\bigr)v(z).
$$
\qed

\section{Appendix 2: Formal Asymptotic Expansion}

{\bf 1)} We consider systems (\ref{IRR1}): 
$$
{dY\over dz}= \left[
\Omega+{D_1\over z} + \sum_{n=2}^\infty {D_n \over z^n}
\right]Y ~:=D(z)~Y,~~~~~
\Omega =\hbox{diag}(\omega_1,\omega_2,...,\omega_n),~~~~~
$$
with $\omega_i\neq \omega_j$, for $i\neq j$. 
We introduce a  gauge transformation $
Y=G(z)\tilde{Y}$, 
such that:
$$
{\tilde{Y}\over d z}=  \left[
G^{-1}(z)D(z)G(z)-G^{-1}(z) {dG(z)\over dz} 
\right]~\tilde{Y},
$$
be in the simple form:
$$
 {d\tilde{Y}\over d z}=  \left[\Omega+{\Omega_1\over z}
\right] \tilde{Y}, ~~~\hbox{ $\Omega$, $\Omega_1$ diagonal}.
$$
Formally, we write  $G(z)$ as:
$$
G(z)=I+{G_1\over z}+{G_2\over z^2}+...,~~~~~~~~z\to \infty.
$$
If the $G_n$'s can be determined, we get the formal solution: 
$$
Y(z)\sim \left[I+\sum_{n=1}^\infty {G_n\over z^n}\right]~\exp\left\{
z~\Omega +\Omega_1\ln(z)
\right\},~~~z\to\infty
$$
For a sector of angular with $\pi+\epsilon$, $\epsilon>0$
sufficiently small (but finite and non zero), there exists a  unique
 solution $Y(z)$ with the above asymptotic expansion \cite{BJL1}.

In order to determine $G_n$ and $\Omega_1$,
 we solve  $D(z)G(z)-\partial_zG(z)=G(z)(\Omega_0+\Omega_1 z^{-1})$:
$$
\left(
\Omega+\sum_{n=1}^\infty {D_n\over z^n}
\right)
\left(
I+\sum_{n=1}^\infty {G_n\over z^n}
\right)
+
\sum_{n=2}^\infty {(n-1)G_{n-1}\over z^n} 
=
\left( I+ \sum_{n=1}^\infty {G_n\over z^n}
\right)
\left(
\Omega+{\Omega_1\over z}
\right).
$$
We identify equal powers of $z^{-1}$. From the power $1/z$ we get: 
$$
\Omega_1=\hbox{ diag}(\omega_1^{(1)},...,\omega_n^{(1)}),~~~
\omega_i^{(1)}:=(\Omega_1)_{ii}= (D_1)_{ii},~~~~~~~
(G_1)_{ij}=-{(D_1)_{ij}\over \omega_i-\omega_j}.
$$
From the power $1/z^2$, for $i\neq j$,    
we compute $(G_2)_{ij}$,  and for $i=j$ we compute: 
$$
(G_1)_{ii}= -(D_2)_{ii}-\sum_{k\neq i} (D_1)_{ik}(G_1)_{ki}.
$$
From the power $1/z^n$ we get:
$$
(G_{n-1})_{ii}={1\over n-1} \left\{
-\bigl(
D_n+D_{n-1}G_1+...+D_2G_{n-2}
\bigr)_{ii}
-
\sum_{k\neq i} (D_1)_{ik}(G_{n-1})_{ki}
\right\}, 
$$
$$
(G_n)_{ij}= {1\over \omega_i-\omega_j }
\left\{
[\omega_j^{(1)}-\omega_i^{(1)} -(n-1)]
(G_{n-1})_{ij}- 
\sum_{k\neq i}(D_1)_{ik}(G_{n-1})_{kj}-
\right.
$$
$$
~~~~~~~~~~~~~~~-(D_n+D_{n-1}G_1+...+D_2G_{n-2})_{ij}
\Bigr\},~~~~~~i\neq j
$$

\vskip 0.5 cm
\noindent
{\bf 2)} We consider the system (\ref{IRR2}): 
$$
{dY\over dz}= \left[
x^2\Lambda~z+x\Lambda+{E_1\over z} + \sum_{n=2}^\infty {E_n \over z^n}
\right]Y ~:=E(z)~Y,~~~~~
\Lambda=\hbox{ diag}~(\lambda_1,\lambda_2,...,\lambda_n).
$$
with $\lambda_i\neq \lambda_j$ for $i\neq j$. 
We introduce a gauge transformation: 
$$
Y(z)=K(z)~\tilde{Y}(z)
$$
in order  to reach the simple form:
$$
{d\tilde{Y}\over d z}= \left[
K^{-1}(z)E(z)K(z)-K^{-1}(z) {dK(z)\over dz} 
\right]~\tilde{Y}
$$
$$
\equiv
\left[
x^2\Lambda~z+x\Lambda+{\Lambda_1\over z}
\right]~\tilde{Y},~~~~~\Lambda_1 \hbox{ diagonal.}
$$
Formally, we write:
$$
K(z)\sim I+{K_1\over z} + {K_2\over z}+...= I+\sum_{n=1}^\infty
{K_n\over z}.
$$
Provided that we can determine the matrices $K_n$, we obtain the
formal solution:
$$
Y(z)\sim \left[I+\sum_{n=1}^\infty{K_n\over z^n}\right] ~\exp\left\{
{x^2\over 2} \Lambda~z^2 ~+x\Lambda ~z~\Lambda_1\ln~x
\right\}
$$
For a sector of angular with ${\pi\over 2}+\epsilon$, $\epsilon>0$
sufficiently small (but finite and non zero), there exists a  unique
 solution $Y(z)$ with the above asymptotic expansion \cite{BJL1}.

In order to determine $K_n$ and $\Lambda_1$,
 we solve  $E(z)K(z)-\partial_z K(z)=K(z)(x^2\Lambda z^2+\Lambda z
+\Omega_1 z^{-1})$:
$$
\left(x^2\Lambda z^2+\Lambda z+\sum_{n=1}^\infty {E_n\over z^n}
\right)
\left(
I+\sum_{n=1}^\infty {K_n\over z^n}
\right)
+
\sum_{n=2}^\infty {(n-1)K_{n-1}\over z^n} =~~~~~~~~~~~
$$
$$
~~~~~~~~~~~~~~~~~~~~~~~~~~~~~~~~~~~~~~
=
\left( I+ \sum_{n=1}^\infty {K_n\over z^n}
\right)
\left(x^2\Lambda z^2+\Lambda z+{\Omega_1\over z}
\right).
$$
We identify  equal powers of $z^{-1}$. 

\vskip 0.2 cm 
\noindent
Power $z$. It is an identity: $x^2\Lambda = x^2 \Lambda$.

\vskip 0.2 cm 
\noindent
Power $z^0$. 
$$
  x\Lambda+x^2\Lambda K_1= x\Lambda + x^2 K_1 \Lambda~~~ \Longrightarrow ~~~
[\Lambda,K_1]=0
$$
This means that $K_1$ is a diagonal matrix. 

\vskip 0.2 cm 
\noindent
Power $1/z$. We obtain the equation: 
$$ 
x^2~[\Lambda,K_2]= \Lambda_1-E_1 ~~~ \Longrightarrow ~~~
(\Lambda_1)_{ii}= (E_1)_{ii},~~~(G_2)_{ij}= -{(E_1)_{ij} \over
  x^2(\lambda_i-\lambda_j)},~~i\neq j.
$$

  \vskip 0.2 cm 
\noindent
Power $1/z^2$. We obtain the equation:
$$
x^2[\Lambda,K_3]= x[K_2,\Lambda]+K_1\Lambda_1-E_1K_1-K_1-E_2
\equiv {E_1-\Lambda_1\over x} +K_1\Lambda_1-E_1K_1-K_1-E_2
$$
Thus:
$$
(K_1)_{ii}= -(E_2)_{ii}
$$
$$
(K_3)_{ij} = {1\over x^2(\lambda_i-\lambda_j)} 
\left[
(E_1)_{ij} \left(
{1\over x} + (E_2)_{jj}
\right)
-(E_2)_{ij}
\right],~~~i\neq j.
$$

\vskip 0.2 cm 
\noindent
Power $1/z^2$. We obtain the equation:
$$
x^2[\Lambda,K_4]= x[K_3,\Lambda] +K_2\Lambda_1-E_1K_2-E_2 K_1 -2 K_2
-E_3
$$
The diagonal part gives:
$$
2 (K_2)_{ij}= \bigl[(E_2)_{ii}\bigr]^2-(E_3)_{ii}-\sum_{k\neq i}
(E_1)_{ik}(K_2)_{ki}.
$$
The non diagonal part gives $(K_4)_{ij}$, $i\neq j$.

\vskip 0.2 cm 
We content ourselves with these results, namely the determination of
$K_1$ and $K_2$. With the same procedure, we can determine all the
$K_n$'s.

\section{Appendix 3: Birational Transformations}

All the solutions of (PVI)  
of the form:  
\be
y(x)=b_0+b_1x+b_2x^2+...+b_Nx^N +...
\label{seriea0}
\ee
are  obtained from the matching procedure of sections
\ref{irregular0}, \ref{semireducible}, \ref{supersemireducible}. By
this we mean that solutions of type (\ref{seriea0})  are the solutions
given by  the matching procedure, or they can be obtained from
solutions given by  the matching procedure via one of the  birational
transformations of  \cite{Okamoto} and the transformation (\ref{nuovasimmetria}). 

\vskip 0.2 cm
Birational transformations are
symmetries of (PVI), namely invertible transformations: 
$$
y^{\prime}(x)={P(x,y(x))\over Q(x,y(x))}, ~~~x^{\prime}={p(x)\over q(x)},~~~~~
(\theta_0,\theta_x,\theta_1,\theta_\infty)\mapsto 
(\theta_0^{\prime},\theta_x^{\prime},\theta_1^{\prime},\theta_\infty^{\prime})
$$ 
such that $y(x)$ satisfies (PVI) with coefficients $\theta_0,\theta_x,\theta_1,\theta_\infty$ and variable $x$,  if and only if $y^{\prime}(x^{\prime})$ satisfies (PVI) with coefficients $\theta_0^{\prime},\theta_x^{\prime},\theta_1^{\prime},\theta_\infty^{\prime}$ and variable $x^\prime$. The functions 
 $P,Q$ are polynomials; $p,q$ are linear; the transformation of the $\theta_\mu$'s
  is an element of a linear representation of one of the following groups. Permutation group;  the Weyl group of the root system $D_4$; the group of translations $v:=(v_1,v_2,v_3,v_4)\mapsto 
v+e_j$, j=1,2,3,4 (where  $e_1=(1,0,0,0)$, ..., $e_4=(0,0,0,1)$).
\footnote{
 We  note that $\theta_1$, $\theta_x$, $\theta_0$ are defined up to
the sign. Moreover,  $\theta_\infty$ is defined up to
$\theta_\infty\mapsto 2-\theta_\infty$, and for this reason symmetries
are described in \cite{Okamoto}  
in terms of $\chi_\infty:=\theta_\infty-1$ . 
}

\vskip 0.2 cm
* Permutation group:
$$
x^1:~~~\theta_1^{\prime}=\theta_0,~~\theta_0^{\prime}=\theta_1;~~~~~~~~
\theta_x^{\prime}=\theta_x,~~\theta_\infty^{\prime}=\theta_\infty;~~~~~~~~~~y^{\prime}(x)=1-y(x),~~~x=1-x^\prime. 
$$
$$
x^2:~~~\theta_0^{\prime}=\theta_\infty-1,~~\theta_\infty^{\prime}=\theta_0+1;~~~~~~~~
\theta_1^{\prime}=\theta_1,~~\theta_x^{\prime}=\theta_x;~~~~~~~~~~
y^{\prime}(x)={1\over y(x)},~~~x={1\over x^{\prime}}.
$$
$$
x^3:~~~\theta_0^{\prime}=\theta_x,~~\theta_x^{\prime}=\theta_0;~~~~~~~~
\theta_1^{\prime}=\theta_1,~~\theta_\infty^{\prime}=\theta_\infty;~~~~~~~~~~y^{\prime}(x)= {x-y(x)\over x-1},~~~x={x^{\prime}\over x^{\prime}-1}.
$$

\vskip 0.3 cm
* Weyl Group: 
\vskip 0.2 cm
\noindent
$w_1$:
$$\theta_1^\prime=-\theta_1;~~~~~\theta_0^\prime=\theta_0,~~ \theta_x^\prime= \theta_x,~~ \theta_\infty^\prime =\theta_\infty.
$$
$w_2$:
$$
\theta_0^\prime={\theta_0+\theta_1+\theta_x+\theta_\infty\over 2} -1,~~~
\theta_1^\prime={\theta_0+\theta_1-\theta_x-\theta_{\infty}\over 2}+ 1,
$$
$$
\theta_x^\prime={\theta_0-\theta_1+\theta_x-\theta_{\infty}\over 2} +1,~~~
\theta_{\infty}^\prime={\theta_0-\theta_1-\theta_x+\theta_\infty\over 2}+1
$$
$
w_3$:
$$
\theta_\infty^\prime=2-\theta_\infty;~~~~~\theta_0^\prime=\theta_0,~~ \theta_x^\prime= \theta_x,~~ \theta_1^\prime =\theta_1.
$$
$
w_4$:
$$
\theta_\infty^\prime=2-\theta_\infty;~~\theta_x^\prime=2-\theta_x;~~~~~ \theta_0^\prime= \theta_0,~~ \theta_1^\prime =\theta_1.
$$
\vskip 0.2 cm
\noindent
The variable $x^\prime=x$, but $y^{\prime}(x)$ is quite complicated
and will not be given here (see \cite{Okamoto}).

\vskip 0.3 cm

* Shift $l_j:~v\mapsto v+e_j$:
$$
l_1:~~~~~~~~~~\theta_0^\prime=\theta_0+1,~~\theta_1^\prime=\theta_1+1;~~~~~
 \theta_x^\prime=\theta_x,~~\theta_\infty^\prime=\theta_\infty.
$$
$$
l_2:~~~~~~~~~~\theta_0^\prime=\theta_0+1,~~\theta_1^\prime=\theta_1-1;~~~~~
 \theta_x^\prime=\theta_x,~~\theta_\infty^\prime=\theta_\infty.
$$
$$
l_3:~~~~~~~~~~\theta_x^\prime=\theta_x+1,~~\theta_\infty^\prime=\theta_\infty+1;~~~~~
 \theta_0^\prime=\theta_0,~~\theta_1^\prime=\theta_1.
$$
$$
l_4:~~~~~~~~~~\theta_x^\prime=\theta_x+1,~~\theta_\infty^\prime=\theta_\infty-1;~~~~~
 \theta_0^\prime=\theta_0,~~\theta_1^\prime=\theta_1.
$$
The variable $x^\prime=x$, but $y^{\prime}(x)$ is quite complicated and will not be given here (see \cite{Okamoto}).

\vskip 0.3 cm
For the Taylor solutions, we have  $
\sigma= \theta_1-\theta_\infty
$. Denote $ \sigma^\prime:=\theta_1^\prime-\theta_\infty^\prime$. 
If we start from a Taylor solution constructed in sections  \ref{irregular0},
\ref{semireducible}, \ref{supersemireducible} by means of  the
matching procedure, developed for   $\sigma=\theta_1-\theta_\infty$,
then the birational transformations allow to  obtain the solutions  
defined for $\sigma= \pm(\theta_1\pm\backslash \mp\theta_\infty)+n$,
$n\in \bf Z$. This is a consequence of the following actions: 
$$
l_1\hbox{ and }l_4:~\sigma\mapsto \sigma^\prime=\sigma+1.~~~~~
l_2\hbox{ and }l_3:~\sigma\mapsto \sigma^\prime=\sigma-1 
$$
$$
w_3  \cdot l_1 \cdot l_1:~\sigma\mapsto \sigma^\prime=
\theta_1+\theta_\infty.~~~~~~~ 
w_1 \cdot w_3 \cdot l_1 \cdot l_1:~\sigma\mapsto \sigma^\prime= -\sigma.
$$
Note also that (PVI) is invariant for $\theta_1\mapsto
-\theta_1$. This maps $\sigma\mapsto -\sigma$. Other actions are: 
$$
w_1:~~~\sigma^\prime= -(\theta_1+\theta_\infty);~~~~~~~~w_2:~~~\sigma^\prime= \sigma;~~~~~~~~w_3:~~~\sigma^\prime= \theta_1+\theta_\infty-2;~~~~~~~~w_4:~~~\sigma^\prime= \theta_1+\theta_\infty-2.
$$
$$
x^1:~~~\sigma^\prime= \theta_0-\theta_\infty;~~~~~~~~~~x^2:~~~\sigma^\prime= \theta_1-\theta_0;~~~~~~~~~~x^3:~~~\sigma^\prime= \theta_1-\theta_\infty=\sigma.
$$

\vskip 0.2 cm

\section{Appendix 4:  Examples of Taylor Expansions}
\label{EXAMPLES}
 We give the solutions of (PVI)  of the form
 $y(x)=\sum_{n=0}^\infty b_n~x^n$, $b_0\neq 0$, depending on the
 value of the coefficients
 $\theta_0,\theta_1,\theta_x,\theta_\infty$. Similar examples can be
 constructed for solutions  $y(x)=\sum_{n=1}^\infty b_n~x^n$,
 $b_1\neq 0$, by the symmetry (\ref{nuovasimmetria}). 

We observe that, in general, the free parameter appears starting from 
some power $x^m$. If we truncate the series at $x^{m-1}$, we cannot see it.
   
We always denote $b_n$ the coefficients, though they change case by
case.

\vskip 0.3 cm 
{\bf * Example 1:} (PVI) always has one solution  (\ref{form1}) when 
 $\theta_\infty-\theta_1\not\in{\bf Z}$, and one solution  (\ref{riuffa})  when $\theta_\infty+\theta_1\not\in{\bf Z}$.

\vskip 0.3 cm 
{\bf * Example 2:} Case    $\theta_1+\theta_\infty=0$.

\vskip 0.2 cm
\noindent
i) There are  solutions (\ref{form2}), if $\theta_x=\pm\theta_0$, $\theta_\infty\neq 1$.

\vskip 0.2 cm 
\noindent
ii) There is a solution  (\ref{form1}),  
%\be
%y(x)= {2\theta_\infty-1\over \theta_\infty-1}
%-{4\theta_\infty(\theta_\infty-1)+\theta_x^2-\theta_0^2\over
%  8(\theta_\infty-1)}~x~+\sum_{n=3}^\infty
%b_n(-\theta_\infty,\theta_\infty,\theta_0,\theta_x)~x^n,
%\label{uffauffa1}
%\ee
defined for $\theta_\infty\neq1,~{2n+1\over 2}$, where $n\in{\bf Z}$. 
The condition $\theta_\infty\neq {2n+1\over 2}$ follows from the condition 
$\theta_1-\theta_\infty\neq n$ in (\ref{form1}), when
$\theta_1=-\theta_\infty$. The denominators of the
coefficients  $ b_n$  vanish  for half-integer $\theta_\infty$, and
for $\theta_\infty=1$.

\vskip 0.2 cm
\noindent
 iii) For $\theta_\infty={2n+1\over 2}$, solutions $y(x)=\sum_{n=0}^\infty
b_nx^n$ in ii) are not defined. On the other hand, solutions  (\ref{form2}) are defined, provided that  $\theta_x=\pm\theta_0$. Moreover, for a given $\theta_\infty={2n+1\over 2}$, there may be solutions equivalent to (\ref{form2}), provided that  $\theta_x\pm\theta_0$ assumes some integer values. For example, consider $\theta_\infty={3\over 2}$. We have (\ref{form2}):
$$
y(x)=-2~+a~x~+\left(\theta_0^2-1+{3\over 2}~a-{1\over 2}~a^3\right)~x^2~+~
\sum_{n=3}^\infty b_n(a;\theta_0,{3\over 2})~x^n,~~~~
\theta_x=\pm \theta_0;
$$
and a solution equivalent to (\ref{form2}): 
$$
y(x)=4 -(2\pm \theta_0)~x~+ a~x^2~+\sum_{n=3}^\infty b_n(a;\theta_0)~x^n,
           ~~~~~
\theta_x=\pm(1+\theta_0),~~\pm(1-\theta_0). 
$$ 

\vskip 0.2 cm
\noindent
iv) If $\theta_\infty=0,1$, we do not have any solution of the desired
form, except for the singular solutions $y=0,1$. If $\theta_\infty=2$ we have a solution equivalent to    (\ref{form3}). The parameter $a$ is the coefficient
 of $x^3$  (coefficients of $x^0$, $x$ and $x^2$ have no parameter).

\vskip 0.3 cm
{\bf * Example 3:} Case 
$\theta_\infty-\theta_1=0$. 

\vskip 0.2 cm
\noindent
i) We have solutions (\ref{form2})

\vskip 0.2 cm
\noindent
ii)
 We  have the  solution  (\ref{riuffa}), with the  substitution: $\theta_1=\theta_\infty$. This is  defined for $\theta_\infty \neq 1,{2n+1\over 2}$. 

\vskip 0.2 cm
\noindent
iii) For  $\theta_\infty={2n+1\over 2}$ we find solutions of exactly the form  (\ref{form2}), when $\theta_x=\pm\theta_0$. Moreover, for any given
 $\theta_\infty={2n+1\over 2}$, there may be solutions equivalent to
 (\ref{form2}), provided that  $\theta_x\pm\theta_0$ has some integer
 value. For example, consider $\theta_\infty=-{1\over 2}$. There are
 solutions  (\ref{form2}): 
$$
y(x)={2\over 3}~+a~x~+\sum_{n=2}^\infty b_n(a;\theta_0,-{1\over
  2})~x^n,~~~~~\theta_x=\pm\theta_0; 
$$
and the equivalent solutions:
$$
y(x)={4\over 3}-{2\pm \theta_0\over 9}~x~+a~x^2~+ \sum_{n=3}^\infty
b_n(a;\theta_0)x^n, ~~~~~\theta_x=\pm(1+\theta_0)~\hbox{ or } \pm(1-\theta_0).  
$$

Another example: consider  
 $\theta_\infty={1\over 2}$. In this case we just have  (\ref{form2}),
 or the singular solution $y=1$.  

\vskip 0.2 cm
\noindent
iv) For $\theta_\infty=0,1$ we don't have solutions of the desired
form, except for the singular solutions $y=0,1$.   
If $\theta_\infty=2$ we have a solution equivalent to  
  (\ref{form3}). The parameter $a$ is the coefficient
 of $x^3$  (coefficients of $x^0$, $x$ and $x^2$ have no parameter).

\vskip 0.3 cm 
{\bf * Example 4:} Case $\theta_1-\theta_\infty=-1$. 
\vskip 0.2 cm
\noindent

\vskip 0.2 cm
\noindent
i) We have the solution (\ref{riuffa}),
% which becomes: 
%$$
%   y(x)~=2~-{3+4\theta_\infty^2-8\theta_\infty-\theta_0^2+\theta_x^2\over 
%2(2\theta_\infty-1)(2\theta_\infty-3)}~x~+\sum_{n=3}^\infty
%   b_n(1-\theta_\infty,\theta_\infty,\theta_0,\theta_x)~x^n,
%$$
 defined for $\theta_\infty\neq 1,~{2n+1\over 2}$, where $n\in{\bf Z}$.

\vskip 0.2 cm
\noindent
ii) For any $\theta_\infty={2n+1\over 2}$,  solutions (\ref{form2})
are defined, when $\theta_0\pm\theta_x=0$. Moreover, solutions may
exist  equivalent to (\ref{form2}) by symmetry, provided that
$\theta_0\pm\theta_x$ is some integer. 

\vskip 0.2 cm
\noindent
iii) For $\theta_\infty=1$ we are exactly in the case (\ref{form3}).

\vskip 0.2 cm
\noindent
iv) For $\theta_\infty=0,2$ we have solutions equivalent to (\ref{form3}):
$$
y(x)=2~+{\theta_0^2-\theta_x^2-3\over
  6}~x~+a~x^2~+\left\{{[(\theta_0-\theta_x)^2-9] 
[(\theta_0+\theta_x)^2-9](\theta_0^2-\theta_x^2)\over 4320}+
~a~{\theta_0^2-\theta_x^2+3\over 6}\right\}~x^3~+...   
$$

\vskip 0.3 cm 
{\bf * Example 5:} Case $\theta_1-\theta_\infty=1$. 
\vskip 0.2 cm 
\noindent
i) We have solution (\ref{riuffa}), with the substitution of
$\theta_1=\theta_\infty+1$, 
%$$
%y(x)= {2\theta_\infty\over
%  \theta_\infty-1}~-{(\theta_\infty+1)(4\theta_\infty^2-1+\theta_x^2-
%\theta_0^2)\over 2(\theta_\infty-1)(2\theta_\infty
%-1)(2\theta_\infty+1)}~x
%~  +\sum_{n=3}^\infty 
% b_n(-\theta_\infty-1,\theta_\infty,\theta_0,\theta_x)~x^n$$
 and defined for $\theta_\infty\neq 1,~{2n+1\over 2}$, $n\in{\bf Z}$. 

\vskip 0.2 cm
\noindent
ii) Equivalent to (\ref{form2}) by symmetry,  we have the solution: 
$$
y(x)= {2\over 1-\theta_\infty}~+{(\theta_\infty+1)(2\pm \theta_0)\over 
3(\theta_\infty -1)}~x~+ax^2+ \sum_{n=3}^\infty
b_n(a,\theta_\infty,\theta_0)x^n,$$ 
$$
   \theta_\infty\neq 0,1;~~~~\theta_x= \pm(1-\theta_0)\hbox{ or } 
\pm(1+\theta_0)
$$
The two signs in the coefficient of $x$ depend on the choice
$\theta_x= \pm(1-\theta_0)$ or $  
\pm(1+\theta_0)$ respectively.  Similar change of signs occur in all
the coefficients $b_n$. 

\vskip 0.2 cm 
\noindent
iii) If $\theta_\infty=0$,  the  solution in case ii) 
is not defined (denominators in the coefficients $b_n$  vanish). We
have anyway  a solution equivalent 
 to (\ref{form3}): 
$$
y(x)= 2~+{\theta_0^2-\theta_x^2-3\over 6}~x~+ a~x^2~+
\sum_{n=3}^\infty
 b_n(a;\theta_0,\theta_\infty) x^n.
$$
If $\theta_\infty=2$, we have a solution equivalent to (\ref{form3}). The 
parameter $a$ is the coefficient of $x^4$ 
(no parameter in lower powers of $x$). 

\vskip 0.2 cm 

\noindent
iv) If $\theta_\infty={2n+1\over 2} $ the solution in i) is not
defined. Solutions exist equivalent to (\ref{form2}), provided that
$\theta_x\pm\theta_0$ is some integer.

\vskip 0.2 cm
\noindent
v) For $\theta_\infty=-1$ we have a solution equivalent to (\ref{form3}):
$$
y(x)=1~+a~x^2~+~{a(\theta_x^2-\theta_0^2+3)\over 6}~x^3~
+\sum_{n=4}^\infty b_{n}(a;\theta_0,\theta_x)~x^n.
$$
\vskip 0.2 cm
\noindent
vi) For $\theta_\infty=1$, solutions of the desired form do not exist,
except for $y=1$.  
\vskip 0.3 cm

\vskip0.3 cm

We could proceed at our pleasure, choosing
 any value $\theta_\infty\pm\theta_1$ integer. We would always have solutions 
 of three kinds. 1) {\it One} out of the two solutions  (\ref{form1}) 
 and  ({\ref{riuffa}). 2) Solutions equivalent to (\ref{form2}) --
 at least when $\theta_\infty $ is 
half integer -- provide that $\theta_x\pm\theta_0$ is some integer. 3)  
Solutions equivalent to (\ref{form3}), for $\theta_\infty$ equal to some integer. 

%%%%%%%%%%%%
%%%%%%       END OF \small
%%%%%%%%%%%%
}

\end{document}